\documentclass[11pt,twoside]{article}
\usepackage[hmargin=1.25in,vmargin=1.25in, a4paper, centering]{geometry}
\usepackage[nottoc]{tocbibind}

\usepackage[biblabel,nosort,nocompress]{cite}

\usepackage{graphicx}
\usepackage{amsmath,amsthm}
\usepackage{amssymb}
\usepackage{rotate}

\usepackage{upgreek}
\usepackage{biolinum}
\usepackage{garamondlibre}
\usepackage[notext,noDcommand,partialup]{kpfonts}
\usepackage{courierten}

\usepackage{physics}

\usepackage{graphicx}
\usepackage{mathtools,slashed}
\usepackage{cases}
\usepackage{xfrac}
\usepackage{cancel}

\usepackage{enumitem}

\usepackage[dvipsnames]{xcolor}
\usepackage[colorlinks=true, pdfstartview=FitV,linkcolor=RoyalBlue,citecolor=ForestGreen, urlcolor=BurntOrange]{hyperref}

\usepackage{tikz}
\usetikzlibrary{hobby}

\usepackage{fancyhdr}
\renewcommand\footnotemark{}
\usepackage[sf,sl,outermarks]{titlesec}
\titleformat{\section}
{\Large\bfseries\sffamily}{\filcenter\Large{\thesection.}}{1em}{\filcenter}
\titleformat{\subsection}
{\large\bfseries\sffamily}{\thesubsection.}{1em}{}
\titleformat{\subsubsection}[runin]
{\normalsize\bfseries\itshape\sffamily}{\thesubsubsection.}{0.5em}{}[.\hspace*{0.5ex}]
\titleformat{\paragraph}[runin]
{\normalsize\itshape\sffamily}{\theparagraph.}{0.5em}{}[.\hspace*{0.5ex}]

\usepackage{abstract}

\theoremstyle{definition}
\newtheorem{defi}{\sffamily Definition}[section]

\theoremstyle{plain}
\newtheorem{theorem}[defi]{\sffamily Theorem}
\newtheorem{prop}[defi]{\sffamily Proposition}
\newtheorem{lemma}[defi]{\sffamily Lemma}
\newtheorem{cor}[defi]{\sffamily Corollary}

\theoremstyle{remark}
\newtheorem*{remark}{\sffamily Remark}

\makeatletter
\providecommand{\proofnamestyle}{\itshape\sffamily\bfseries}

\makeatother
\newcommand*{\wt}[1]{\widetilde{#1}}
\newcommand*{\wh}[1]{\widehat{#1}}
\newcommand*{\R}{\mathbb{R}}

\providecommand\given{} 
\newcommand\SetSymbol[1][]{\nonscript\:#1\vert \allowbreak \nonscript\: \mathopen{}}
\DeclarePairedDelimiterX\Set[1]\{\}{ \renewcommand\given{\SetSymbol[\delimsize]} #1 }

\newcommand*{\lap}{\mathop{}\!\operatorname{\vb*{\Delta}}}

\DeclareMathOperator*{\Curl}{\operatorname{curl}}
\newcommand*{\const}{\operatorname{const.}}
\newcommand*{\Dt}{\mathop{}\!\textbf{D}_t}

\newcommand*{\pd}{\partial}

\newcommand*{\vv}{{\vb{v}}}
\newcommand*{\vw}{{\vb{w}}}
\newcommand*{\vh}{{\vb{h}}}
\newcommand*{\Om}{{\Omega}}
\newcommand*{\Omt}{\Omega_t}
\newcommand*{\Gmt}{\Gamma_t}
\newcommand*{\vtau}{{\vb*{\tau}}}

\newcommand*{\cE}{\mathcal{E}}
\newcommand*{\fks}{\mathfrak{s}}

\newcommand*{\cD}{\mathcal{D}}

\newcommand*{\cH}{{\mathcal{H}}}
\newcommand*{\vxi}{\vb*{\xi}}
\newcommand*{\vn}{{\vb{n}}}
\newcommand*{\vbu}{{\vb{u}}}
\newcommand*{\fkR}{\mathfrak{R}}

\newcommand*{\cN}{\mathop{}\!{\mathcal{N}}}
\newcommand*{\vPhi}{\vb*{\Phi}}
\newcommand*{\Oms}{\Omega_*}
\newcommand*{\Gms}{\Gamma_*}

\newcommand*{\scL}{\mathscr{L}}
\newcommand*{\cX}{\mathcal{X}}
\newcommand*{\cDt}{\mathop{}\!{\cD_t}}
\newcommand*{\cP}[1]{\mathcal{P}\qty(#1)}
\newcommand*{\cQ}{\mathcal{Q}}

\newcommand*{\Elin}{E_{\textup{lin}}}
\newcommand*{\Mlin}{M_{\textup{lin}}}
\newcommand*{\err}{\textup{err. }}
\newcommand*{\cL}{\mathcal{L}}

\newcommand*{\ka}{{\varkappa_a}}
\newcommand*{\fkk}{\mathfrak{K}}
\newcommand*{\udl}{\underline}
\newcommand*{\Dbt}{\mathop{}\!\mathcal{D}_{\beta}}

\newcommand*{\scA}{\mathop{}\!\mathscr{A}}
\newcommand*{\scB}{\mathop{}\!\mathscr{B}}
\newcommand*{\scF}{\mathop{}\!\mathscr{F}}
\newcommand*{\scG}{\mathscr{G}}
\newcommand*{\fkX}{\mathfrak{X}}
\newcommand*{\fkI}{\mathfrak{I}}
\newcommand*{\fkV}{\mathfrak{V}}
\newcommand*{\fkH}{\mathfrak{H}}

\newcommand*{\er}{{\vu{e}_r}}
\newcommand*{\etheta}{{\vu{e}_\theta}}
\newcommand*{\vX}{\vb{X}}
\newcommand*{\scZ}{\mathscr{Z}}
\newcommand*{\vY}{\vb{Y}}
\newcommand*{\vps}{{\vb*{\uppsi}}}

\newcommand*{\kaz}{{\varkappa_a^0}}
\newcommand*{\omz}{{\omega_*^0}}
\newcommand*{\jz}{{\jmath_*^0}}
\newcommand*{\xiz}{{\vb*{\xi}_*^0}}
\newcommand*{\vz}{{\vb{v}^0}}
\newcommand*{\hz}{\vb{h}^0}
\newcommand*{\Hz}{\vb{H}^0}
\newcommand*{\kao}{{\varkappa_a^1}}
\newcommand*{\omo}{{\omega_*^1}}
\newcommand*{\jo}{{\jmath_*^1}}
\newcommand*{\cB}{\mathcal{B}}
\newcommand*{\fkT}{\mathfrak{T}}
\newcommand*{\vbz}{{\vb{b}^0}}
\newcommand*{\kal}{{\varkappa_a^l}}
\newcommand*{\vH}{{\vb{H}}}
\newcommand*{\cU}{\mathcal{U}}
\newcommand*{\cS}{\mathcal{S}}
\newcommand*{\vN}{{\vb{N}}}
\newcommand*{\smallwedge}{\operatorname{\scalebox{0.8}{$\vb*{\wedge}$}}}
\newcommand*{\cV}{\mathcal{V}}
\newcommand*{\cJ}{\mathcal{J}}

\newcommand*{\vnu}{{\vb*{\upnu}}}
\newcommand*{\vphi}{{\vb*{\upphi}}}

\newcommand*{\fkE}{\mathfrak{E}}

\numberwithin{equation}{section}

\newcommand*{\dive}{\divergence}

\begin{document}
	\title{%
		\bfseries \biolinum
		On the 2D Plasma-Vacuum Interface Problems for Ideal Incompressible MHD%
		\footnote{{\bf\biolinum Date: }{\today}.}
		\footnote{{\bf\biolinum MSC(2020): }Primary 35Q35, 76W05, 76B03, 76E25.}%
		\footnote{{\bf\biolinum Keywords: }ideal MHD, free boundary problems, local well-posedness, surface tension, non-degenerate magnetic fields, the Rayleigh-Taylor sign condition.}%
	}
	\author{%
		\scshape Sicheng LIU
		\thanks{%
			{\biolinum \textbf{Sicheng LIU} (\textit{Corresponding author}): } Department of Mathematics, Faculty of Science and Technology, University of Macau, Taipa, Macao SAR, China.
			\textit{E-mail}: \href{mailto:scliu@link.cuhk.edu.hk}{\nolinkurl{scliu@link.cuhk.edu.hk}}. 
		}
		\and
		\scshape Tao LUO
		\thanks{%
			\textbf{\biolinum Tao LUO: } Department of Mathematics, City University of Hong Kong, Kowloon, Hong Kong SAR, China.
			\textit{E-mail}: \href{mailto:taoluo@cityu.edu.hk}{\nolinkurl{taoluo@cityu.edu.hk}}.
		}
	}
	\date{} 
	\maketitle
	\begin{abstract}  
		This manuscript concerns the stability conditions for the well-posedness of the two-dimensional plasma-vacuum interface problems for ideal incompressible magnetohydrodynamics (MHD) equations, which describe the dynamics of conducting perfect fluids in a vacuum region under the influence of magnetic fields. Due to the counterexamples constructed by {\sc C.~Hao} and the second author [{\it Comm. Math. Phys.} {\bf 376} (2020), 259-286], the plasma-vacuum problems violating the Rayleigh-Taylor sign condition and without surface tension are believed to be highly unstable or even ill-posed. However, based on a geometric perspective, we demonstrate that, although the initial data provided therein are ill-posed in the Lagrangian framework characterizing the regularity of flow maps, plasma-vacuum problems with such initial data can still be stable/well-posed in the Eulerian setting without involving flow maps, under the hypothesis that the total magnetic fields are non-degenerate on the free boundary. This discloses the essential difference between the Eulerian and Lagrangian frameworks of MHD free boundary problems and the stabilizing effect of the magnetic field. We also consider the surface tension effect and establish the local well-posedness theories in standard Sobolev spaces in this case. Furthermore, under either the non-degeneracy assumption on magnetic fields or the Rayleigh-Taylor sign condition on the effective pressure, we prove vanishing surface tension limits. These results indicate that both capillary forces and non-degenerate tangential magnetic fields can stabilize the motion of the plasma-vacuum interface, and the stabilizing effect of surface tension is stronger than that of the tangential magnetic fields.
	\end{abstract}
	
	{\sffamily\tableofcontents}
	
	\addtocontents{toc}{\protect\setcounter{tocdepth}{2}}
	\section{Introduction}
	In this manuscript, we shall study the vacuum free boundary problems for ideal incompressible magnetohydrodynamics (MHD) equations. The problem can be written as
	\begin{subequations}\label{PV problem}
		\begin{equation}\label{MHD}
			(\text{MHD})\, \begin{cases*}
				\pd_t \vv + \qty(\vv \vdot \grad) \vv + \grad p = \qty(\vh \vdot \grad) \vh &in $ \Omt$, \\
				\pd_t \vh + \qty(\vv \vdot \grad) \vh = \qty(\vh \vdot \grad) \vv &in $ \Omt $, \\ 
				\dive \vv = \dive \vh = 0 &in $ \Omt $;
			\end{cases*}
		\end{equation}
		\begin{equation}\label{pre maxwell}
			(\text{pre-Maxwell})\,\begin{cases*}
				\dive\vH = 0 &in $\cV_t$, \\
				\Curl \vH = 0 &in $\cV_t$;
			\end{cases*}
		\end{equation}
		with boundary conditions:
		\begin{equation}\label{BC}
			(\text{BC}) \begin{cases*}
				\vv \vdot \vn = \fks &on $ \Gmt $, \\
				p = \alpha \varkappa + \frac{1}{2}\abs{\vH}^2 &on $ \Gmt $, \\
				\vh \vdot \vn = \vH \vdot \vn = 0  &on $ \Gmt $, \\
				\vN \smallwedge \vH = \cJ &on $\cS$,
			\end{cases*}
		\end{equation}
		where $\vN\smallwedge\vH = N^1H^2-N^2H^1$ in Cartesian's coordinates.
	\end{subequations}
	
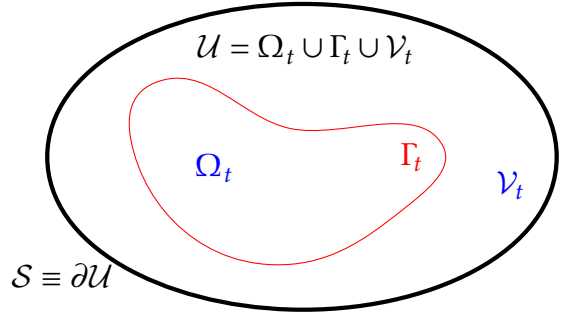
\begin{center}
		\begin{tikzpicture}[use Hobby shortcut, yscale=0.75]
			\path
			(3.5,0) coordinate (a0)
			(0,2.7) coordinate (a1)
			(-3.2,0) coordinate (a2)
			(0,-2.7) coordinate (a3);
			\draw[closed,ultra thick] (a0) .. (a1) .. (a2) .. (a3);
			\node[blue] at (2.9,-0.5) {\large $ \cV_t $};
			\node at (-3, -2.1) {\large $\cS \equiv \pd\cU$};
			\node at (0.2, 1.9) {\large ${\cU} = \Omt \cup \Gmt \cup\cV_t$};
			
			\path
			(0,0.5) coordinate (b0)
			(2,0.2) coordinate (b1)
			(1.6,-0.8) coordinate (b2)
			(0,-1.9) coordinate (b3)
			(-2.1,0.3) coordinate (b4)
			(-1.8,1.3) coordinate (b5);
			\draw[red,closed] (b0) .. (b1) .. (b2) .. (b3) .. (b4) .. (b5);
			\node[red] at (1.6, 0) {\large$ \Gamma_t $};
			\node[blue] at (-1,-0.2) {\large$ \Omega_t $};
		\end{tikzpicture}
	\end{center}
	
	Here $ \Omt \subset \R^2 $ is the plasma region evolving with the conducting fluid, $\Gmt \coloneqq \pd\Omt$ is the free interface, $\cV_t$ is the vacuum region, $\cU = \Omt\cup\Gmt\cup\cV_t$ is a fixed bounded domain, and $\cS \coloneqq \pd\cU$ is the fixed boundary, which represents a perfect conducting wall. Vector fields $\vv$ and $\vh$ are respectively the velocity and magnetic fields of the plasma, $\vH$ represents the magnetic field in the vacuum region, and $p$ is the effective pressure of the plasma, which is the sum of the fluid one and the magnetic one. Furthermore, $ \vn $ is the outer normal of $ \pd\Omt $, $ \fks $ is the normal speed of $ \Gmt $ in the direction of $ \vn $, $\varkappa$ is the signed curvature of $\Gmt$ with respect to $\vn$, $\alpha \ge 0$ is the constant surface tension coefficient, $\vN$ is the outer unit normal of $\pd\cU$, and $\cJ$ represents the prescribed surface current on the conducting wall $\cS$.
	
	The system \eqref{MHD} is derived from the Euler and Maxwell equations ignoring the effect of displacement currents, which is legitimate for non-relativistic motions. \eqref{pre maxwell} is called the pre-Maxwell system, which is the simplification of the full Maxwell one under the hypothesis that the amplitude of plasma velocity is negligible when compared to the speed of light. The first boundary condition in \eqref{BC} states that the free boundary evolves with the fluid particles lying on it. The second one can be deduced from the balance of momenta. The third condition follows from the fact that magnetic fields are always divergence free. The forth one is derived from Ampere's law. For the detailed discussions, one can refer to \cite[Ch. 8]{Landau-Lifshitz-Vol8} and \cite[Ch. 4]{goedbloed_poedts_2004}.
	
	\subsection{Backgrounds and Related Works}
	
	MHD combines principles of fluid dynamics and electromagnetism, providing a comprehensive framework to analyze the interactions between magnetic fields and conductive fluids. The mathematical study of MHD involves solving complicated differential equations to predict fluid behavior under different conditions, which is essential for accurate simulations and practical applications in physics, medical treatments, and industrial processes.	
	Vacuum free boundary problems focus on the dynamics of the interface separating a plasma and a vacuum, which have plenty of applications, such as engineering (e.g. the design of nuclear reactors, heat exchangers, and electronic devices) and astrophysics (e.g. understanding solar flares). The mathematical study of free boundary problems for MHD is indispensable for advancing technology, improving medical treatments, and enhancing our understanding of natural phenomena.
	
	Free boundary problems in hydrodynamics, such as well-posedness for water waves and vortex sheets, have made great advances in recent decades. Concerning the local well-posedness of free-boundary incompressible Euler equations in Sobolev spaces, one can refer to \cite{Wu1997,Wu1999,Alazard-Burq-Zuily2014,Hunter-Ifrim-Tataru2016} for irrotational water waves, and \cite{Christodoulou-Lindblad2000,Lindblad2005,Coutand-Shkoller2007,zhangzhang2008,Shatah-Zeng2008-Geo,Shatah-Zeng2011,Wang_etal2021,ifrim2023sharp} for one-phase (fluid-vacuum) flows with vorticity. (See also \cite{Ebin1987} for an ill-posed example with initial data violating the Rayleigh-Taylor sign condition.) The well-posedness results for vortex sheet problems with surface tension were established in \cite{Cheng-Coutand-Shkoller2008,Shatah-Zeng2008-vortex,Shatah-Zeng2011}. (See also \cite{Ebin1988} for an ill-posed example without capillary forces.)
	
	The mathematical analysis for the free-boundary MHD problems are more involved than that for the Euler equations, as one needs to handle two deeply interacted hyperbolic systems simultaneously.	Due to rich physical phenomena, the free interface problems for MHD can be categorized into several scenarios (we focus on the discussions of results on incompressible ideal MHD closely related to this paper here):
	\begin{itemize}
		\item \textbf{Current-vortex sheets}
		
		Such problems describe the motions of two plasmas separated by a free interface, along which both magnetic and velocity fields admit tangential jumps. {\sc Syrovatskij} \cite{Syrovatskij} and {\sc Axford} \cite{axford1962note} raised the stability conditions for planer incompressible current-vortex sheets (cf. \cite[eq. (71.14)]{Landau-Lifshitz-Vol8}; see also \cite{Li-Li2021,Liu-Xin2023} for general versions):
		\begin{equation}\label{Syr 1"}
			2\abs{\vh_+}^2 + 2\abs{\vh_-}^2 > \abs{\llbracket\vv\rrbracket}^2,
		\end{equation}
		\begin{equation}\label{Syr 2"}
			2\abs{\vh_+ \cp \vh_-}^2 \ge  \abs{\vh_+ \cp \llbracket\vv\rrbracket}^2 +  \abs{\vh_-\cp\llbracket\vv\rrbracket}^2,
		\end{equation}
		where $ \llbracket\vv\rrbracket \coloneqq \vv_+ - \vv_- $ is the velocity jump. If \eqref{Syr 2"} were replaced by
		\begin{equation}\label{Syro 3"}
			2\abs{\vh_+ \cp \vh_-}^2 > \abs{\vh_+ \cp \llbracket\vv\rrbracket}^2 + \abs{\vh_-\cp\llbracket\vv\rrbracket}^2,
		\end{equation}
		from which \eqref{Syr 1"} follows, one can refer to {\sc Morando, Trakhinin, and Trebeschi} \cite{Morando-Trakhinin-Trebeschi2008}, {\sc Sun, Wang, and Zhang} \cite{Sun-Wang-Zhang2018}, and the first author and {\sc Z. Xin} \cite{Liu-Xin2023} for the local well-posedness under \eqref{Syro 3"}. (See also \cite{Trakinin2005} and \cite{Coulombel-Morando-Secchi-Trebischi2012} for energy estimates with a stronger stability condition.) For the results under the effect of capillary forces, one can consult {\sc Li and Li} \cite{Li-Li2022}. (See also \cite{Liu-Xin2023} for results without graph assumptions on the interface.) A very recent work by {\sc Cai and Lei} \cite{cai2024global} demonstrated the 2D global existence under the Syrovatskij-type condition, which can be viewed as small perturbations to non-zero constant magnetic fields.
		
		\item \textbf{Plasma-vacuum problems}
		
		They characterize the dynamics of plasma bulks surrounded by vacuum states, which can be further separated into several cases (here we denote by $\alpha$ and $\vH$ respectively the surface tension coefficient and the magnetic field in the vacuum region):
		\begin{itemize}
			
			\item {\it $\alpha=0$, $\vH \equiv \vb{0}$, and under the Rayleigh-Taylor sign condition:}
			\begin{equation}\label{RT}
				-\grad_\vn p > 0 \qq{on} \Gmt.
			\end{equation}
			
			The first advancement on nonlinear problems might be the a priori estimate given by {\sc C. Hao} and the second author in \cite{Hao-Luo2014}. Local well-posedness results were established by {\sc Gu and Wang} \cite{Gu-Wang2019} in the 3D model domain $\Omega = \mathbb{T}^2 \times (0, 1)$ with flat initial interface.
			Without the graph assumption on the free interface, the well-posedness of linearized problems were shown by {\sc Hao} and the second author in \cite{Hao-Luo2021}. The 3D local well-posedness with graph assumptions was also proven by {\sc Zhao} \cite{zhao2024local}. The Hadamard-style local well-posedness and continuation criteria for general fluid domains were recently demonstrated by {\sc Ifrim, Pineau, Tataru, and Taylor} \cite{ifrim2024sharp} in low-regularity Sobolev spaces.  These well-posedness results are based on the Rayleigh-Taylor sign condition \eqref{RT}.  Here we remark that, if \eqref{RT} is strictly violated on the free boundary, the examples given by {\sc Hao} and the second author in \cite{Hao-Luo2020} indicate the ill-posedness for 2D problems in  Sobolev spaces of the Lagrangian variable.
			
			\item {\it 3D problems, $\alpha=0$, and $\abs{\vh\cp\vH}>0$ on $\Gmt$ (non-collinearity condition):}
			
			One can refer to {\sc Morando, Trakhinin, and Trebeschi} \cite{Morando-Trakhinin-Trebeschi2014} for the liner well-posedness and {\sc Sun, Wang, and Zhang} \cite{Sun-Wang-Zhang2019} for the nonlinear local well-posedness, both under graph assumptions on the interface. Nonlinear results without graph assumptions were shown by the first author and {\sc Z. Xin} in \cite{liu-Xin2025free}.
			
			\item {\it 3D problems, $\alpha > 0$:}
			
			The local well-posedness in the model domain $\Omega=\mathbb{T}^2 \times (0, 1)$ with flat initial surface and vanishing $\vH$ was established by {\sc Gu, Luo, and Zhang} \cite{Gu-Luo-Zhang2021} (see also their vanishing surface-tension limit results \cite{Gu-Luo-Zhang22} under \eqref{RT}). Without the graph assumptions on free interface, the local well-posedness and vanishing surface-tension limit were shown in \cite{liu-Xin2025free} for possibly non-zero $\vH$. See also \cite{hao2023motion} for the a priori estimates and blow-up criteria with zero $\vH$.
		\end{itemize}
	\end{itemize}
	Additionally, when the effects of magnetic diffusion (which induce certain dissipation), surface tension, and transversality of the magnetic field to the interface are taken into account, global solutions are obtained in  \cite{Wang-Xin2021}.
	
	Given the counter examples presented by {\sc Hao} and the second author,  the 2D plasma-vacuum problem in the case when $\alpha=0$ and  $\vH \equiv \vb{0}$ is proved to be  ill-posed in some Sobolev spaces of Lagrangian variable (cf. \cite{Hao-Luo2020}). It is believed that, without capillary forces nor the Rayleigh-Taylor sign condition, the 2D plasma-vacuum problems are highly unstable or even ill-posed  in terms of the Largrangian variables of the flow map. However, based on geometrical observations on the interface evolutions, we find out that  the problems can be  well-posed in suitably chosen standard Sobolev spaces  in Eulerian coordinates without involving Lagrangian variables even when the  the Rayleigh-Taylor sign condition is violated when the tangential magnetic fields are non-degenerate on the free interfaces. This discloses somewhat surprising difference of the well-posedness in the Eulerian and Lagrangian frameworks for MHD free boundary problems.  It worths to point out,   unlike the three-dimensional cases, the well-posedness theories for two-dimensional plasma-vacuum interface problems (even under the graph assumption on the interface) without assuming the Rayleigh-Taylor sign condition \eqref{RT} are still open. Identifying suitable stability conditions to  establish such a  well-posedness theory is another motivation of this work. The results in this paper indicate that, besides the capillary forces, non-degenerate tangential magnetic fields can stabilize the motion of the plasma-vacuum interface, distinguishing MHD free interface problems from those for hydrodynamics. 
	
	\subsection{Conservation of the Physical Energy}
	
	In order to further indicate the physical background of the problem \eqref{PV problem}, we now consider the $L^2$-type energy defined by
	\begin{equation}\label{def Ephy}
		E_{\text{phy}} \coloneqq \frac{1}{2}\int_{\Omt} \abs{\vv}^2 + \abs{\vh}^2 \dd{x} + \frac{1}{2}\int_{\cV_t} \abs{\vH}^2 \dd{x} + \int_{\Gmt} \alpha \dd{\ell},
	\end{equation}
	which consists in the kinetic, magnetic, and surface energies. Then, it follows from \eqref{MHD} and the transport formula in fluid dynamics (see, for example, \cite[p. 10]{chorin-marsden1993}) that
	\begin{equation*}
		\begin{split}
			\dv{t} \frac{1}{2}\int_{\Omt} \abs{\vv}^2 + \abs{\vh}^2 \dd{x} &= \int_{\Omt} \vv\vdot(\pd_t + \vv\vdot\grad)\vv + \vh\vdot(\pd_t+\vv\vdot\grad)\vh \dd{x} \\
			&= \int_{\Omt} -\grad p \vdot \vv + (\vh\vdot\grad)(\vv\vdot\vh) \dd{x} \\
			&= \int_{\Gmt} - p (\vv\vdot\vn) + (\vh\vdot\vn)(\vv\vdot\vh) \dd{\ell}.
		\end{split}
	\end{equation*}
	The variational formula for the curve length implies that
	\begin{equation*}
		\dv{t} \int_{\Gmt} \alpha \dd{\ell} = \int_{\Gmt} \alpha\fks\varkappa \dd{\ell}.
	\end{equation*}
	Furthermore, since $\cV_t$ is a moving domain, the transport formula yields that
	\begin{equation*}
		\dv{t} \frac{1}{2}\int_{\cV_t} \abs{\vH}^2 \dd{x} = \int_{\cV_t} \vH\vdot\pd_t\vH \dd{x} - \frac{1}{2} \int_{\Gmt} \fks \abs{\vH}^2 \dd{\ell}.
	\end{equation*}
	Therefore, the boundary conditions \eqref{BC} yield that
	\begin{equation}\label{dt Ephy pre}
		\dv{t} E_{\text{phy}} = \int_{\cV_t} \vH\vdot\pd_t\vH \dd{x} - \int_{\Gmt} \fks \abs{\vH}^2 \dd{\ell}.
	\end{equation}
	Denote by  $\upepsilon$ the (intensity of) electric field pointing upward of the plane. Then, Faraday's law states that
	\begin{equation*}
		\pd_t \vH = \grad^\perp \upepsilon,
	\end{equation*}
	where $\grad^\perp \equiv (-\pd_2, \pd_1)$ in Cartesian's coordinates.
	If we denote by $e$ the (scalar) electric field in the plasma region $\Omt$, then the Rankine–Hugoniot-type jump condition on the free interface can be written as
	\begin{equation*}
		(\upepsilon - e) \vn^{\perp} = -\fks (\vH - \vh) \qq{on} \Gmt,
	\end{equation*}
	where $\vn^\perp = (-n^2, n^1)$ in Cartesian's coordinates. On the other hand, since the plasma is assumed to be a perfect conductor, Ohm's law reads that
	\begin{equation*}
		e + \vv\smallwedge\vh = 0,
	\end{equation*} 
	which yield 
	\begin{equation*}
		\upepsilon \vn^{\perp} = e\vn^{\perp} - \fks (\vH-\vh) = - (\vv\smallwedge\vh)\vn^\perp - \fks (\vH-\vh) = -\fks\vH + (\vh\vdot\vn)\vv \qq{on} \Gmt.
	\end{equation*}
	Hence, one obtains (recall that $\Curl \vH \equiv 0$ in $\cV_t$)
	\begin{equation*}
		\begin{split}
			\int_{\cV_t} \vH\vdot\pd_t\vH \dd{x} &= \int_{\cV_t} \vH \vdot \grad^\perp\upepsilon \dd{x} = \int_{\Gmt} - \upepsilon\vn^{\perp} \vdot\vH \dd{\ell} + \int_{\cS} \upepsilon \vN \smallwedge \vH \dd{\ell} \\
			&= \int_{\Gmt} \fks\abs{\vH}^2 - (\vh\vdot\vn)(\vv\vdot\vH) \dd{\ell} + \int_{\cS} \cJ\upepsilon \dd{\ell},
		\end{split}
	\end{equation*}
	which, together with \eqref{dt Ephy pre} and \eqref{BC}, yield that
	\begin{equation}\label{dt Ephy}
		\dv{t} E_{\text{phy}} = \int_{\cS} \cJ\upepsilon \dd{\ell}.
	\end{equation}
	Here we remark that the integral on the right hand side of \eqref{dt Ephy pre} is exactly the power of energy input caused by the surface current $\cJ$. Particularly, when $\cJ \equiv 0$, the physical energy is conserved for all time.
	
	\subsection{Analysis of Circular Flows}\label{sec circ flow}
	
	\subsubsection{Ill-posedness in the Lagrangian framework}
	We now briefly review the counterexamples given in \cite{Ebin1987} and \cite{Hao-Luo2020}. First, through introducing the flow map $\vX$:
	\begin{equation*}
		\dv{t} \vX(t, y) \coloneqq \vv\qty\big[t, \vX(t, y)] \qc \vX(0, y) = y \qfor y \in \Om_{\restriction_{t=0}},
	\end{equation*}
	the plasma-vacuum problem with $\vH\equiv\vb{0}$ can be written as
	\begin{equation}\label{Lag form}
		\begin{cases*}\displaystyle
			\dv[2]t \vX = \qty{\grad\lap^{-1}\qty[\tr(\grad\vv)^2-\tr(\grad\vh)^2]}\circ\vX + \qty\big[\qty(\vh\vdot \grad)\vh]\circ\vX, \\
			\vh\qty[t, \vX(t, y)] = \grad_y\vX(t, y) \vdot \vh(0, y).
		\end{cases*}
	\end{equation} 
	On the other hand, it is clear that the following data give a (stationary) solution to the plasma-vacuum problem \eqref{PV problem}:
	\begin{equation}\label{cir back sol}
		\Omt \equiv B_1 \subset \R^2 \qc \vH \equiv \vb{0} \qc \vv = r\fkV\etheta, \qand \vh = r\fkH\etheta,
	\end{equation}
	where $r \coloneqq \abs{x}$, $B_1=\{r < 1\}$ is the unit disk, $(\er, \etheta)$ is the standard orthonormal frame in polar coordinates, and $\fkV, \fkH$ are both constants. Particularly, for such circular background solutions, the equations \eqref{Lag form} can be rewritten as (here $\dot{\vX} \coloneqq \dv*{\vX}{t}$)
	\begin{equation}\label{Lag eqn}
		\ddot{\vX} = \scZ(\dot{\vX}, \vX),
	\end{equation}
	where the (multilinear) functional $\scZ$ depends on the parameters $\fkV$ and $\fkH$. From now on, for the simplicity of notations, we denote $\vX$ to be the flow induced by $\vv\equiv r\fkV\etheta$. Consider the linearized equation for \eqref{Lag eqn}
	\begin{equation}\label{Lag lin}
		\ddot{\vY} = (\var{\scZ})_{\restriction{(\vX, \dot{\vX})}}(\vY, \dot{\vY}),
	\end{equation}
	with initial data
	\begin{equation*}
		\vY_n(0, y) = y \qc \dot\vY(0, y) = \vv(y) + \fkV e^{-n^{\frac{1}{4}}} \vw_n(y), \qfor y \in B_1,
	\end{equation*}
	where in the polar coordinates, $\vw_n$ is given by
	\begin{equation*}
		\vw_n(r, \theta) = \mqty(r^n\cos(n\theta) \\ -r^n\sin(n\theta)).
	\end{equation*}
	It follows from \cite[Corollary 5.50]{Ebin1987} and \cite[Theorems 6.1-6.2]{Hao-Luo2020} that (here $H^k$ represents the standard Sobolev space, and $k \ge 2$ is an integer)
	\begin{equation*}
		\norm{(\vX-\vY_n,\; \dot{\vX}-\dot{\vY_n})(t)}_{H^{k}(B_1)} \xrightarrow{n\to\infty} \infty \qq{for any} t > 0,
	\end{equation*}
	while
	\begin{equation*}
		\mqty(\vY_n(t=0) \\ \dot{\vY}_n(t=0)) \xrightarrow{n\to\infty} \mqty(\vX(t=0) \\ \dot{\vX}(t=0)) \qin C^\infty(B_1) \times C^{\infty}(B_1),
	\end{equation*}
	provided that $\fkV > 0$ and $0 \le \fkH \ll \fkV$. Such constructions are derived from the analysis on the eigenvalues and eigenfunctions of the linear (differential) operator $(\var{\scZ})_{\restriction(\vX, \dot{\vX})}$, which has unstable eigenfunctions with arbitrarily large eigenvalues. Particularly, this implies the high instability of the linear problem \eqref{Lag lin}, and hence the non-continuous dependence on initial data for the nonlinear plasma-vacuum interface problem \eqref{Lag form} in the standard Sobolev spaces, which leads to the ill-posedness of \eqref{PV problem} (when $\vH\equiv\vb{0}$) in the Lagrangian framework.
	
	\subsubsection{Revisit through perturbing waves}
	In addition to the Lagrangian approach, we presently consider the linearized problem of \eqref{PV problem} in the Eulerian framework. For the sake of uniformity, we still assume that $\vH \equiv \vb{0}$. The detailed derivations will be given in Appendix \ref{sec App A}. First, the linearized equations and boundary conditions for \eqref{MHD} and \eqref{BC} can be written as
	\begin{subequations}\label{lin PV}
		\begin{equation}\label{lin MHD}
			\begin{cases*}
				\pd_t\wt{\vv}+(\vv\vdot\grad)\wt{\vv}+(\wt{\vv}\vdot\grad)\vv+\grad\wt{p} = (\wt{\vh}\vdot\grad)\vh+(\vh\vdot\grad)\wt{\vh} &in $\Omt$, \\
				\pd_t\wt{\vh} + (\vv\vdot\grad)\wt{\vh}+(\wt{\vv}\vdot\grad)\vh = (\wt{\vh}\vdot\grad)\vv+(\vh\vdot\grad)\wt{\vv} &in $\Omt$, \\
				\dive\wt{\vv} = 0 = \dive\wt{\vh} &in $\Omt$,
			\end{cases*}
		\end{equation}
		and
		\begin{equation}\label{lin BC}
			\begin{cases*}
				\pd_t\wt{\fks} + (\vv\vdot\grad)\wt{\fks} = \wt{\vv}\vdot\vn + \wt{\fks} (\grad_\vn\vv\vdot\vn) &on $\Gmt$, \\
				\wt{p} = -\alpha \grad_\vtau\grad_\vtau\wt{\fks} - \qty(\alpha\varkappa^2 + \grad_\vn p)\wt{\fks} &on $\Gmt$, \\
				\wt{\vh}\vdot\vn = \grad_\vtau\qty[\wt{\fks}(\vh\vdot\vtau)] &on $\Gmt$,
			\end{cases*}		
		\end{equation}
		where $(\wt{\vv}, \wt{\vh}, \wt{p}, \wt{\fks})$ are respectively the linearized variables of $(\vv, \vh, p, \fks)$, $\vtau$ is the unit tangential fields of $\Gmt$, $\grad$ is the covariant derivative on $\R^2$, and the new variables $(\wt{\vv}, \wt{\vh}, \wt{p}, \wt{\fks})$ are regarded as functions defined in $\Omt$ and $\Gmt$ for each fixed time moment.
	\end{subequations}
	Take the background solutions as \eqref{cir back sol}, and consider the perturbation variables in the form (here $k$ is a non-zero integer, $c \in \mathbb{C}$ is a complex number representing the phase velocity of $k$-wave)
	\begin{equation}\label{ptb profile}
		\wt{\vv} = \wh{\vv}(r)e^{ik(\theta-ct)}\qc \wt{\vh}=\wh{\vh}(r)e^{ik(\theta-ct)} \qc \wt{p} = \wh{p}(r) e^{ik(\theta-ct)}, \qand \wt{s} = \wh{s}e^{ik(\theta-ct)},
	\end{equation}
	where $\wh{s} \in \mathbb{C}\setminus\{0\}$ is a non-zero constant. We assume here that $\wh{s}\neq 0$ since we plan to study the surface wave, otherwise the problem will be reduced to a fixed boundary one. Since $\wt{\vv}$ and $\wt{\vh}$ are both solenoidal, the linear equations \eqref{lin MHD} can be reduced to (cf. Appendix \ref{sec App A})
	\begin{equation}\label{lin wave pur}
		\begin{cases*}\displaystyle
			(c-\fkV)\pd_r\qty[r\pd_r(r\wh{v}^r)] + \fkH \pd_r\qty[r\pd_r(r\wh{h}^r)] = k^2(c-\fkV)\wh{v}^r + k^2\fkH\wh{h}^r, \\
			(\fkV-c)\wh{h}^r = \fkH\wh{v}^r.
		\end{cases*}
	\end{equation}
	Taking the new variables 
	\begin{equation}\label{def z}
		s \coloneqq \log{r} \qand z (s)\coloneqq (r\wh{v}^r)_{\restriction r=e^s},
	\end{equation}
	the equations \eqref{lin wave pur} can be rewritten as
	\begin{equation}\label{ODE}
		\dv[2]{s}z = k^2 z, \qif c-\fkV\neq 0.
	\end{equation}
	Concerning the boundary values of $z$, it is natural to impose the compatibility condition
	\begin{equation*}
		z(-\infty) \coloneqq \lim_{s\to-\infty} z(s) = 0 \times \wh{v}^r(0) = 0.
	\end{equation*}
	On the other hand, it follows from \eqref{lin BC}\textsubscript{1} that
	\begin{equation*}
		ik\wh{s}(\fkV - c) = \wh{v}^r \qfor r = 1.
	\end{equation*}
	Thus, we may take the renormalization
	\begin{equation}\label{renormalize s}
		ik\wh{s} = 1,
	\end{equation}
	which yields the boundary values for \eqref{ODE}:
	\begin{equation}\label{ODE BC}
	 	 z(0) = \fkV - c \qand z(-\infty) = 0.
	\end{equation}
	Particularly, the ODE \eqref{ODE} with boundary condition \eqref{ODE BC} admit a unique solution
	\begin{equation}\label{sol z}
		z(s) = (\fkV-c) \exp(\abs{k}s) \qfor s\le 0.
	\end{equation}
	Moreover, \eqref{lin BC}\textsubscript{2} implies the dispersive relation (cf. Appendix \ref{sec App A}):
	\begin{equation}\label{disp}
		k\qty[c-\frac{\abs{k}-1}{\abs{k}}\fkV]^2 = (\abs{k}-1)\qty[\alpha(\abs{k}+1)+\fkH^2-\frac{1}{\abs{k}}\fkV^2],
	\end{equation}
	which is an algebraic equation for $c$.
	
	In view of the perturbation profile \eqref{ptb profile}, the $k$-wave is unstable iff the solution $c_k$ to \eqref{disp} has non-zero imaginary part. Particularly, when $\alpha = 0$ and $\fkH = 0$, the equation \eqref{disp} has non-real roots for all $k$-waves with $\abs{k} > 1$. Namely, all large waves are unstable, which principally leads to the Rayleigh-Taylor instability (ill-posedness) of the one-phase free-boundary Euler equations violating the Rayleigh-Taylor sign condition \eqref{RT} (see \cite{Ebin1987}). When $\alpha > 0$, i.e., there exists surface tension, all large waves are stable, which leads to the well-posedness of capillary surface waves (cf. \cite{Cheng-Coutand-Shkoller2008,Shatah-Zeng2011}). Concerning the MHD flows, if $\alpha = 0$ and $\fkH \neq 0$, it can be easily derived from \eqref{disp} that the perturbing waves with sufficiently large wave numbers are all stable, which implies that the well-posedness is hopefully true. Such analysis indicates that tangential magnetic fields (no matter how tiny they are) can stabilize the motion of surface waves, which were not obtained in the linear analysis of flow maps in the Lagrangian framework (cf. \cite{Hao-Luo2020}).
	
	As a separate note, for large wave numbers, it holds that $\alpha(1+\abs{k}) > \fkH^2$, which yields that the stabilization effect of capillary forces is stronger than that of the tangential magnetic fields.
	
	\subsubsection{Revisit through geometrical analysis}
	Due to the fundamental theorem of plane curves, each curve is uniquely characterized by its curvature. Thus, heuristically speaking, the evolution equations for the curvature largely reflect the stability/well-posedness of the free boundary problems. 
	
	To proceed, we denote by $\Dt\coloneqq \pd_t + (\vv\vdot\grad)$ the material derivative along the fluid particle trajectory. Then, through plugging \eqref{PV problem} (with $\vH \equiv \vb{0}$) into the second order evolution equation for $\varkappa$, one will obtain a wave-type equation (cf. \S\ref{sec geo anal} for details):
	\begin{equation}\label{wave eqn formal}
		\Dt\Dt\varkappa + \alpha (-\slashed{\lap})^\frac{3}{2} \varkappa + \abs{\vh}^2 (-\slashed{\lap})\varkappa + (-\grad_\vn p) (-\slashed{\lap})^{\frac{1}{2}}\varkappa \approx 0,
	\end{equation}
	where $\slashed{\lap}$ is the Laplace-Beltrami operator on $\Gmt$. It is straightforward to see that surface tension ($\alpha>0$), non-degenerate tangential magnetic fields ($\abs{\vh} > 0$), and the Rayleigh-Taylor sign condition \eqref{RT} can all ensure that the equation \eqref{wave eqn formal} is strictly hyperbolic, whence hopefully solvable. It is clear that the stabilizing effect of the surface tension is the strongest, as it corresponds to the leading order term in \eqref{wave eqn formal}. Moreover, provided that $\alpha = 0$ and $\inf_{\Gmt}\abs{\vh}>0$, the Rayleigh-Taylor sign condition \eqref{RT} will only be reflected by a lower-order term, meaning the sign of $(\grad_\vn p)$ is insignificant for such scenario. Specifically, the problem \eqref{PV problem} with initial data \eqref{cir back sol} can be well-posed/stable, at least for a short period.
	
	It is interesting to note that the initial data \eqref{cir back sol} can lead to high-instability/ill-posedness in the Lagrangian framework \eqref{Lag eqn}, while the interface problem \eqref{PV problem} with initial data \eqref{cir back sol} is stable/well-posed in the Eulerian scheme. Such phenomena imply that only the flow maps is insufficient to capture the whole information of surface waves. Indeed, the dispersive relation \eqref{disp} indicates that the magnetic field can stabilize the large waves, while the linear (differential) operator in \eqref{Lag lin} still has eigenvalues diverge to infinity. This implies that the effect of tangential magnetic fields on the interfacial waves is not fully reflected through analyzing merely the evolution of flow maps. For the free-boundary Euler equations, it is shown in \cite[\S5]{Shatah-Zeng2008-Geo} that occasionally the free interface is smoother than the flow map. These observations suggest that, in contrast to the classical and prevalent Lagrangian framework, the Eulerian scheme is potentially better suited to the study of free interface problems for incompressible flows.
	
	\subsection{Main Results}
	
	In the following, we aim at establishing the local well-posedness theories in the standard Sobolev spaces \emph{without} assuming the Rayleigh-Taylor sign condition \eqref{RT}. Instead, we will show them under the hypothesis that either \emph{total} magnetic fields are uniformly non-degenerate, i.e.,
	\begin{equation}\label{bc h}
		\abs{\vh} + \abs{\vH} > 0 \qq{on} \Gmt,
	\end{equation}
	or the surface tension prevails (i.e., $\alpha>0$).
	Regarding the physical energy \eqref{def Ephy} and the conservation law \eqref{dt Ephy}, the $L^2$-based standard Sobolev spaces could be adapted to the free boundary problem \eqref{PV problem}. For the sake of technical convenience, throughout this manuscript, it is always assume that the plasma region $\Omt$ is simply connected and the vacuum region $\cV_t$ has only one ``hole'' (i.e., the Betti number of $\cV_t$ is one).
	
	\subsubsection{Well-posedness without surface tension}
	
	If $\alpha = 0$ and the total magnetic fields are non-degenerate on the free boundary, it holds that
	
	\begin{theorem}\label{thm non-deg}
		Let $s \ge 3$ be a constant with $2s\in \mathbb{N}$. Assume that $\Gamma_0 \in C^2 \cap H^{s+\frac{1}{2}}$, $\cS \in H^{s-\frac{1}{2}}$, and the prescribed surface current $\cJ$ satisfies
		\begin{equation*}
			\cJ \in C^0_t H^{s-\frac{1}{2}}(\cS) \cap C^1_t H^{s-\frac{3}{2}}(\cS).
		\end{equation*}
		Let $(\vv_0, \vh_0) \in H^s(\Omega_0)$ be two solenoidal vector fields, for which $\vh_0$ is tangential to $\pd\Omega_0$. Denote by $\vH_0$ the magnetic field in the initial vacuum region $\cV_0$ induced by the surface current $\cJ_{\restriction_{t=0}}$. If there exists some positive constant $\lambda_0$ such that
		\begin{equation*}
			\abs{\vH_0} + \abs{\vh_0} \ge 2\lambda_0 \qq{on} \Gamma_0,
		\end{equation*}
		then, there is a positive constant $T$, so that the problem \eqref{PV problem} with $\alpha = 0$ admits a solution in the following space for $t \in [0, T]$:
		\begin{equation*}
			\vv \in C^0_t H^s(\Omt) \qc \vh \in C^0_t H^s(\Omt) \qc \Gmt \in C^0_t H^{s+\frac{1}{2}}, \qand \vH \in C^0_t H^s(\cV_t),
		\end{equation*}
		which also satisfies the relation:
		\begin{equation*}
			\inf_{\Gmt} \qty(\abs{\vh} + \abs{\vH}) \ge \lambda_0 \qfor t \in [0, T].
		\end{equation*}
		Furthermore, if $s \ge 4$, the solution is unique and the data-to-solution map is Lipschitz in the space
		\begin{equation*}
			\vv \in C_0^t H^{s-1}(\Omt) \qc \vh \in C^0_t H^{s-1}(\Omt) \qc \Gmt \in C^0_t H^{s-\frac{1}{2}}, \qand \vH \in C^0_t H^{s-1}(\cV_t).
		\end{equation*}
		Namely, the plasma-vacuum problem \eqref{PV problem} is locally well-posed and weakly Lipschitz with respect to initial data.
	\end{theorem}
	
	\begin{remark}
		The above theorem indicates that non-trivial tangential magnetic fields can indeed stabilize the motion of the plasma-vacuum interface, which verifies the heuristic arguments presented in \S\ref{sec circ flow}. In contrast to the three-dimensional case (cf. \cite{Sun-Wang-Zhang2019,liu-Xin2025free}), whose well-posedness requires the tangential magnetic fields to be non-collinear, the two-dimensional scenario only needs the non-degeneracy of total magnetic fields.
	\end{remark}
	
	\subsubsection{Well-posedness with surface tension}
	When $\alpha > 0$, i.e., the surface tension is non-vanishing, there holds the following result:
	\begin{theorem}
		Let $s \ge 3$ be a constant with $2s \in \mathbb{N}$ ($s \neq 3.5$). Assume that $ \Gamma_0 \in C^2 \cap H^{s+1}$, the fixed wall $\cS \in H^{s-\frac{1}{2}}$, and the prescribed surface current $\cJ$ satisfies 
		\begin{equation*}
			\cJ \in C^0_t H^{s-\frac{1}{2}}(\cS) \cap C^1_t H^{s-\frac{3}{2}}(\cS).
		\end{equation*}
		Given two solenoidal vector fields $(\vv_0, \vh_0) \in H^{s}(\Omega_0)$, where $\vh_0$ is tangential to $\pd\Omega_0$, there exists a positive constant $T$, so that the plasma-vacuum interface problem \eqref{PV problem} admits a solution in the following space for $t \in [0, T]$:
		\begin{equation*}
			\vv \in C^0_t H^{s}(\Omt) \qc \vh \in C^0_t H^{s}(\Omt) \qc \Gmt \in C^0_t H^{s+1}, \qand \vH \in C^0_t H^s(\cV_t).
		\end{equation*}
		Moreover, if $s \ge 4.5$, the solution is unique, and the data-to-solution map is Lipschitz in the space
		\begin{equation*}
			\vv \in C_0^t H^{s-\frac{3}{2}}(\Omt) \qc \vh \in C^0_t H^{s-\frac{3}{2}}(\Omt) \qc \Gmt \in C^0_t H^{s-\frac{1}{2}}, \qand \vH \in C^0_t H^{s-\frac{3}{2}}(\cV_t).
		\end{equation*}
		Namely, the problem \eqref{PV problem} is locally well-posed and Lipschitz (with respect to weaker norms) in initial data.
	\end{theorem}
	
	\begin{remark}
		The Lipschitz continuity here is interpreted through the equivalent formulations of the problem \eqref{PV problem} defined in a fixed reference frame (cf. the discussions in \S\ref{sec it} and \S\ref{sec var est} for details). The requirement that $s \neq 3.5$ is due to technical reasons.
	\end{remark}
	
	\begin{remark}
		Under the presence of surface tension, there is no requirement on the tangential magnetic fields, as capillary forces have stronger stabilizing effect on the interfacial waves. This result coincides with the stability and geometrical analysis presented in \S\ref{sec circ flow}.
	\end{remark}

	\subsubsection{Uniform energy estimates and vanishing surface-tension limits}
	Concerning the a priori estimates, there holds the following theorem:
	\begin{theorem}\label{thm eng est}
		Let $m \ge 0$ be an integer and $(\vv, \vh, \vH, \Gmt)$ a solution to \eqref{PV problem} for $t \in [0, T]$. Define the $H^{m+3}$-size as:
		\begin{equation*}
			\begin{split}
				M^m_\alpha \coloneqq\, &\norm{\vv}_{H^{m+3}(\Omt)}^2 + \norm{\vh}_{H^{m+3}(\Omt)}^2 + \qty(\abs{\pd_t \cJ}_{H^{m+\frac{3}{2}}(\cS)}^2 + \abs{\cJ}_{H^{m+\frac{5}{2}}(\cS)}^2)\cdot \qty(1 + \abs{\varkappa}_{H^{m+\frac{3}{2}}(\Gmt)}^2) + \\
				&+ \alpha\abs{\varkappa}_{H^{m+2}(\Gmt)}^2 + \abs{(\vh\vdot\grad)\varkappa}_{H^{m+\frac{1}{2}}(\Gmt)}^2 + \abs{\varkappa}_{H^{m+1}(\Gmt)}^2,
			\end{split}
		\end{equation*}
		where $\varkappa$ is the signed curvature of $\Gmt$. Then, for each $m \ge 0$, there exists an energy functional $E^m_\alpha$ satisfying the following properties (here $\cQ$ represents generic polynomials):
		\begin{itemize}
			\item {\bf Propagation estimate:}
			\begin{equation*}
				\abs{\dv{t}E^m_\alpha} \le \cQ\qty(M^m_\alpha).
			\end{equation*}
			\item {\bf Energy bound:}
			\begin{equation*}
				\abs{E^m_\alpha} \le \cQ\qty(M^m_\alpha).
			\end{equation*}
			\item {\bf Energy coercivity:}
			\begin{itemize}
				\item if $\alpha > 0$, there holds
				\begin{equation*}
					M^m_\alpha \le \cQ_\alpha\qty(E^m_\alpha, \abs{\pd_t \cJ}_{H^{m+\frac{3}{2}}(\cS)}, \abs{\cJ}_{H^{m+\frac{5}{2}}(\cS)});
				\end{equation*}
				\item if $\alpha = 0$ and there exists a constant $\lambda_0 > 0$ so that
				\begin{equation*}
					\inf_{\Gmt}\qty(\abs{\vh}+\abs{\vH}) \ge \lambda_0 \qq{for all} t \in [0, T],
				\end{equation*}
				then it holds that
				\begin{equation*}
					M^m_{\alpha=0} \le  \cQ_{\lambda_0}\qty(E^m_{\alpha=0}, \abs{\pd_t \cJ}_{H^{m+\frac{3}{2}}(\cS)}, \abs{\cJ}_{H^{m+\frac{5}{2}}(\cS)});
				\end{equation*}
				\item if $\alpha = 0$, $\cJ \equiv 0$, and there is a constant $c_0 > 0$ so that
				\begin{equation*}
					- (\grad_\vn p) > c_0 \, \text{ on } \, \Gmt \qq{for all} t \in [0, T],
				\end{equation*}
				then it follows that
				\begin{equation*}
					M^m_{\alpha=0} \le \cQ_{c_0}\qty(E^m_{\alpha=0}).
				\end{equation*}
			\end{itemize}
		\end{itemize}
	\end{theorem}
	
	As corollaries, there hold the following vanishing surface tension limit results:
	\begin{cor}[$\alpha \to 0$ limit under non-degenerate magnetic fields]
		Assume that $\alpha \ge 0$, $s \ge 3$, $2s \in \mathbb{N} (s \neq 3.5)$, and $\cJ \in C^0_tH^{s-\frac{1}{2}}(\cS) \cap C^1_t H^{s-\frac{3}{2}}(\cS)$. For each initial data $(\vv_0, \vh_0, \Gamma_0)$ satisfying $(\vv_0, \vh_0) \in H^{s}(\Omega_0)$, $\dive\vv_0 = \dive\vh_0 = 0$, $\Gamma_0 \in H^{s+1}$, and $\inf_{\Gamma_0}\qty(\abs{\vh_0}+\abs{\vH_0}) \ge 2\lambda_0$ for some positive constant $\lambda_0$, there exists a constant $T > 0$ independent of $\alpha$, so that the problem \eqref{PV problem} is solvable for $t \in [0, T]$. Moreover, as $\alpha \to 0$, by passing to a subsequence, the solutions to \eqref{PV problem} with surface tension weakly converges to a solution with $\alpha = 0$ in the space $\Gmt \in H^{s+\frac{1}{2}} $, $\vv \in H^{s}(\Omt) $, $\vh \in H^{s}(\Omt) $, and $\vH \in H^s(\cV_t)$.
	\end{cor}
	\begin{remark}
		When the magnetic fields are uniformly non-degenerate on $\Gmt$, the uniform energy estimates can be refined so that the regularity indices are not necessarily being integral (cf. \S\S\ref{sec est alpha>0}-\ref{sec est h non-degenerate} for details).
	\end{remark}
	
	\begin{cor}[$\alpha \to 0$ limit under the Rayleigh-Taylor sign condition]\label{cor 2}
		Assume that $\alpha \ge 0$, $s \ge 3$, $s \in \mathbb{N}$, and $\cJ \equiv 0$. Let $(\vv_0, \vh_0, \Gamma_0)$ be initial data satisfying $(\vv_0, \vh_0) \in H^{s}(\Omega_0)$, $\dive\vv_0 = \dive\vh_0 = 0$, $\Gamma_0 \in H^{s+1}$, and $(-\grad_\vn q_0) \ge 2c_0$ for some positive constant $c_0$, here $q_0$ is a multiplier-type pressure defined through solving the Dirichlet problem
		\begin{equation*}
			\begin{cases*}
				-\lap q_0 = \tr{(\grad\vv_0)^2- (\grad\vh_0)^2} &in $\Omega_0$, \\
				q = 0 &on $\Gamma_0$.
			\end{cases*}
		\end{equation*}
		Then, for each $\alpha > 0$, there is a constant $T > 0$, which is independent of $\alpha$, so that the problem \eqref{PV problem} is solvable for $t \in [0, T]$. Moreover, as $\alpha \to 0$, by taking a subsequence, the solutions to \eqref{PV problem} with positive surface tension converges weakly to a solution without surface tension in the space $\Gmt \in H^s$, $\vv \in H^s(\Omt)$, and $\vh \in H^{s}(\Omt)$.
	\end{cor}
	
	\begin{remark}
		The well-posedness theory of \eqref{PV problem} under merely the R-T condition is beyond the scope of this manuscript. For results under the additional assumption that $\cJ \equiv 0$, one can refer to \cite{zhao2024local} for the 3D case in a modal domain, and \cite{ifrim2024sharp} for results applicable to general bounded plasma domains in high space dimensions. (Note that the state spaces in \cite{ifrim2024sharp} is different from those in this manuscript.)
	\end{remark}
	
	\subsection{Crucial Strategies}
	
	Inspired by the geometric perspectives of {\sc Shatah and Zeng} \cite{Shatah-Zeng2008-Geo,Shatah-Zeng2008-vortex,Shatah-Zeng2011}, we utilize the evolution of signed curvature of $\Gmt$ to characterize the motion of the free interface. Indeed, the evolution equation for the curvature $\varkappa$ can be written as
	\begin{equation*}
		\Dt\varkappa = -\vn \vdot \slashed{\lap} \vv - 2\varkappa \slashed{\grad}\vdot\vv,
	\end{equation*}
	where $\slashed{\lap}$ and $(\slashed{\grad}\vdot)$ are respectively the Laplace-Beltrami operator and tangential divergence on $\Gmt$. Thus, one can plug the system \eqref{PV problem} into the second order evolution equation for $\varkappa$, and, at the principal level, obtain a wave-type equation:
	\begin{equation}\label{evo eqn kp formal}
		\Dt\Dt\varkappa + \alpha (-\slashed{\lap})^\frac{3}{2} \varkappa + \qty(\abs{\vh}^2+\abs{\vH}^2) (-\slashed{\lap})\varkappa + (\grad_\vn\wt{q}-\grad_\vn q)(-\slashed{\lap})^{\frac{1}{2}}\varkappa = \boxed{\text{L.O.T.}} \qq{on} \Gmt,
	\end{equation}
	where $q$ and $\wt{q}$ are both multiplier-type pressures (whose definitions will be given by \eqref{def q} and \eqref{def wtq} respectively), and ``L.O.T.'' stands for lower order terms. It is routine to see that the three additional assumptions in the ``energy coercivity'' part of Theorem \ref{thm eng est} are all to ensure that \eqref{evo eqn kp formal} is strictly hyperbolic. Particularly, concerning the stabilizing effect on the interfacial waves, there hold the following orders: \emph{surface tension} $>$ \emph{non-degenerate tangential magnetic fields} $>$ \emph{the Rayleigh-Taylor sign condition}. 
	
	Note that both $\vv$ and $\vh$ are solenoidal vector fields, standard div-curl analysis implies that they can be uniquely determined by their boundary behaviors and curls. The evolution of the curvature yields the information about the interface dynamics, and the vorticity and current satisfy transport-type equations. Therefore, the whole problem \eqref{PV problem} can be resolved through such divide-and-conquer strategies.
	
	This manuscript is structured as follows. In \S\ref{sec geo anal}, we derive the adapted evolution equations for curvatures through geometrical perspectives. In \S\ref{sec tech preparation}, we introduce reference frames for the convenience of analysis. After that, we collect some useful technical results. \S\ref{sec lin est} is devoted to the energy estimates for linearized problems, which are crucial for later arguments. In \S\ref{sec a priori est}, we show energy estimates for the non-linear problem \eqref{PV problem} under multiple stability conditions.  \S\ref{sec construct ite} is dedicated to the solution construction schemes. We first pull-back the evolution problem in moving domains to those defined in the preassigned reference frame, in order that one can naturally utilize the Picard iteration scheme. We then demonstrate the equivalence between solutions to the problems defined in the reference frame and those to the original plasma-vacuum problem \eqref{PV problem}. The estimates for iteration maps will be discussed in \S\ref{sec est ite}.
	
	\section{Geometrical Analysis}\label{sec geo anal}
	Suppose that $\Gmt = \pd\Omt$ is a $C^2$-curve, with $\vn$ being the outward unit normal. Denote by $\vtau$ the unit tangent vector field of $\Gmt$, with the orientation that $\vtau$ is the counterclockwise rotation of $\vn$ with angle $\frac{\pi}{2}$. Such an orientation is chosen so that for $\Omt = B_1 \subset \R^2$, one has that $\vn = \hat{\vb{e}}_r$ and $\vtau = \hat{\vb{e}}_{\theta}$. Then, the curvature of $\Gmt$ is given by:
	\begin{equation}
		\varkappa \coloneqq \grad_{\vtau} \vn \vdot \vtau.
	\end{equation}
	Thus, the fact that $\vtau$ and $\vn$ both have unit lengths yields 
	\begin{equation}\label{deriv gmt}
		\grad_\vtau \vn = \varkappa \vtau \qand \grad_{\vtau}\vtau = - \varkappa\vn.
	\end{equation}
	
	\subsection{Evolution of Curves}

	Now, assume that $\Gmt$ evolves with velocity $\vbu : \Gmt \to \R^2$, and we shall derive the evolution equations for $\vtau$ and $\vn$. Note here that $\vbu$ is not necessarily being the fluid velocity $\vv$. (Of course, $\vv$ is one of the evolution velocities in the context of free boundary problems, but $(\vv + \sigma\vb*{\tau})$ is still a candidate for any function $\sigma$ defined on $\Gmt$.)
	
	Denote by $\cDt$ the evolution velocity with respect to $\vbu$:
	\begin{equation}
		\cDt \coloneqq \pd_t + \grad_\vbu.
	\end{equation}
	Then, the fact that $\vtau$ has unit length implies
	\begin{equation*}
		\cDt\vtau \vdot \vtau \equiv 0.
	\end{equation*}
	Indeed, by invoking the coordinate maps, it is routine to derive that (see also \cite[eq. (3.3)]{Shatah-Zeng2008-Geo})
	\begin{equation}\label{dt tau}
		\cDt\vtau = (\grad_\vtau \vbu \vdot \vn)\vn.
	\end{equation}
	In particular, there holds the following commutator formula:
	\begin{equation}\label{comm Dt d slash}
		\comm{\cDt}{\grad_\vtau} \coloneqq \cDt\grad_\vtau - \grad_\vtau\cDt = \grad_{(\cDt\vtau-\grad_\vtau\vbu)} = -(\grad_\vtau\vbu\vdot\vtau)\grad_\vtau
	\end{equation}
	The fact that $\vn$ is the unit normal of $\Gmt$ yields
	\begin{equation}\label{dt vn}
		\cDt \vn = - (\grad_\vtau\vbu \vdot \vn)\vtau.
	\end{equation}
	Moreover, for the length element of $\Gmt$, one has
	\begin{equation}\label{Dt dl}
		\cDt \dd{\ell} = (\grad_\vtau \vbu \vdot \vtau) \dd{\ell}.
	\end{equation}
	One may refer to \cite[\S3.1]{Liu-Xin2023} for the detailed derivations using local coordinate maps. 
	
	Due to the boundary condition for the effective pressure $p$ and the fundamental theorem of plane curves, it would be useful to compute the evolution of the curvature. Indeed, it is routine to derive that
	\begin{equation}\label{Dt kappa}
		\cDt\varkappa = \cDt(\grad_{\vtau}\vn \vdot \vtau) = -\vn \vdot \grad_\vtau\grad_\vtau\vbu - 2\varkappa(\grad_\vtau\vbu\vdot\vtau).
	\end{equation}
	In order to plug the system \eqref{PV problem} into the boundary evolution, one needs to consider the second order evolution equation for $\varkappa$. It follows from the standard calculations that
	\begin{equation}\label{Dt2 kappa}
		\begin{split}
			\cDt\cDt \varkappa = &-\grad_\vtau\grad_\vtau(\cDt\vbu\vdot\vn) + (\grad_\vtau\varkappa)\vtau\vdot\cDt\vbu - \varkappa^2 \vn\vdot\cDt\vbu + \\
			&+ 2 (\grad_\vtau\vbu \vdot\vn)(\vtau\vdot\grad_\vtau\grad_\vtau\vbu) + 4(\grad_\vtau\vbu\vdot\vtau)(\vn\vdot\grad_\vtau\grad_\vtau\vbu) + \\
			&+ 6\varkappa\abs{\grad_\vtau\vbu\vdot\vtau}^2 - 3\varkappa\abs{\grad_\vtau\vbu\vdot\vn}^2 \\
			\eqqcolon &-\grad_\vtau\grad_\vtau(\cDt\vbu\vdot\vn) + (\grad_\vtau\varkappa)\vtau\vdot\cDt\vbu - \varkappa^2 \vn\vdot\cDt\vbu + \fkR_0 (\vbu, \varkappa, \vtau, \vn).
		\end{split}
	\end{equation}
	
	\subsection{Evolution of the Interface}
	In the context of plasma-vacuum interface problems, one may take the evolution velocity of $\Gmt$ being exactly that of the plasma, i.e.,
	\begin{equation*}
		\vbu \equiv \vv \qq{on} \Gmt.
	\end{equation*}
	In other words, when restricted to $\Gmt$, the trajectory derivative $\cDt$ is now identical to the material derivative $\Dt$ along the fluid flows, where $\Dt$ is defined by
	\begin{equation}
		\Dt \coloneqq \pdv{t} + \grad_\vv.
	\end{equation}
	Then, by invoking \eqref{PV problem}, there hold
	\begin{equation}
		\Dt\vv = -\grad p + \grad_\vh \vh \qq{with} p = q + \alpha \cH \varkappa + \cH\qty(\tfrac{1}{2}\abs{\vH}^2),
	\end{equation}
	where $\cH$ is the harmonic extension into $\Omt$, i.e., $\cH f \equiv f_{\cH}$ satisfies
	\begin{equation}\label{def harmonic ext}
		\begin{cases*}
			-\lap (\cH f) = 0 &in $\Omt$, \\
			\cH f = f &on $\Gmt$;
		\end{cases*}
	\end{equation}
	and $q$ is the Lagrange multiplier through solving the Poisson equation:
	\begin{equation}\label{def q}
		\begin{cases*}
			-\lap q = \tr{(\grad\vv)^2 - (\grad\vh)^2} &in $\Omt$, \\
			q = 0 &on $\Gmt$.
		\end{cases*}
	\end{equation}
	Namely, when restricted to the free boundary $\Gmt$, the acceleration $\Dt\vv$ can be expressed as follows:
	\begin{equation*}
		\begin{split}
			(\Dt\vv)_{\restriction_{\Gmt}} = &- \grad\qty[q + \alpha\cH\varkappa + \cH\qty(\tfrac{1}{2}\abs{\vH}^2)] + \grad_{\vh}\vh   \\
			= &-\alpha(\grad_\vtau\varkappa)\vtau + \tfrac{1}{2}\grad_\vtau(\abs{\vh}^2)\vtau -\tfrac{1}{2}\grad_\vtau(\abs{\vH}^2)\vtau + \\
			&-(\grad_\vn q)\vn - \alpha(\grad_\vn\varkappa_{\cH})\vn - \tfrac{1}{2}\grad_\vn\cH(\abs{\vH}^2)\vn - \abs{\vh}^2\varkappa\vn.
		\end{split}
	\end{equation*}
	Denote by $\cN \coloneqq \grad_\vn \cH$ the Dirichlet-Neumann operator, i.e., for a function $f$ defined on $\Gmt$,
	\begin{equation*}
		\cN f \coloneqq (\grad_\vn f_{\cH})_{\restriction_{\Gmt}}.
	\end{equation*}
	Note that $\cN$ can be viewed as a first order differential operator, whose properties will be discussed in \S\ref{sec DN op}.	
	
	In particular, it follows that \eqref{Dt2 kappa} can be further computed as:
	\begin{equation}\label{DtDt kappa 1}
		\begin{split}
			\Dt\Dt\varkappa =\, & \alpha\grad_\vtau\grad_\vtau\cN\varkappa + \abs{\vh}^2\grad_\vtau\grad_\vtau\varkappa + \grad_\vtau\grad_\vtau\grad_\vn q + \tfrac{1}{2}\grad_\vtau\grad_\vtau\cN(\abs{\vH}^2) + \\
			&-\alpha\abs{\grad_\vtau\varkappa}^2 + \alpha\varkappa^2\cN\varkappa + \tfrac{5}{2}\grad_\vtau(\abs{\vh}^2)\grad_\vtau\varkappa - \tfrac{1}{2}\grad_\vtau(\abs{\vH}^2)\grad_\vtau\varkappa + \\
			&+ \varkappa^3\abs{\vh}^2 + 2\varkappa\grad_\vtau\grad_\vtau(\abs{\vh}^2) + \varkappa^2\grad_\vn q + \tfrac{1}{2}\varkappa^2 \cN(\abs{\vH}^2) +  \fkR_0(\vv, \varkappa, \vtau, \vn),
		\end{split}
	\end{equation}
	where $\fkR_0$ is given by \eqref{Dt2 kappa}. One can actually further simplify the seemingly complicated factor $\grad_\vtau\grad_\vtau\grad_\vn q$ as follows. First recall that the coordinate independent definition of Laplacian in $\R^2$ is the trace of the Hessian tensor. Hence, the restriction of $\lap f$ onto $\Gmt$ can be written as:
	\begin{equation}\label{lap rel}
		\begin{split}
			\lap f =\, &\grad^2_{(\vtau, \vtau)} f + \grad^2_{(\vn, \vn)} f \\
			=\, &\grad_\vtau\grad_\vtau f - \grad_{\grad_\vtau\vtau}f + \grad^2 f(\vn, \vn) \\
			=\, &\grad_\vtau\grad_\vtau f + \grad^2 f(\vn, \vn) + \varkappa\grad_\vn f,
		\end{split}
	\end{equation} 
	for a $C^2$-function $f$ defined in a neighborhood of $\Gmt$. By considering the harmonic extension, one may regard $\vn$ as a well defined vector field in $\Omt$. In other words, it can be computed that
	\begin{equation}\label{q1}
		\begin{split}
			\grad_\vtau\grad_\vtau\grad_\vn q =\, & \lap(\vn_{\cH} \vdot \grad q) - \grad^2(\grad_{\vn_{\cH}}q)(\vn, \vn) + \varkappa\grad_\vn(\grad_{\vn_{\cH}}q) \\
			=\, & \vn_{\cH} \vdot \lap\grad q + 2 \grad{\vn_{\cH}}\vb{\colon}\grad^2 q - \grad^2\vn_{\cH}(\grad q, \vn, \vn) - \grad^3 q(\vn, \vn, \vn) + \\
			&-2\grad_\vn \vn_{\cH} \vdot \grad_\vn\grad q - \varkappa\grad_\vn \vn_{\cH} \vdot \grad q - \varkappa\vn_{\cH}\vdot\grad_\vn\grad q \\
			=\, &\grad_{\vn}\qty[\lap q - \grad^2 q(\vn_\cH, \vn_\cH)] + 2\grad^2 q (\vn, \cN\vn) -\grad q \vdot \qty[\grad^2\vn_{\cH}(\vn, \vn)] \\
			&+ 2 \grad{\vn_{\cH}}\vb{\colon}\grad^2 q - 2 \grad^2 q(\cN\vn, \vn) - \varkappa\grad_{\cN\vn}q - \varkappa \grad^2 q (\vn, \vn).
		\end{split}
	\end{equation}
	It follows from \eqref{lap rel} that
	\begin{equation}\label{q2}
		\begin{split}
			\grad q \vdot \grad^2\vn_{\cH}(\vn, \vn) &= (\grad_\vn q) \vn \vdot \qty(\lap\vn_\cH - \varkappa\grad_{\vn}\vn_\cH - \grad_\vtau\grad_\vtau\vn) \\
			&= (\grad_\vn{q}) \vn \vdot \qty[\varkappa\cN\vn - (\grad_\vtau\varkappa)\vtau + \varkappa^2\vn] \\
			&=  (\grad_\vn{q})\qty[\varkappa \vn \vdot\cN\vn + \varkappa^2],
		\end{split}
	\end{equation}
	and
	\begin{equation}\label{q3}
		\begin{split}
			\grad_{\vn}\qty[\lap q - \grad^2 q(\vn_\cH, \vn_\cH)]  &= \grad_\vn (\varrho + \varkappa_{\cH}\grad_{\vn_{\cH}}q ) \\
			&= (\grad_\vn q)\cN\varkappa + \varkappa \grad^2 q(\vn, \vn) + \grad_\vn \varrho,
		\end{split}
	\end{equation}
	where $\varrho$ is an ancillary function defined by
	\begin{equation}\label{def varrho}
			\varrho \coloneqq \lap q - \grad^2 q(\vn_\cH, \vn_\cH) - \varkappa_{\cH}\grad_{\vn_{\cH}}q.
	\end{equation}
	In particular, $\varrho$ satisfies the property that
	\begin{equation}
		\varrho = \grad_\vtau \grad_\vtau q = 0 \qq{on} \Gmt.
	\end{equation}
	Therefore, one can calculate that
	\begin{equation}\label{tan lap dn q}
		\grad_\vtau\grad_\vtau\grad_\vn q = (\grad_\vn q)\cN\varkappa - (\grad_\vn q)\varkappa(\varkappa+\vn\vdot\cN\vn) - (\grad_{\cN\vn}q)\varkappa + 2\grad\vn_\cH\vb{\colon}\grad^2q + \grad_\vn \varrho.
	\end{equation}
	To handle the term $\grad_\vtau\grad_\vtau\cN(\abs{\vH}^2)$, one may first consider the harmonic extension of functions on $\Gmt$ into $\cV_t$ defined as follows:
	\begin{equation}\label{def har ext'}
		\begin{cases*}
			\lap \wt{\cH}f = 0 &in $\cV_t$, \\
			\wt{\cH} f = f &on $\Gmt$, \\
			\grad_\vN \wt{\cH} f = 0 &on $\cS$.
		\end{cases*}
	\end{equation}
	Define the Dirichlet-Neumann operator with respect to $\wt{\cH}$ as:
	\begin{equation}
		\wt{\cN} f \coloneqq -\vn\cdot\grad\wt{\cH}f.
	\end{equation}
	Indeed, it will be seen in \S\ref{sec DN op} that the two D-N operators $\cN$ and $\wt{\cN}$ coincide at the principle level. Thus, it would be natural to replace $\cN(\abs{\vH}^2)$ by $\wt{\cN}(\abs{\vH}^2)$. Namely, one can compute that
	\begin{equation}
			\cN\qty(\tfrac{1}{2}\abs{\vH}^2) = (\cN-\wt{\cN})\qty(\tfrac{1}{2}\abs{\vH}^2) + \vn \vdot \grad\qty[\tfrac{1}{2}\abs{\vH}^2 - \wt{\cH}\qty(\tfrac{1}{2}\abs{\vH}^2)_{\restriction_{\Gmt}}] - \vn\vdot\grad\qty(\tfrac{1}{2}\abs{\vH}^2).
	\end{equation}
	It follows from \eqref{pre maxwell}\textsubscript{2} that
	\begin{equation*}
		\grad(\tfrac{1}{2}\abs{\vH}^2) = \grad_\vH \vH,
	\end{equation*}
	which, together with \eqref{BC}\textsubscript{3}, yield that
	\begin{equation}
		-\vn\vdot\grad(\tfrac{1}{2}\abs{\vH}^2) = -\vn\vdot\grad_\vH\vH = \abs{\vH}^2\varkappa.
	\end{equation}
	Denote by
	\begin{equation}
		\wt{q} \coloneqq \tfrac{1}{2}\abs{\vH}^2 - \wt{\cH}(\tfrac{1}{2}\abs{\vH}^2)_{\restriction_{\Gmt}} \colon \cV_t \to \R.
	\end{equation}
	Then, it is clear that $\wt{q}$ uniquely solves
	\begin{equation}\label{def wtq}
		\begin{cases*}
			\lap\wt{q} = \abs{\grad\vH}^2 &in $\cV_t$, \\
			\wt{q} = 0 &on $\Gmt$, \\
			\grad_{\vN} \wt{q} = \vH \vdot \grad_\vN\vH &on $\cS$.
		\end{cases*}
	\end{equation}
	Therefore, one obtains that
	\begin{equation}\label{H2 bd term}
		\begin{split}
			\grad_\vtau\grad_\vtau\cN\qty(\tfrac{1}{2}\abs{\vH}^2) =\, &\abs{\vH}^2\grad_\vtau\grad_\vtau\varkappa + 2\grad_\vtau(\abs{\vH}^2) \grad_{\vtau}\varkappa + \varkappa\grad_\vtau\grad_\vtau(\abs{\vH}^2) + \grad_\vtau\grad_\vtau\grad_{\vn}\wt{q} + \\
			&+ \grad_\vtau\grad_\vtau(\cN-\wt{\cN})\qty(\tfrac{1}{2}\abs{\vH}^2).
		\end{split}
	\end{equation}
	Similar to \eqref{q1}-\eqref{def varrho}, one can conclude that, for $\wt{\varrho}$ defined by
	\begin{equation}\label{def wt varrho}
		\wt{\varrho} \coloneqq \lap\wt{q} - \grad^2 \wt{q}\qty(\vn_{\wt{\cH}}, \vn_{\wt{\cH}}) - \varkappa_{\wt{\cH}}\vn_{\wt{\cH}}\vdot\grad\wt{q},
	\end{equation}
	there holds
	\begin{equation}\label{wtq bd term}
		\grad_\vtau\grad_\vtau\grad_\vn\wt{q} = -(\grad_\vn \wt{q})\wt{\cN}\varkappa - (\grad_\vn\wt{q})\varkappa(\varkappa-\vn\vdot\wt{\cN}\vn) + (\grad_{\wt{\cN}\vn}\wt{q})\varkappa + 2\grad\vn_{\wt{\cH}}\vb{\colon}\grad^2\wt{q} + \grad_\vn\wt{\varrho}.
	\end{equation}
	
	Through combining \eqref{DtDt kappa 1}, \eqref{tan lap dn q}, \eqref{H2 bd term}, and \eqref{wtq bd term}, one arrives at
	\begin{equation}\label{DtDtk final}
		\begin{split}
			\Dt \Dt \varkappa =\, & \alpha\grad_\vtau\grad_\vtau\cN\varkappa -\alpha\abs{\grad_\vtau\varkappa}^2 + \alpha\varkappa^2\cN\varkappa + \abs{\vh}^2\grad_\vtau\grad_\vtau\varkappa + \abs{\vH}^2\grad_\vtau\grad_\vtau\varkappa + \\
			&+ \tfrac{5}{2}\grad_\vtau(\abs{\vh}^2)\grad_\vtau\varkappa + \tfrac{3}{2}\grad_\vtau(\abs{\vH}^2)\grad_\vtau\varkappa + \qty(\grad_\vn q - \grad_\vn \wt{q})\cN\varkappa + \\
			&+ \varkappa^3\abs{\vh}^2 + \tfrac{1}{2}\varkappa^2\cN(\abs{\vH}^2) + 2\varkappa\grad_\vtau\grad_\vtau(\abs{\vh}^2) + \varkappa\grad_\vtau\grad_\vtau\abs{\vH}^2 + (\grad_\vn\wt{q})(\cN - \wt{\cN})\varkappa  + \\
			&-(\grad_\vn q)\varkappa(\varkappa + \vn\vdot\cN\vn)   - (\grad_{\cN\vn}q)\varkappa- (\grad_\vn\wt{q})\varkappa(\varkappa-\vn\vdot\wt{\cN}\vn) + (\grad_{\wt{\cN}\vn}\wt{q})\varkappa + \\
			&+  2\grad\vn_\cH\vb{\colon}\grad^2q + \grad_\vn \varrho  + 2\grad\vn_{\wt{\cH}}\vb{\colon}\grad^2\wt{q} + \grad_\vn\wt{\varrho}+ \grad_\vtau\grad_\vtau(\cN-\wt{\cN})\qty(\tfrac{1}{2}\abs{\vH}^2) +  \\
			& + \fkR_0(\vv, \varkappa, \vtau, \vn) \\
			\eqqcolon\, & \alpha\grad_\vtau\grad_\vtau\cN\varkappa -\alpha\abs{\grad_\vtau\varkappa}^2 + \alpha\varkappa^2\cN\varkappa \quad \dasharrow \boxed{\text{ST terms}}  \\
			&+ \abs{\vh}^2\grad_\vtau\grad_\vtau\varkappa + \abs{\vH}^2\grad_\vtau\grad_\vtau\varkappa + \tfrac{5}{2}\grad_\vtau(\abs{\vh}^2)\grad_\vtau\varkappa + \tfrac{3}{2}\grad_\vtau(\abs{\vH}^2)\grad_\vtau\varkappa \quad \dasharrow \boxed{\text{MF terms}} \\
			&+ (\grad_\vn q - \grad_\vn\wt{q})\cN\varkappa \quad \dasharrow \boxed{\text{RT term}} \\
			& + \fkR_1(\vv, \vh, \vH, \Gmt). \quad \dasharrow \boxed{\text{Remainders}}
		\end{split}
	\end{equation}
	Here `ST', `MF', and `RT' stand for `surface tension', `magnetic fields', and `Rayleigh-Taylor sign' respectively. The naming rule follows directly from the coefficients of the derivatives of $\varkappa$. The factors in $\fkR_1$ can be viewed as remainders because, for geometric quantities, they only contain the first order (non-local) derivatives of $\vn$, which have the same regularity as $\varkappa$ does.
	
	\begin{remark}
		As we shall see in \S \ref{sec DN op}, the Dirichlet-Neumann operator $\cN$ is a self-adjoint, first order differential operator, which is, at the principal level, equivalent to $(-\slashed{\lap})^{\frac{1}{2}}$. (Here $\slashed{\lap}$ is the Laplace-Beltrami operator on $\Gmt$, and, for the present case, $\slashed{\lap} = \grad_\vtau\grad_\vtau$.) Thus, when counting the order of differential operators, one can infer straightforwardly from \eqref{DtDtk final} that the surface tension contributes to the highest order spacial derivatives, secondly does the tangential magnetic fields, and the Rayleigh-Taylor sign condition only contributes to a lower order term under the presence of surface tension or non-degenerate magnetic fields.
	\end{remark}
	
	\section{Technical Preparations}
	\label{sec tech preparation}
	
	\subsection{Reference Frames}\label{sec ref frame}
	As we aim to address the free boundary problems in Sobolev spaces, it would be more convenient that the implicit constants for elliptic estimates and the Sobolev embeddings are uniform under small variations of the moving domain. Inspired by \cite{Shatah-Zeng2008-Geo}, one may regard the fluid region $\Omt$ as the perturbation of a fixed one $\Oms$ with a smooth boundary $\Gms$. More precisely, for given $\Gms$, let $\vn_*$ be the outer unit normal of $\pd\Oms$. Then, there is a small constant $\delta > 0$, so that the map defined via:
	\begin{equation*}
		\begin{split}
			\vPhi \colon  & \Gms \times (-\delta, \delta) \to \R^2 \\
			& (z, \varphi) \mapsto z + \varphi \vn_*(z)
		\end{split}
	\end{equation*}
	is a diffeomorphism to a small neighborhood of $\Gms$. Namely, for each curve $\Gamma$ close to $\Gms$ in the $C^1$-topology, there exists a unique scalar coordinate function $\varphi_{\Gamma}$ defined on $\Gms$, so that the map:
	\begin{equation}\label{def vPhi}
		\begin{split}
			\vPhi_\Gamma \colon & \Gms \to \R^2 \\
			& z \mapsto z + \varphi_{\Gamma}(z)\vn_*(z)
		\end{split}
	\end{equation}
	is a $C^1$ diffeomorphism between $\Gms$ and $\Gamma$. 
	
	Since $\Gms$ is a fixed smooth curve, one can define the standard Sobolev spaces $H^s(\Gms)$ through a partition of unity (cf. \cite[\S 5.3]{ifrim2023sharp}). Such a partition still works for any $\Omt$ with $\Gmt$ sufficiently close to $\Gms$ in the $C^1$-topology. Thus, one can denote by
	\begin{equation}
		\abs{\Gamma}_{H^{s}} \coloneqq \abs{\varphi_{\Gamma}}_{H^{s}(\Gms)}.
	\end{equation}
	
	For the simplicity of notations, we shall use the following notion (cf. \cite{Shatah-Zeng2008-Geo}):
	\begin{defi}
		Let $s > 2$ and $0<\delta \ll 1$ be two constants. Define $\Lambda(\Gms, s-\frac{1}{2}, \delta)$ to be the collection of $C^1$ curves near $\Gms$, whose coordinate function $\varphi_{\Gamma}$ satisfies $\abs{\varphi_{\Gamma}}_{H^{s-\frac{1}{2}}(\Gms)}<\delta$.
	\end{defi}
	
	\begin{remark}
		By the Sobolev embedding $H^{1.5+}(\Gms) \hookrightarrow C^{1+}(\Gms)$, one can conclude that each curve in $\Lambda(\Gms, s-\frac{1}{2}, \delta) $ is close to $\Gms$ in the $C^1$-topology. In particular, so long as $\Gms$ does not self-intersect, the same property also holds for each $\Gamma \in \Lambda(\Gms, s-\frac{1}{2}, \delta)$ provided that $\delta \ll 1$. To avoid the clutter of notations, we shall denote $\Lambda(\Gms, s-\frac{1}{2}, \delta)$ simply by $\Lambda_s$.
	\end{remark}
	
	Indeed, any function defined on $\Gamma \in \Lambda_s$ can be pulled-back to the reference frame $\Gms$ through the coordinate map $\vPhi_\Gamma$. Moreover, there holds the following estimate (cf. \cite[Appendix A]{Shatah-Zeng2008-Geo}):
	\begin{equation}
		\abs{\psi\circ\vPhi_{\Gamma}}_{H^{\sigma}(\Gms)} \simeq_{\Lambda_s} \abs{\psi}_{H^{\sigma}(\Gamma)} \qfor \frac{1}{2}-s \le \sigma \le s-\frac{1}{2}.
	\end{equation}
	
	Concerning the regularity of boundaries, one has the following lemma, which can be derived from the definition of curvature and the adapted elliptic estimates (cf. \cite[Proposition 5.22]{ifrim2023sharp}):
	\begin{lemma}\label{lem reg Gam}
		Suppose that $\Gamma \in \Lambda_s$, and $\varkappa$ is the signed curvature of $\Gamma$. Then, it holds that
		\begin{equation*}
			\abs{\Gamma}_{H^{\sigma+2}} + \abs{\vtau}_{H^{\sigma+1}(\Gamma)} + \abs{\vn}_{H^{\sigma+1}(\Gamma)} \lesssim_{\Lambda_s} 1 + \abs{\varkappa}_{H^{\sigma}(\Gamma)} \qfor \sigma\ge 0.
		\end{equation*}
		Thus, $\varkappa \in H^{\sigma}(\Gamma)$ implies that $\Gamma \in H^{\sigma+2}$, for any $\Gamma \in \Lambda_s$ and $\sigma \ge 0$.
	\end{lemma}
	
	Now, we turn to the correspondence between interiors. More precisely, given a bounded domain $\Omega$ with $\pd\Omega \in \Lambda_s$, one can also define a bijection between $\Omega$ and $\Oms$ via the harmonic extension:
	\begin{equation}
		\begin{cases*}
			\lap \cX_{\Gamma} = 0 &in $\Oms$,\\
			\cX_{\Gamma} = \vPhi_{\Gamma} &on $\pd\Oms$.
		\end{cases*}
	\end{equation}
	Then, it is clear that
	\begin{equation}
		\norm{\grad\cX_{\Gamma} - \textup{Id}_{\restriction_{\Oms}}}_{H^{s-1}(\Oms)} \lesssim_{\Oms} \abs{\phi_{\Gamma}}_{H^{s-\frac{1}{2}}(\Gms)} \lesssim_{\Oms} \delta.
	\end{equation}
	In particular, there is a generic constant $\delta_0 \ll 1$ depending only on the geometry of $\Oms$, so that, whenever $\delta \le \delta_0 \ll 1$, the map $\cX_{\Gamma} \colon \Oms \to \Omega$ is a $C^1$-diffeomorphism. From now on, we will always implicitly assume that the bound $\delta$ in $\Lambda_s$ is necessarily small. 
	
	Similarly, one also has the following equivalence of norms:
	\begin{equation}
		\norm{f \circ \cX_{\Gamma}}_{H^{\sigma}(\Oms)} \simeq_{\Lambda_s} \norm{f}_{H^{\sigma}(\Omega)} \qfor -s \le \sigma \le s.
	\end{equation}
		
	Throughout this manuscript, we will always assume that $\Omega$ is a bounded domain, and its boundary $\pd\Omega \in \Lambda_s = \Lambda(\Gms, s-\frac{1}{2}, \delta)$ with $s > 2$ and $\delta \ll 1$. We shall not repeat this hypothesis later.

	\subsection{Product Estimates}
		
	Suppose that $\sigma_1 \le \sigma_2 \le s-\frac{1}{2}$, then there hold  (cf. \cite[Appendix A]{Shatah-Zeng2008-Geo}):
	\begin{equation}
			\abs{\phi\psi}_{H^{\sigma_1 + \sigma_2 - \frac{1}{2}}(\Gamma)} \lesssim_{\Lambda_s} \abs{\phi}_{H^{\sigma_1}(\Gamma)} \abs{\psi}_{H^{\sigma_2}(\Gamma)} \qif \sigma_2 < \frac{1}{2} \text{ and }  \sigma_1 + \sigma_2 > 0;
	\end{equation}
		and
	\begin{equation}
			\abs{\phi\psi}_{H^{\sigma_1}(\Gamma)} \lesssim_{\Lambda_s} \abs{\phi}_{H^{\sigma_1}(\Gamma)} \abs{\psi}_{H^{\sigma_2}(\Gamma)} \qif \sigma_2  > \frac{1}{2} \text{ and } \sigma_1 + \sigma_2 \ge 0.
	\end{equation}
	Furthermore, one has the following refinements for any $\sigma_3 \ge 0$ (cf. \cite[\S 5.3.2]{ifrim2023sharp}):
	\begin{equation}
		\begin{split}
			\abs{\phi\psi}_{H^{\sigma_3}(\Gamma)} \lesssim_{\Lambda_s}\, & \qty(\abs{\phi}_{L^\infty(\Gamma)}\abs{\psi}_{W^{1, \infty}(\Gamma)} + \abs{\phi}_{W^{1, \infty}(\Gamma)}\abs{\psi}_{L^\infty(\Gamma)})\cdot \abs{\Gamma}_{H^{\sigma_3}} + \\
			&+ \abs{\phi}_{L^\infty(\Gamma)}\abs{\psi}_{H^{\sigma_3}(\Gamma)} + \abs{\phi}_{H^{\sigma_3}(\Gamma)}\abs{\psi}_{L^{\infty}(\Gamma)},
		\end{split}
	\end{equation}
	and
	\begin{equation}
		\abs{\phi\psi}_{H^{\sigma_3}(\Gamma)} \lesssim_{\Lambda_s} \abs{\phi}_{L^\infty(\Gamma)}\abs{\psi}_{L^\infty(\Gamma)}\abs{\Gamma}_{H^{1+\sigma_3}} + \abs{\phi}_{L^\infty(\Gamma)}\abs{\psi}_{H^{\sigma_3}(\Gamma)} + \abs{\phi}_{H^{\sigma_3}(\Gamma)}\abs{\psi}_{L^\infty(\Gamma)}.
	\end{equation}

	For the estimates in the interior, there hold for $s_1 \le s_2$:
	\begin{equation}
		\norm{fg}_{H^{s_1 + s_2 -1}(\Omega)} \lesssim_{\Lambda_s} \norm{f}_{H^{s_1}(\Omega)} \norm{g}_{H^{s_2}(\Omega)} \qif s_2 < 1 \text{ and } s_1 + s_2 > 0;
	\end{equation}
	and
	\begin{equation}
		\norm{fg}_{H^{s_1}(\Omega)} \lesssim_{\Lambda_s} \norm{f}_{H^{s_1}(\Omega)} \norm{g}_{H^{s_2}(\Omega)} \qif s_2 > 1 \text{ and } s_1 + s_2 \ge 0,
	\end{equation}
	which can be achieved through the combination of Stein's extension operators and the standard para-product decompositions (see also \cite[\S\S 5.1-5.2]{ifrim2023sharp}).
	
	\subsection{Uniform Elliptic Estimates}
	After the introduction of reference frames, one can unify the implicit constants in the Sobolev embeddings and elliptic estimates. We first state the following result (cf. \cite[Proposition 5.19]{ifrim2023sharp}):
	\begin{lemma}
		Let $\sigma \ge 2$ be a constant and $u$ a solution to the Dirichlet problem:
		\begin{equation*}
			\begin{cases*}
				\lap u = f &in $\Omega$,\\
				u = g &on $\pd\Omega$.
			\end{cases*}
		\end{equation*}
		Then, for parameters $\lambda \ge 0$, $0 \le \beta, \gamma \le 1$, $0 < \epsilon < s-2$, and any sequence of partitions $u = u_j^1 + u_j^2$, there holds
		\begin{equation*}
			\begin{split}
			 \norm{u}_{H^{\sigma}(\Omega)} \lesssim_{\Lambda_s, \epsilon}\, &\norm{f}_{H^{\sigma-2}(\Omega)} + \abs{g}_{H^{\sigma-\frac{1}{2}}(\pd\Omega)} + \abs{\pd\Omega}_{H^{\sigma+\lambda-\frac{1}{2}}} \sup_{j > 0} 2^{-j(\lambda+\beta-1)}\norm*{u_j^1}_{C^{\beta}(\Omega)} + \\
			 &+\sup_{j > 0} 2^{j(\sigma+\gamma-1-\epsilon)}\norm*{u_j^2}_{H^{1-\gamma}(\Omega)}.
			\end{split}
		\end{equation*}
	\end{lemma}
	In particular, by taking $u_j^1 \coloneqq P_{\le j} \cE_{\Omega} u$, $u_j^2 \coloneqq P_{>j}\cE_{\Omega} u$ (here $P_j$ is the standard Littlewood-Paley projector and $\cE_\Omega$ is Stein's extension operator), and $\lambda = 0$, standard interpolations yield the following
	\begin{equation}
			\norm{u}_{H^\sigma(\Omega)} \lesssim_{\Lambda_s} \norm{\lap u}_{H^{\sigma-2}(\Omega)} + \abs{u}_{H^{\sigma-\frac{1}{2}}(\pd\Omega)} + \cP{\abs{\pd\Omega}_{H^{\sigma-\frac{1}{2}}}}\cdot \norm{u}_{L^2(\Omega)},
	\end{equation}
	where $\mathcal{P}$ is a generic polynomial.
	
	Moreover, one has the bound (cf. \cite[\S A.2]{Shatah-Zeng2008-Geo} and \cite[\S 5.5.1]{ifrim2023sharp}):
	\begin{equation}
		\norm{\lap^{-1} f}_{H^{\sigma_1+2}(\Omega)} \lesssim_{\Lambda_s} \norm{f}_{H^{\sigma_1}(\Omega)} \qfor -1 \le \sigma_1 \le s- 2,
	\end{equation}
	here and henceforward $\lap^{-1}\psi \eqqcolon \phi$ represents the solution to Poisson's equation with a homogeneous Dirichlet boundary condition:
	\begin{equation*}
		\begin{cases*}
			\lap \phi = \psi &in $\Om$, \\
			\phi = 0 &on $\pd\Omega$.
		\end{cases*}
	\end{equation*}
	For the harmonic extensions, one has
	\begin{equation}
		\norm{\cH g}_{H^{\sigma_2 + \frac{1}{2}}(\Omega)} + \norm{\wt{\cH} g}_{H^{\sigma_2 + \frac{1}{2}}(\cV)} \lesssim_{\Lambda_s} \abs{g}_{H^{\sigma_2}(\Gamma)} \qfor 0 < \sigma_2 \le s - \frac{1}{2}.
	\end{equation}
	
	At the end of this subsection, we sate the trace estimates (cf. \cite[Propositions 5.11 \& 5.28]{ifrim2023sharp}):
	\begin{lemma}
		Let $\sigma > 0$ be a constant and $f$ a function defined in $\Omega$. Then, there hold the following estimates:
		\begin{equation*}
			\abs{f_{\restriction_{\pd\Omega}}}_{H^{\sigma}(\pd\Omega)} \lesssim_{\Lambda_s} \qty\big(1 + \abs{\pd\Omega}_{H^{\sigma}}) \cdot \norm{f}_{H^{\sigma+\frac{1}{2}}(\Omega)},
		\end{equation*}
		and
		\begin{equation*}
			\abs{\grad_\vn f}_{H^{\sigma}(\pd\Omega)} \lesssim_{\Lambda_s} \norm{f}_{H^{\sigma+\frac{3}{2}}(\Omega)} \times \begin{cases*}
				\qty\big(1 + \abs{\pd\Omega}_{H^{\sigma+\frac{3}{2}}}) &if $\sigma < \frac{1}{2}$, \\
				\qty\big(1 + \abs{\pd\Omega}_{H^{\sigma+1}}) &if $\sigma \ge \frac{1}{2}$.
			\end{cases*}
		\end{equation*}
	\end{lemma}
	
	\subsection{Div-Curl Systems}
	First, we state the following div-curl estimate with Neumann-type boundary data (cf. \cite[Proposition 5.27]{ifrim2023sharp}):
	\begin{lemma}\label{lem div-curl}
		Let $\vw \in H^{\sigma}(\Omega)$ be a vector field defined in $\Omega$ and $\sigma > \frac{3}{2}$. Then, $\vw$ satisfies the following estimate:
		\begin{equation*}
			\begin{split}
				\norm{\vw}_{H^{\sigma}(\Omega)} \lesssim_{\Lambda_s}\, &\norm{\dive\vw}_{H^{\sigma-1}(\Omega)} + \norm{\Curl \vw}_{H^{\sigma-1}(\Omega)} + \abs{\vn\vdot\grad_\vtau\vw}_{H^{\sigma-\frac{3}{2}}(\pd\Omega)} + \\
				&+ \begin{cases*}
					\cP{\abs{\pd\Omega}_{H^{\sigma}}} \cdot \norm{\vw}_{L^2(\Omega)} &if $ \sigma < 2$, \\
					\cP{\abs{\pd\Omega}_{H^{\sigma-\frac{1}{2}}}} \cdot \norm{\vw}_{L^2(\Omega)} &if $\sigma \ge 2$,
				\end{cases*} 
			\end{split}
		\end{equation*}
		where $\mathcal{P}$ is a generic polynomial.
	\end{lemma}
	
	Concerning the solvability of div-curl problems, since the plasma region $\Om$ is always assumed to be simply connected, the following div-curl problem admits a unique solution (cf. \cite[\S7]{auchmuty2001})
	\begin{equation}
		\begin{cases*}
			\dive\vw = f &in $\Omega$, \\
			\Curl \vw = g &in $\Omega$, \\
			\vw\vdot\vn = \theta &on $\Gamma$,
		\end{cases*}
	\end{equation}
	for which there holds the compatibility condition
	\begin{equation}
		\int_{\Omega} f \dd{x} = \int_{\Gamma} \theta \dd{\ell}.
	\end{equation}
	Furthermore, there holds the estimate (see also \cite[\S4]{Sun-Wang-Zhang2019}):
	\begin{equation}
		\norm{\vw}_{H^{s'}(\Omega)} \le \cP{\abs{\Gamma}_{H^{s'+\frac{1}{2}}}} \cdot \qty(\norm{f}_{H^{s'-1}(\Omega)} + \norm{g}_{H^{s'-1}(\Omega)} + \abs{\theta}_{H^{s'-\frac{1}{2}}(\Gamma)}) \qfor s'\ge 1,
	\end{equation}
	where $\mathcal{P}$ is a generic polynomial determined by $\Lambda_s$.
	
	Although the vacuum region $\cV$ is not simply-connected, the following div-curl problem with mixed type boundary data is still \emph{uniquely} solvable (c.f \cite[\S\S14-15]{auchmuty2001}):
	\begin{equation}\label{div-curl'}
		\begin{cases*}
			\dive\vH = 0 &in $\cV$, \\
			\Curl \vH = 0&in $\cV$, \\
			\vH\vdot\vn = 0 &on $\Gamma$, \\
			\vN \smallwedge \vH = \cJ &on $\cS$.
		\end{cases*}
	\end{equation}
	Moreover, the solution $\vH$ satisfies the estimate (see also \cite[\S4]{Sun-Wang-Zhang2019}):
	\begin{equation}
		\norm{\vH}_{H^{\sigma}(\cV)} \le \cP{\abs{\Gamma}_{H^{\sigma+\frac{1}{2}}}} \cdot \abs{\cJ}_{H^{\sigma-\frac{1}{2}}(\cS)} \qfor \sigma \ge 1,
	\end{equation}
	where $\mathcal{P}$ is a generic polynomial determined by $\Lambda_s$ and the fixed boundary $\cS$.
	
	\begin{remark}
		If the Betti number of $\cV$ exceeds one, the problem \eqref{div-curl'} may exist non-trivial solution even when $\cJ = 0$. Thus, it is always assumed in this manuscript that $\cV$ has only one ``hole''.
	\end{remark}
	
	\subsection{Dirichlet-Neumann Operators}\label{sec DN op}
	The Dirichlet-to-Neumann operator on $\Gamma$ is defined as:
	\begin{equation}
		\cN f \coloneqq \vn_{\Gamma} \vdot (\grad\cH f)_{\restriction_{\Gamma}}.
	\end{equation}
	Assume that $ \Gamma \in \Lambda_s \subset H^{s-\frac{1}{2}} $ and $ \frac{3}{2}-s \le \sigma \le s-\frac{1}{2} $. Then, the Dirichlet-Neumann operator 
	\begin{equation*}
		\cN : H^{\sigma}(\Gamma) \to H^{\sigma-1}(\Gamma)
	\end{equation*} 
	satisfies the following properties (cf. \cite[\S A.2]{Shatah-Zeng2008-Geo}):
	\begin{enumerate}
		\item $ \cN $ is self-adjoint and non-negative on $ L^2(\Gamma) $ with compact resolvents.
		\item $\ker(\cN) = \qty{\const}$.
		\item It holds for any $ f $ satisfying $ \int_\Gamma f \dd{\ell} = 0 $ that
		\begin{equation}\label{property1 n}
			\abs{f}_{H^{\sigma}(\Gamma)} \simeq_{\Lambda_s} \abs{\cN f }_{H^{\sigma-1}(\Gamma)}.
		\end{equation} 
		\item For $ \frac{1}{2}-s \le s_1 \le s-\frac{1}{2} $, there holds
		\begin{equation}\label{equiv n lap}
			\qty(\textup{I}-\slashed{\lap})^{\frac{s_1}{2}} \simeq_{\Lambda_s} \qty(\textup{I}+\cN)^{s_1},
		\end{equation}
		i.e., the norms on $ H^{s_1}(\Gamma) $ defined by interpolating $ \qty(\textup{I}-\slashed{\lap})^\frac{1}{2} $ and $ \qty(\textup{I}+\cN) $ are equivalent.
		\item For $ \frac{1}{2}-s \le s_2 \le s-\frac{3}{2} $,
		\begin{equation*}
			(\cN)^{-1} : H^{s_2}_{0}(\Gamma) \to H^{s_2+1}_{0}(\Gamma) \qq{with}
			H^{s_2}_{0}(\Gamma) \coloneqq \Set*{f\in H^{s_2}(\Gamma) \given \int_\Gamma f \dd{S} = 0 }
		\end{equation*}
		is well-defined and bounded uniformly in $ \Gamma \in \Lambda_s $.
	\end{enumerate}
	It is clear that all the above properties also hold for $\wt{\cN}$.
	
	In addition to \eqref{property1 n}-\eqref{equiv n lap}, one has the following quantitative characterization of Sobolev norms on $\Gamma \in \Lambda_s$ (cf. \cite[\S 5.5.4]{ifrim2023sharp}):
	\begin{lemma}\label{lem est N}
		Let $s' \ge \frac{1}{2}$ be a real number and $k \ge 1$ an integer. Then, for any $\Gamma \in \Lambda_s$, there hold the following bounds:
		\begin{equation*}
			\abs{f}_{H^{s'+k}(\Gamma)} \lesssim_{\Lambda_s} \abs{\cN^k f}_{H^{s'}(\Gamma)} + \cP{\abs{\Gamma}_{H^{s'+k}}} \cdot \abs{f}_{L^2(\Gamma)},
		\end{equation*}
		and
		\begin{equation*}
			\abs{\cN^k f}_{H^{s'}(\Gamma)} \lesssim_{\Lambda_s} \qty(1 + \abs{\Gamma}_{H^{s'+k}}) \cdot \abs{f}_{H^{s'+k}(\Gamma)},
		\end{equation*}
		where $\mathcal{P}$ represents a generic polynomial.
	\end{lemma}

	To further indicate that $\cN$ behaves like a derivative, we state the following Leibniz rule (cf. \cite[eq. (A.11)]{Shatah-Zeng2008-Geo}):
	\begin{equation}\label{Leibniz N}
		\cN(fg) = f \cN g + g \cN f - 2\grad_\vn\lap^{-1}\qty(\grad f_{\cH} \vdot \grad g_{\cH}).
	\end{equation}
	In particular, the product and elliptic estimates yield the following commutator estimates:
	\begin{equation}
		\abs{\comm{\cN}{f}}_{\scL \colon H^{\sigma_1}(\Gamma) \to H^{\sigma_1} (\Gamma)} \lesssim_{\Lambda_s} \abs{f}_{H^{s-\frac{1}{2}}(\Gamma)} \qfor -\frac{1}{2} \le \sigma_1 \le \frac{1}{2},
	\end{equation}
	and
	\begin{equation}
		\abs{\comm{\cN}{f}}_{\scL \colon H^{\sigma_2}(\Gamma) \to H^{\sigma_2} (\Gamma)} \lesssim_{\Lambda_s} \cP{\abs{\Gamma}_{H^{1+\sigma_2}}} \cdot \abs{f}_{H^{1+\sigma_2}(\Gamma)} \qfor\sigma_2\ge\frac{1}{2},
	\end{equation}
	where $\mathcal{P}$ is a generic polynomial.
	
	Next, in view of the relation \eqref{lap rel}, one has (cf. \cite[eq. (A.13)]{Shatah-Zeng2008-Geo})
	\begin{equation*}
		(\grad_\vtau\grad_\vtau + \cN^2)f =( \vn\vdot\cN \vn -\varkappa)\cN f - 2\grad_\vn \lap^{-1}\qty(\grad\vn_{\cH} \vb{\colon} \grad^2f_{\cH}) + (\grad_\vtau f) \vtau\vdot\cN \vn ,
	\end{equation*}
	Thus, product and elliptic estimates, together with Lemma \ref{lem reg Gam}, yield that (cf. \cite[Proposition A.6]{Shatah-Zeng2008-Geo}):
	\begin{lemma}\label{lem lap + N2}
		Assume that $\Gamma \in \Lambda_s$ and $\Gamma \in H^{\sigma}$ for some $\sigma > \frac{5}{2}$. Then, there holds the following estimate:
		\begin{equation*}
			\abs{\grad_\vtau\grad_\vtau + \cN^2}_{\scL\colon H^{\sigma'}(\Gamma) \to H^{\sigma'-1}(\Gamma)} \lesssim_{\Lambda_s} 1 + \abs{\varkappa}_{H^{\sigma-2}(\Gamma)} \qfor 2-\sigma \le \sigma' \le \sigma-1.
		\end{equation*}
	\end{lemma}
	To estimate the operator $\wt{\cN}$, one can first extend the normal field $\vn$ into $\cV$ through solving the elliptic problem
	\begin{equation*}
		\begin{cases*}
			\lap\wt{\vnu} = \vb{0} &in $\cV$,\\
			\wt{\vnu} = \vn &on $\Gamma$,\\
			\wt{\vnu} = \vb{0} &on $\cS$.
		\end{cases*}
	\end{equation*}
	Additionally, by utilizing the Bogovskii map (cf. \cite[\S2.8]{Tsai2018book} and the references therein), one can take a vector field $\wh{\vnu}$, so that
	\begin{equation*}
		\begin{cases*}
			\dive\wh{\vnu} = \dive\wt{\vnu} - \frac{1}{\abs{\cV}}\int_{\cV}(\dive\wt{\vnu}) \dd{x} &in $\cV$,\\
			\wh{\vnu} = \vb{0} &on $\Gamma$, \\
			\wh{\vnu} = \vb{0} &on $\cS$,
		\end{cases*}
	\end{equation*}
	where $\abs{\cV}$ represents the volume (area) of $\cV$. Thus, by defining the vector field
	\begin{equation*}
		\vnu \coloneqq \wt{\vnu}-\wh{\vnu} \colon \cV \to \R^2,
	\end{equation*}
	one obtains that
	\begin{equation*}
		\begin{cases*}
			\dive\vnu = \const &in $\cV$, \\
			\vnu = \vn &on $\Gamma$, \\
			\vnu = \vb{0} &on $\cS$.
		\end{cases*}
	\end{equation*}
	Moreover, there holds the following estimate
	\begin{equation*}
		\norm{\vnu}_{H^{s'}(\cV)} \lesssim_{\Lambda_s} 1 + \abs{\varkappa}_{H^{s'-\frac{3}{2}}(\Gamma)}.
	\end{equation*}
	Hence, one can now calculate that
	\begin{equation*}
		\begin{split}
			\grad^2 f_{\wt{\cH}}(\vn, \vn) &= \grad_\vn\grad_\vnu f_{\wt{\cH}} - \grad_\vn\vnu\vdot\grad f_{\wt{\cH}}\\
			&= \grad_\vn\wt{\cH}(\grad_\vnu f_{\wt{\cH}})_{\restriction_{\Gamma}} + \grad_\vn\qty[\grad_\vnu f_{\wt{\cH}} - \wt{\cH}(\grad_\vnu f_{\wt{\cH}})_{\restriction_\Gamma}] - \grad_\vn\vnu\vdot\grad f_{\wt{\cH}} \\
			&= \wt{\cN}^2 f + \grad_\vn \uppsi - \grad_\vn\vnu\vdot\grad f_{\wt{\cH}},
		\end{split}
	\end{equation*}
	where $\uppsi$ satisfies the following elliptic system:
	\begin{equation*}
		\begin{cases*}
			\lap\uppsi = \lap(\vnu\vdot\grad f_{\wt{\cH}}) &in $\cV$, \\
			\uppsi = 0 &on $\Gamma$, \\
			\grad_\vN\uppsi = \grad_\vN(\vnu\vdot\grad f_{\wt{\cH}}) &on $\cS$.
		\end{cases*}
	\end{equation*}
	It follows from the construction of $\vnu$ that
	\begin{equation*}
		\begin{split}
			\lap(\vnu\vdot\grad f_{\wt{\cH}}) &= 2\grad\vnu\vb{\colon}\grad^2 f_{\wt{\cH}} + \grad f_{\wt{\cH}} \vdot\lap\vnu \\
			&= \dive(\grad f_{\wt{\cH}}\vdot\grad\vnu + \grad\vnu\vdot\grad f_{\wt{\cH}}).
		\end{split}
	\end{equation*}
	For the simplicity of notations, we denote by
	\begin{equation*}
		\vphi \coloneqq \grad f_{\wt{\cH}}\vdot\grad\vnu + \grad\vnu\vdot\grad f_{\wt{\cH}}.
	\end{equation*}
	Thus, it is clear that
	\begin{equation*}
		\grad_\vN\uppsi = \vphi \vdot \vN \qq{on} \cS.
	\end{equation*}
	Consider the following Neumann-type problem:
	\begin{equation*}
		\begin{cases*}
			\lap\uppsi_1 = \dive\vphi &in $\cV$, \\
			\grad_\vn \uppsi_1 = \vphi\vdot\vn &on $\Gamma$, \\
			\grad_\vN \uppsi_1 = \vphi\vdot\vN &on $\cS$, \\
			\int_{\cV} \uppsi_1 \dd{x} = 0,
		\end{cases*}
	\end{equation*}
	which is uniquely solvable with estimate (cf. \cite[Theorem 1.2]{Geng2012}):
	\begin{equation*}
		\norm{\uppsi_1}_{H^1(\cV)} \lesssim_{\Lambda_s} \norm{\vphi}_{L^2(\cV)}.
	\end{equation*}
	It is then obvious that $\uppsi$ is given by
	\begin{equation*}
		\uppsi = \uppsi_1 - \wt{\cH}(\uppsi_1)_{\restriction_{\Gamma}}.
	\end{equation*}
	Next, one can observe that
	\begin{equation*}
		\grad_\vn\uppsi - \grad_\vn\vnu\vdot\grad f_{\wt{\cH}} = \wt{\cN}(\uppsi_1) + \vn\vdot\grad\vnu\vdot\grad f_{\wt{\cH}} = \wt{\cN}(\uppsi_1) -\grad\vnu(\vn, \vn)\wt{\cN}(f).
	\end{equation*}
	It follows from \eqref{lap rel} that
	\begin{equation*}
		\qty\big(\grad_\vtau\grad_\vtau + \wt{\cN}^2)f = \qty\big[\varkappa + \grad\vnu(\vn, \vn)]\wt{\cN}(f) - \wt{\cN}(\uppsi_1),
	\end{equation*}
	which implies the estimate
	\begin{equation*}
		\abs{\qty\big(\grad_\vtau\grad_\vtau + \wt{\cN}^2)f}_{H^{-\frac{1}{2}}(\Gamma)} \lesssim_{\Lambda_s} \qty(1+\abs{\varkappa}_{H^{\frac{1}{2}+}(\Gamma)})\cdot \abs{f}_{H^{\frac{1}{2}}(\Gamma)}.
	\end{equation*}
	Standard interpolations yield that, under the same hypothesis in Lemma \ref{lem lap + N2}, there holds
	\begin{equation}
		\abs{\grad_\vtau\grad_\vtau + \wt{\cN}^2}_{\scL\colon H^{\sigma''}(\Gamma) \to H^{\sigma''-1}(\Gamma)} \lesssim_{\Lambda_s} 1 + \abs{\varkappa}_{H^{\sigma - 2}(\Gamma)} \qfor \frac{1}{2} \le \sigma'' \le \sigma - 1.
	\end{equation}

	Furthermore, a functional analysis result \cite[Proposition A.7]{Shatah-Zeng2008-Geo} yields the following (cf. \cite[Theorem A.8]{Shatah-Zeng2008-Geo} and \cite[Lemma 4.6]{Shatah-Zeng2008-vortex}):
	\begin{lemma}\label{lem N lap}
		Suppose that $\Gamma \in \Lambda_s \cap H^{\sigma}$ with $\sigma > \frac{5}{2}$. Then, there hold
		\begin{equation*}
			\abs{(-\slashed{\lap})^{\frac{1}{2}} - \cN}_{\scL\colon H^{\sigma_1}(\Gamma) \to H^{\sigma_1}(\Gamma)} \lesssim_{\Lambda_s}  1 + \abs{\varkappa}_{H^{\sigma-2}} \qfor 1-\sigma \le \sigma_1 \le \sigma -1
		\end{equation*}
		and
		\begin{equation*}
				\abs{(-\slashed{\lap})^{\frac{1}{2}} - \wt{\cN}}_{\scL\colon H^{\sigma_2}(\Gamma) \to H^{\sigma_2}(\Gamma)} \lesssim_{\Lambda_s}  1 + \abs{\varkappa}_{H^{\sigma-2}} \qfor -\frac{1}{2} \le \sigma_1 \le \sigma -1.
		\end{equation*}
	\end{lemma}
	Particularly, one can regard both $\cN$ and $\wt{\cN}$ as $(-\slashed{\lap})^{\frac{1}{2}}$ modulo an operator of order zero.

	\section{Linear Propagation Estimates}\label{sec lin est}
	
	\subsection{Evolution of the Curvature}\label{sec evo lin k}
	In view of the evolution equation \eqref{DtDtk final}, one can further simplify it as (note that the magnetic field $\vh$ is assumed to be tangential to $\Gmt$):
	\begin{equation}\label{DtDt k new}
		\begin{split}
			\Dt\Dt\varkappa =\, &\alpha\grad_\vtau\grad_\vtau\cN\varkappa -\alpha\abs{\grad_\vtau\varkappa}^2 + \alpha\varkappa^2\cN\varkappa + \abs{\vh}^2\grad_\vtau\grad_\vtau\varkappa + \abs{\vH}^2\grad_\vtau\grad_\vtau\varkappa +\\
			&+   \tfrac{5}{2}\grad_\vtau(\abs{\vh}^2)\grad_\vtau\varkappa + \tfrac{3}{2}\grad_\vtau(\abs{\vH}^2)\grad_\vtau\varkappa + (\grad_\vn q-\grad_\vn\wt{q})\cN\varkappa  + \fkR_1(\vv, \vh, \vH, \Gmt) \\
			=\, & \underbrace{\alpha\grad_\vtau\grad_\vtau\cN\varkappa +\grad_\vh\grad_\vh\varkappa + \grad_\vH\grad_\vH\varkappa - (\grad_\vn\wt{q}-\grad_\vn q)\cN\varkappa}_{\boxed{\text{main terms}}} \, + \, \underbrace{\fkR_2(\vv, \vh, \vH, \Gmt)}_{\boxed{\text{remainder}}}.
		\end{split}
	\end{equation}
	We shall see later in \S\ref{sec a pri est} that the new remainder $\fkR_2(\vv, \vh, \vH, \Gmt)$ is controllable through the energies induced by the main terms.
	
	\subsubsection{Simplified linear equation}
	In order to clarify the boundary evolutions, we now assume that $\{\Gmt\}_{t} \subset \Lambda_s \cap H^{s'}$ (here $s'$ could be as large as necessary) is a family of closed curves evolving with velocity $\vbu \colon \Gmt \to \R^2$. Suppose that $\mu$, $\nu$, and $f$ are scalar functions defined on $\Gmt$. Denote the trajectory derivative along $\vbu$ by
	\begin{equation*}
		\cDt \coloneqq \pdv{t} + \grad_\vbu.
	\end{equation*}
	Consider the following linear evolution equation on $\Gmt$, which can be viewed as a simplified linearized equation for \eqref{DtDt k new}:
	\begin{equation}\label{lin eqn kp}
		\cDt\cDt \vartheta - \alpha \grad_\vtau\grad_\vtau \cN\vartheta - \grad_{\mu\vtau} \grad_{\mu\vtau}\vartheta + \nu\cN\vartheta = f,
	\end{equation}
	where $\vartheta(t)$ is interpreted as a scalar function defined on $\Gmt$, and $\cN$ is the Dirichlet-Neumann operator on $\Gmt$. 
	
	The main goal of \S\ref{sec evo lin k} is to derive the following:
	\begin{prop}\label{prop lin est theta}
		For each integer $m \ge 0$, define the energy functional:
		\begin{equation*}
			\Elin^m(\vartheta) \coloneqq \int_{\Gmt} \abs*{(\slashed{\lap}^m\cN)^{\frac{1}{2}}\cDt\vartheta}^2 + \alpha\abs*{\slashed{\grad}^{1+m}\cN\vartheta}^2 + \abs*{(\slashed{\lap}^m\cN)^{\frac{1}{2}}(\mu\grad_\vtau\vartheta)}^2 + \nu\abs*{\slashed{\grad}^m\cN\vartheta}^2 \dd{\ell}
		\end{equation*}
		and its corresponding bound
		\begin{equation*}
			\Mlin^m(\vartheta) \coloneqq \abs{\cDt\vartheta}_{H^{m+\frac{1}{2}}(\Gmt)}^2 + \alpha\abs{\vartheta}_{H^{m+2}(\Gmt)}^2 + \abs{\mu\grad_\vtau\vartheta}_{H^{m+\frac{1}{2}}(\Gmt)}^2 + \abs{\vartheta}_{H^{m+1}(\Gmt)}^2,
		\end{equation*}
		where $\slashed{\grad} \equiv \grad_\vtau$ is the tangential gradient, and the fractional power of $(\slashed{\lap}^m\cN)$ is defined through functional calculus on the Hilbert space $L^2(\Gmt)$. 
		
		Then, each solution $\vartheta$ to the linear evolution equation \eqref{lin eqn kp} satisfies the propagation estimate
		\begin{equation*}
			\begin{split}
				\abs{\dv{t}\Elin^m(\vartheta)} \le\, & \cP{\abs{\varkappa}_{H^{m+\frac{1}{2}}}, \alpha\abs{\varkappa}_{H^{m+1}}^2, \abs*{(\vbu, \mu, \nu)}_{H^{m+2}}, \abs{\cDt\nu}_{L^\infty}} \cdot \Mlin^m(\vartheta) + \\
				&\quad+ \cP{\abs{\Gmt}_{H^{m+1}}} \cdot \qty(\abs*{\cL_t\mu \cdot \grad_\vtau\vartheta}_{H^{m + \frac{1}{2}}(\Gmt)}^2 + \abs{f}_{H^{m+\frac{1}{2}}(\Gmt)}^2),
			\end{split}
		\end{equation*}
		where $\cL_t\mu \coloneqq \cDt\mu-(\grad_\vtau\vbu\vdot\vtau)\mu$ and $\mathcal{P}$ represents generic polynomials determined by $\Lambda_s$ and $m$.
	\end{prop}
	
	\subsubsection{Preliminary estimates}\label{sec pre lin est}
	To show the energy estimates, one may first note that the transport formula on $\Gmt$ is expressed as:
	\begin{equation}
		\dv{t}\int_{\Gmt} \psi \dd{\ell} = \int_{\Gmt} \cDt\psi + \qty(\vtau\vdot\grad_\vtau\vbu)\psi \dd{\ell},
	\end{equation}
	which directly follows from the relation \eqref{Dt dl}. Since $\Gmt$ is a closed curve, the integration by parts formula can be written as:
	\begin{equation}
		\int_{\Gmt} \psi\grad_\vtau\phi + \phi\grad_\vtau\psi \dd{\ell} = \int_{\Gmt} \grad_\vtau(\phi\psi) \dd{\ell} = 0.
	\end{equation}
	For the commutator formulas on $\Gmt$, one can first recall from \eqref{comm Dt d slash} that
	\begin{equation}\label{comm Dt dtau}
		\comm{\cDt}{\grad_\vtau} = -(\grad_\vtau\vbu \vdot\vtau)\grad_\vtau.
	\end{equation}
	When commuting $\cDt$ with $\cN$, one has the following relation, the derivation of which can be found in \cite[\S 3.1]{Shatah-Zeng2008-Geo}:
	\begin{equation}\label{comm formula dt N}
		\comm{\cDt}{\cN}\psi = 2\grad_\vn\lap^{-1}\dive\qty(\grad\psi_{\cH}\vdot\grad\vbu_{\cH}) - \grad\psi_{\cH}\vdot\cN\vbu - (\grad_\vtau\psi)\vn\vdot\grad_\vtau\vbu.
	\end{equation}
	Thus, there holds following estimate (cf. \cite[eq. (4.24)]{Shatah-Zeng2008-Geo}):
	\begin{equation*}
		\abs{\comm{\cDt}{\cN}}_{\scL\colon H^{\sigma}(\Gmt) \to H^{\sigma-1}(\Gmt)} \lesssim_{\Lambda_s} \abs{\vbu}_{H^{s-\frac{1}{2}}(\Gmt)} \qfor \frac{1}{2} \le \sigma \le s-\frac{1}{2}.
	\end{equation*}
	Furthermore, product and uniform elliptic estimates yield that
	\begin{equation*}
		\abs{\comm{\cDt}{\cN}\psi}_{H^{\sigma'}(\Gmt)} \lesssim_{\Lambda_s} \cP{\abs{\Gamma}_{H^{1+\sigma'}}} \cdot \abs{\vbu}_{H^{1+\sigma'}(\Gmt)}\abs{\psi}_{H^{1+\sigma'}(\Gmt)} \qfor \sigma'>0,
	\end{equation*}
	where $\mathcal{P}$ is a generic polynomial. In view of the following relation:
	\begin{equation*}
		\comm{\cDt}{\cN^k} = \sum_{l+m=k-1} \cN^{l}\comm{\cDt}{\cN}\cN^{m},
	\end{equation*}
	one can derive the estimate (here $k \ge 1$):
	\begin{equation}\label{comm est Dt N}
		\abs{\comm{\cDt}{\cN^k}}_{\scL\colon H^{\sigma'}(\Gmt) \to H^{\sigma'-k}(\Gmt)} \lesssim_{\Lambda_s} \cP{\abs{\Gamma}_{H^{\sigma'}}}\cdot \qty(\abs{\vbu}_{H^{s-\frac{1}{2}}(\Gmt)} + \abs{\vbu}_{H^{\sigma'}(\Gmt)} ) \qfor \sigma'\ge k - \frac{1}{2}.
	\end{equation}
	Similarly, standard calculations yield that:
	\begin{equation}
		\comm{\grad_\vtau}{\cN}\psi = 2\grad_\vn\lap^{-1}\dive(\grad\psi_{\cH}\vdot\grad\vtau_{\cH}) - \grad\psi_{\cH}\vdot\cN\vtau + \varkappa\grad_\vtau\psi,
	\end{equation}
	which leads to
	\begin{equation}
		\abs{\comm{\grad_\vtau}{\cN^k}}_{\scL\colon H^{\sigma''}(\Gmt) \to H^{\sigma''-k}(\Gmt)} \lesssim_{\Lambda_s} \cP{\abs{\Gamma}_{H^{1+\sigma''}}} \qfor \sigma''\ge k- \frac{1}{2}.
	\end{equation}
	In particular, Lemma \ref{lem reg Gam} implies the estimate for $s'\ge\frac{1}{2}$:
	\begin{equation}
		\begin{split}
			&\abs{\comm{\mu\grad_\vtau}{\cN}}_{\scL \colon H^{s'}(\Gmt) \to H^{s'-1}(\Gmt)} \\ 
			&\quad\le \cP{\abs*{\mu\varkappa}_{H^{s'-1}(\Gmt)}, \abs*{\mu\varkappa}_{H^{\frac{1}{2}+}(\Gmt)}, \abs{\varkappa}_{H^{\frac{1}{2}+}(\Gmt)}, \abs{\Gmt}_{H^{s'}},  \abs*{\mu}_{H^{s'}(\Gmt)}, \abs*{\mu}_{H^{\frac{3}{2}+}(\Gmt)}},
		\end{split}
	\end{equation}
	where $\mathcal{P}$ is a generic polynomial determined by $\Lambda_s$.
	
	\subsubsection{Lower order energy estimates}\label{set lin low est}
	To derive the evolution estimates, one may first take $\cN\cDt\vartheta$ as the test function to \eqref{lin eqn kp}. Due to the symmetry of $\cN$ on $L^2(\Gmt)$, it is natural to consider the following $L^2$-type energy (which is actually at the $H^{\frac{1}{2}}$-level):
	\begin{equation}
		\Elin^0(\vartheta) \coloneqq \int_{\Gmt} \abs*{\cN^{\frac{1}{2}}\cDt\vartheta}^2 + \alpha\abs*{\grad_\vtau\cN\vartheta}^2 + \abs*{\cN^{\frac{1}{2}}(\mu\grad_\vtau\vartheta)}^2 + \nu \abs{\cN\vartheta}^2 \dd{\ell},
	\end{equation}
	One can first observe that it satisfies the bound:
	\begin{equation}\label{bound Elin}
		\abs{\Elin^0(\vartheta)} \lesssim_{\Lambda_s} \abs{\cDt\vartheta}_{H^{\frac{1}{2}}(\Gmt)}^2 + \qty(1+\abs{\Gamma}_{H^2}^2)\cdot\alpha\abs{\vartheta}_{H^2(\Gmt)}^2 + \abs{\mu\grad_\vtau\vartheta}_{H^{\frac{1}{2}}(\Gmt)}^2 + \abs{\nu}_{L^\infty}\abs{\vartheta}_{H^1(\Gmt)}^2.
	\end{equation}
	It hints us to introduce the following expression:
	\begin{equation}\label{def Mlin}
		\Mlin^0(\vartheta) \coloneqq \abs{\cDt\vartheta}_{H^{\frac{1}{2}}(\Gmt)}^2 +  \alpha\abs{\vartheta}_{H^2(\Gmt)}^2 + \abs{\mu\grad_\vtau\vartheta}_{H^{\frac{1}{2}}(\Gmt)}^2 + \abs{\vartheta}_{H^1(\Gmt)}^2.
	\end{equation}
	Next, to control the time derivative of the energy, one can first rewrite it by invoking the symmetry of $\cN$ as:
	\begin{equation}\label{Elin'}
		\Elin^0(\vartheta) = \int_{\Gmt} \cDt\vartheta\cdot\cN\cDt\vartheta + \alpha\abs*{\grad_\vtau\cN\vartheta}^2 + \mu\grad_\vtau\vartheta\cdot\cN(\mu\grad_\vtau\vartheta)  + \nu \abs{\cN\vartheta}^2 \dd{\ell}.
	\end{equation}
	Then, it is standard to calculate that:
	\begin{equation}
		\begin{split}
			\dv{t}\Elin^0(\vartheta) =\, &\int_{\Gmt} \cDt\cDt\vartheta\cdot \cN\cDt\vartheta + \cDt\vartheta\cdot\cDt(\cN\cDt\vartheta) \dd{\ell} \quad \dasharrow \boxed{I_1} \\
			&+ \int_{\Gmt} 2\alpha\cDt(\grad_\vtau\cN\vartheta)\cdot\grad_\vtau\cN\vartheta \dd{\ell} \quad \dasharrow \boxed{I_2}\\
			&+ \int_{\Gmt} \cDt(\mu\grad_\vtau\vartheta)\cdot\cN(\mu\grad_\vtau\vartheta) + \mu\grad_\vtau\vartheta\cDt\cN(\mu\grad_\vtau\vartheta) \dd{\ell} \quad \dasharrow \boxed{I_3}\\ 
			&+ \int_{\Gmt} \cDt\nu \cdot \abs{\cN\vartheta}^2 + 2\nu\cDt\cN\vartheta\cdot\cN\vartheta \dd{\ell} \quad \dasharrow \boxed{I_4}\\
			&+ \int_{\Gmt} \boxed{\text{Integrand of \eqref{Elin'}}} \cdot(\grad_\vtau\vbu\vdot\vtau) \dd{\ell}. \quad \dasharrow \boxed{I_5}
		\end{split}
	\end{equation}
	One can first control $I_5$ by using the product estimates and the properties of $\cN$ that
	\begin{equation}
		\begin{split}
			\abs{I_5} \lesssim_{\Lambda_s}\, &\abs{(\grad_\vtau\vbu\vdot\vtau)\cDt\vartheta}_{H^{\frac{1}{2}}(\Gmt)}\abs{\cN\cDt\vartheta}_{H^{-\frac{1}{2}}(\Gmt)} + \alpha\abs{\grad_\vtau\vbu}_{L^\infty}\abs{\cN\vartheta}_{H^1(\Gmt)}^2 + \\
			&+\abs{(\grad_\vtau\vbu\vdot\vtau)\mu\grad_\vtau\vartheta}_{H^{\frac{1}{2}}(\Gmt)}\abs{\cN(\mu\grad_\vtau\vartheta)}_{H^{-\frac{1}{2}}(\Gmt)} + \abs{\grad_\vtau\vbu}_{L^\infty}\abs{\nu}_{L^\infty}\abs{\vartheta}_{H^{1}(\Gmt)}^2 \\
			\lesssim_{\Lambda_s}\, &(1+\abs{\nu}_{L^\infty})\abs{\vbu}_{H^{s-\frac{1}{2}}} \cdot \Mlin^0(\vartheta).
		\end{split}
	\end{equation}
	The integral $I_1$ can be computed as
	\begin{equation}
		\begin{split}
			I_1 =\, & \int_{\Gmt} 2 \cDt\cDt\vartheta \cdot \cN\cDt\vartheta \dd{\ell} \quad \dasharrow \boxed{\text{good}} \\
			&+ \int_{\Gmt} \cDt\vartheta \cdot \comm{\cDt}{\cN} \cDt\vartheta \dd{\ell}, \quad \dasharrow \boxed{\err I_1}
		\end{split}
	\end{equation}
	here the `good' term means that it can be canceled through using the equation \eqref{lin eqn kp}. The error term can be controlled by
	\begin{equation}
		\abs{\err I_1} \lesssim_{\Lambda_s} \abs{\cDt\vartheta}_{H^{\frac{1}{2}}(\Gmt)}\abs{\comm{\cDt}{\cN}\cDt\vartheta}_{H^{-\frac{1}{2}}(\Gmt)} \lesssim_{\Lambda_s} \abs{\cDt\vartheta}_{H^{\frac{1}{2}}(\Gmt)}^2.
	\end{equation}
	The second integral can be calculated through
	\begin{equation}
		\begin{split}
			I_2 =\, &-\int_{\Gmt} 2\alpha\cN\cDt\vartheta \cdot \grad_\vtau\grad_\vtau\cN\vartheta \dd{\ell} \quad \dasharrow \boxed{\text{good}} \\
			&+2\alpha\int_{\Gmt} \comm{\cDt}{\grad_\vtau}\cN\vartheta \cdot \grad_\vtau\cN\vartheta + \grad_\vtau\comm{\cDt}{\cN}\vartheta \cdot \grad_\vtau\cN\vartheta \dd{\ell}, \quad\dasharrow\boxed{\err I_2}
		\end{split}
	\end{equation}
	and the error term is bounded by
	\begin{equation}
			\abs{\err I_2} \lesssim_{\Lambda_s} \qty(1+\abs{\Gamma}_{H^2}^2)\cdot\alpha\abs{\vartheta}_{H^2(\Gmt)}^2.
	\end{equation}
	Next, $I_3$ can be decomposed into:
	\begin{equation}
		\begin{split}
			I_3 =\, & -\int_{\Gmt} 2 \cN\cDt\vartheta \cdot \grad_{\mu\vtau}\grad_{\mu\vtau}\vartheta \dd{\ell} \quad \dasharrow\boxed{\text{good}} \\
			&+ 2\int_{\Gmt} \cN(\mu\grad_\vtau\vartheta) \cdot \qty\big(\cDt\mu\cdot\grad_\vtau\vartheta - \cDt\vartheta\cdot\grad_\vtau\mu + \mu\comm{\cDt}{\grad_\vtau}\vartheta) \dd{\ell} \quad\dasharrow\boxed{\err I_3} \\
			&+ \underbrace{\int_{\Gmt} \mu\grad_\vtau\vartheta\cdot \comm{\cDt}{\cN}(\mu\grad_\vtau\vartheta) - 2\cDt\vartheta\cdot\qty\big(\mu\comm{\grad_\vtau}{\cN}(\mu\grad_\vtau\vartheta) + \comm{\mu}{\cN}\grad_\vtau(\mu\grad_\vtau\vartheta)) \dd{\ell}}_{\boxed{\err I_3}},
		\end{split}
	\end{equation}
	Due to the commutator formula \eqref{comm Dt dtau}, one has
	\begin{equation}\label{def lt mu}
		\cDt\mu\cdot\grad_\vtau\vartheta + \mu\comm{\cDt}{\grad_\vtau}\vartheta = \underbrace{\qty\big(\cDt\mu-(\grad_\vtau\vbu\vdot\vtau)\mu)}_{\eqqcolon \cL_t \mu} \cdot \grad_\vtau\vartheta.
	\end{equation}
	Then, it can be derived from the product and commutator estimates that
	\begin{equation}
		\abs{\err I_3} \lesssim_{\Lambda_s} \abs{\mu\grad_\vtau\vartheta}_{H^{\frac{1}{2}}(\Gmt)} \qty(\abs{(\cL_t\mu)\grad_{\vtau}\vartheta}_{H^{\frac{1}{2}}(\Gmt)} + \abs{\mu\grad_\vtau\vartheta}_{H^{\frac{1}{2}}(\Gmt)} + \abs{\mu}_{H^{s-\frac{1}{2}}}\abs{\cDt\vartheta}_{H^{\frac{1}{2}}(\Gmt)}).
	\end{equation}
	Similarly, one can show that
	\begin{equation}
		\begin{split}
			I_4 = &\int_{\Gmt}2\cN\cDt\vartheta \cdot \nu\cN\vartheta \dd{\ell} \quad \dasharrow\boxed{\text{good}} \\
			&+\int_{\Gmt} \cDt\nu\cdot \abs{\cN\vartheta}^2 + 2\nu \comm{\cDt}{\cN}\vartheta \cdot\cN\vartheta \dd{\ell}, \quad \dasharrow \boxed{\err I_4}
		\end{split}
	\end{equation}
	with
	\begin{equation}
		\abs{\err I_4} \lesssim_{\Lambda_s} \qty\big(\abs{\cDt\nu}_{L^\infty} + \abs{\nu}_{L^\infty}) \cdot \abs{\vartheta}_{H^{1}(\Gmt)}^2.
	\end{equation}
	
	Combining the above estimates, one finally obtains
	\begin{equation}
		\begin{split}
			\abs{\dv{t}\Elin^0(\vartheta)} \lesssim_{\Lambda_s}\, &\qty\big((1+\abs{\nu}_{L^\infty})(1+\abs{\vbu}_{H^{s-\frac{1}{2}}}) + \abs*{\mu}_{H^{s-\frac{1}{2}}} + \abs{\cDt\nu}_{L^\infty}) \Mlin^0(\vartheta) + \\
			&\quad+ \alpha\cP{\abs{\Gmt}_{H^2}}\abs{\vartheta}_{H^2(\Gmt)}^2 + \abs{\cL_t\mu \cdot \grad_\vtau\vartheta}_{H^{\frac{1}{2}}(\Gmt)}^2 + \abs{f}_{H^{\frac{1}{2}}(\Gmt)}^2,
		\end{split}
	\end{equation}
	where $\mathcal{P}$ is a generic polynomial and $\cL_t\mu$ is defined in \eqref{def lt mu}.
	
	\subsubsection{Higher order estimates}
	In order to preserve the symmetry of \eqref{lin eqn kp}, one can now take  $\grad_\vtau\grad_\vtau\cN\cDt\vartheta$ to be its test function. Thus, it is natural to consider the following energy functional:
	\begin{equation}\label{Elin1}
		\Elin^1(\vartheta) \coloneqq \int_{\Gmt} \abs*{(\slashed{\lap}\cN)^\frac{1}{2}\cDt\vartheta}^2 + \alpha \abs{\grad_\vtau\grad_\vtau\cN\vartheta}^2 + \abs*{(\slashed{\lap}\cN)^{\frac{1}{2}}(\mu\grad_\vtau\vartheta)}^2 + \nu \abs{\grad_\vtau\cN\vartheta}^2 \dd{\ell}.
	\end{equation}
	The corresponding bound is defined as
	\begin{equation}
		\Mlin^1(\vartheta) \coloneqq \abs{\cDt\vartheta}_{H^{\frac{3}{2}}(\Gmt)}^2 + \alpha\abs{\vartheta}_{H^{3}(\Gmt)}^2 + \abs{\mu\grad_\vtau\vartheta}_{H^{\frac{3}{2}}(\Gmt)}^2 + \abs{\vartheta}_{H^2(\Gmt)}^2.
	\end{equation}
	Similar to \eqref{Elin'}, the energy \eqref{Elin1} can be rewritten as
	\begin{equation*}
			\Elin^1(\vartheta) = \int_{\Gmt} \cDt\vartheta\cdot\grad_\vtau\grad_\vtau\cN\cDt\vartheta + \alpha \abs{\grad_\vtau\grad_\vtau\cN\vartheta}^2 + \mu\grad_\vtau\vartheta\cdot\grad_\vtau\grad_\vtau\cN(\mu\grad_\vtau\vartheta) + \nu \abs{\grad_\vtau\cN\vartheta}^2 \dd{\ell}.
	\end{equation*}
	Thus, the calculations in \S\ref{set lin low est} yield that
	\begin{equation}
		\begin{split}
			\abs{\dv{t}\Elin^1 (\vartheta)} \le\, & \cP{\abs{\varkappa}_{H^{s-\frac{1}{2}}}, \alpha\abs{\varkappa}_{H^2}^2, \abs*{\mu\varkappa}_{H^{\frac{3}{2}}}, \abs*{(\vbu, \mu, \nu)}_{H^{3}}, \abs{\cDt\nu}_{L^\infty}} \cdot \Mlin^1(\vartheta) + \\
			&\quad + \cP{\abs{\Gmt}_{H^{2}}}\cdot\qty(\abs*{\cL_t\mu \cdot \grad_\vtau\vartheta}_{H^{\frac{3}{2}}(\Gmt)}^2 +\abs{f}_{H^{\frac{3}{2}}(\Gmt)}^2), 
		\end{split}
	\end{equation}
	where $\mathcal{P}$ represents generic polynomials determined by $\Lambda_s$.
	
	Inductively, for each integer $m \ge 0$, one can take the test function $\slashed{\lap}^m \cN\cDt\vartheta$ and consider the corresponding energy functional:
	\begin{equation}
		\Elin^m(\vartheta) \coloneqq \int_{\Gmt} \abs*{(\slashed{\lap}^m\cN)^{\frac{1}{2}}\cDt\vartheta}^2 + \alpha\abs*{(\grad_\vtau)^{1+m}\cN\vartheta}^2 + \abs*{(\slashed{\lap}^m\cN)^{\frac{1}{2}}(\mu\grad_\vtau\vartheta)}^2 + \nu\abs*{(\grad_\vtau)^m\cN\vartheta}^2 \dd{\ell},
	\end{equation}
	together with its bound
	\begin{equation}
		\Mlin^m(\vartheta) \coloneqq \abs{\cDt\vartheta}_{H^{m+\frac{1}{2}}(\Gmt)}^2 + \alpha\abs{\vartheta}_{H^{m+2}(\Gmt)}^2 + \abs{\mu\grad_\vtau\vartheta}_{H^{m+\frac{1}{2}}(\Gmt)}^2 + \abs{\vartheta}_{H^{m+1}(\Gmt)}^2.
	\end{equation}
	Therefore, when $m \ge 2$, there holds the propagation estimate:
	\begin{equation}
		\begin{split}
			\abs{\dv{t}\Elin^m(\vartheta)} \le\, & \cP{\abs{\varkappa}_{H^{m+\frac{1}{2}}}, \alpha\abs{\varkappa}_{H^{m+1}}^2, \abs*{\mu\varkappa}_{H^{m+\frac{1}{2}}}, \abs*{(\vbu, \mu, \nu)}_{H^{m+2}}, \abs{\cDt\nu}_{L^\infty}} \cdot \Mlin^m(\vartheta) + \\
			&\quad+ \cP{\abs{\Gmt}_{H^{m+1}}} \qty(\abs*{\cL_t\mu \cdot \grad_\vtau\vartheta}_{H^{m + \frac{1}{2}}(\Gmt)}^2 + \abs{f}_{H^{m+\frac{1}{2}}(\Gmt)}^2),
		\end{split}
	\end{equation}
	where $\mathcal{P}$ represents generic polynomials determined by $\Lambda_s$ and $m$.
	
	\subsubsection{Modifications when $\alpha > 0$}\label{sec est alpha>0}
	Note that, if $\alpha$ is assumed to be strictly positive, then the term $\nu\cN\vartheta$ can be regarded as a remainder during the energy estimates. Thus, one can take $\cN^m\slashed{\lap}\cN\cDt\vartheta$ to be the test function to \eqref{lin eqn kp} and consider the energy functional
	\begin{equation}
		\wt{\Elin^m}(\vartheta) \coloneqq \int_{\Gmt} \abs*{(\cN^m\slashed{\lap}\cN)^{\frac{1}{2}}\cDt\vartheta}^2 + \alpha \abs{\cN^\frac{m}{2}\slashed{\lap}\cN\vartheta}^2 + \abs*{(\cN^m\slashed{\lap}\cN)^{\frac{1}{2}}(\mu\grad_\vtau\vartheta)}^2 \dd{\ell},
	\end{equation}
	together with its bound:
	\begin{equation}
		\wt{\Mlin^m}(\vartheta) \coloneqq \abs{\cDt\vartheta}_{H^{\frac{m}{2}+\frac{3}{2}}(\Gmt)}^2 + \alpha\abs{\vartheta}_{H^{\frac{m}{2}+3}(\Gmt)}^2 + \abs{\mu\grad_\vtau\vartheta}_{H^{\frac{m}{2}+\frac{3}{2}}(\Gmt)}^2.
	\end{equation}
	Then, the calculations in \S\ref{set lin low est} yield that
	\begin{equation}
			\abs{\dv{t}\wt{\Elin^m}(\vartheta)} \le \cP{\abs{\varkappa}_{H^{\frac{m}{2}+2}}, \abs*{(\vbu, \mu, \nu)}_{H^{\frac{m}{2}+3}}, \abs{\cL_t\mu}_{H^{\frac{m}{2}+\frac{3}{2}}}}\cdot\wt{\Mlin^m}(\vartheta) +  \cP{\abs{\Gmt}_{H^{\frac{m}{2}+\frac{5}{2}}}}\cdot\abs{f}_{H^{\frac{m}{2}+\frac{3}{2}}(\Gmt)}^2,
	\end{equation}
	where $\mathcal{P}$ represents generic polynomials determined by $\Lambda_s$, $m$, and $\alpha > 0$.

	\subsubsection{Modifications when $\alpha = 0$ and $\mu$ is non-degenerate}\label{sec est h non-degenerate}
	Suppose now that $\alpha = 0$, and there exists a constant $\lambda_0 > 0$ so that
	\begin{equation}
		\abs*{\mu} \ge \lambda_0 \qq{for all} t.
	\end{equation}
	Then, the term $\nu\cN\vartheta$ can also be treated as a remainder. Particularly, for each $m \ge 0$, one can take $\cN^m\cDt\vartheta$ to be the test function and consider the energy functional:
	\begin{equation}
		\wh{\Elin^m}(\vartheta) \coloneqq \int_{\Gmt} \abs*{\cN^{\frac{m}{2}}\cDt\vartheta}^2 + \abs*{\cN^{\frac{m}{2}}(\mu\grad_\vtau\vartheta)}^2 \dd{\ell},
	\end{equation}
	together with its bound:
	\begin{equation}
		\wh{\Mlin^m}(\vartheta) \coloneqq \abs{\cDt\vartheta}_{H^{\frac{m}{2}}(\Gmt)}^2 + \abs{\vartheta}_{H^{\frac{m}{2}+1}(\Gmt)}^2.
	\end{equation}
	Similar computations to those in \S\ref{set lin low est} yield that
	\begin{equation}
		\abs{\dv{t}\wh{\Elin^m}(\vartheta)} \le \cP{\abs{\varkappa}_{H^{\frac{m}{2}+\frac{1}{2}}}, \abs*{(\vbu, \mu, \nu)}_{H^{\frac{m}{2}+2}}, \abs*{\cL_t\mu}_{H^{\frac{m}{2}}}}\cdot\wh{\Mlin^m}(\vartheta) + \cP{\abs{\Gmt}_{H^{\frac{m}{2}+1}}}\cdot\abs{f}_{H^{\frac{m}{2}}(\Gmt)}^2.
	\end{equation}

	\subsection{Evolutions of the Vorticity and Current}
	Denote by
	\begin{equation}
		\Dt^\pm \coloneqq \Dt \pm \grad_\vh \qand \varpi_\pm \coloneqq \Curl(\vv\pm\vh).
	\end{equation}
	Then, it is obvious that \eqref{MHD} yields
	\begin{equation}
		\Dt^\pm (\vv\mp\vh) + \grad{p} = \vb{0},
	\end{equation}
	which lead to (in Cartesian's coordinates):
	\begin{equation}\label{evo eqn curl pre}
		\Dt^\pm \varpi_{\mp} + \sum_{k=1}^2 \pd_1 (v^k \pm h^k)\pd_k(v^2 \mp h^2) - \pd_2(v^k \pm h^k)\pd_k(v^1 \mp h^1) = 0.
	\end{equation}
	For the simplicity of notations, we denote by $\cB(\pd\vb{a}, \pd\vb{b})$ a bilinear form:
	\begin{equation}\label{def cB}
		\cB(\pd\vb{a}, \pd\vb{b}) = \sum_{k=1}^2 \pd_1a^k\pd_kb^2 - \pd_2a^k\pd_kb^1.
	\end{equation}
	Then, one can rewrite \eqref{evo eqn curl pre} as
	\begin{equation}\label{evo eqn curl}
		\Dt^\pm \varpi_{\mp} + \cB\qty[\pd(\vv\pm\vh), \pd(\vv\mp\vh)] = 0.
	\end{equation}
	
	\subsubsection{Simplified linearized equations}\label{sec curl lin est}
	As before, assume that $\{\Omt\}_{t}$ is a family of domains with $\pd\Omt \in \Lambda_s $, and $\vw \coloneqq \Omt \to \R^2$ is the evolution velocity of $\Omt$. In particular, $\vw_{\restriction_{\pd\Omt}}$ is one of the evolution velocities of $\pd\Omt \in \Lambda_s$. With a slight abuse of notations, denote by
	\begin{equation}
		\cDt \coloneqq \pdv{t} + \grad_\vw
	\end{equation}
	the material derivative along the trajectory of $\vw$. Consider the following linear evolution (transport) equation with a source term:
	\begin{equation}
		\cDt \zeta = g \qin \Omt.
	\end{equation}
	Then, it is clear that
	\begin{equation}
		\dv{t}\int_{\Omt} \abs{\zeta}^2 \dd{x} = \int_{\Omt} 2\zeta\cdot\cDt\zeta + \abs{\zeta}^2 (\dive\vw) \dd{x},
	\end{equation}
	which yields that
	\begin{equation}
		\abs{\dv{t}\norm{\zeta}_{L^2(\Omt)}^2} \lesssim( 1 + \norm{\grad\vw}_{L^\infty})\norm{\zeta}_{L^2(\Omt)}^2 + \norm{g}_{L^2(\Omega)}^2.
	\end{equation}
	
	For the higher order estimates, it is direct to derive that
	\begin{equation*}
		\cDt\pd^m\zeta  = \pd^mg + \comm{\grad_\vw}{\pd^m}\zeta,
	\end{equation*}
	and
	\begin{equation*}
		\begin{split}
			\norm{\comm{\grad_\vw}{\pd^m}\zeta}_{L^2(\Omt)} &\lesssim_{\Lambda_s} \norm{\pd\vw \cdot \pd\zeta}_{H^{m-1}(\Omt)} \\
			&\lesssim_{\Lambda_s} \qty(\norm{\vw}_{H^{2+}(\Omt)} + \norm{\vw}_{H^m(\Omt)})\cdot\norm{\zeta}_{H^m(\Omt)}.
		\end{split}
	\end{equation*}
	Therefore, $L^2$ energy estimate and interpolations imply that
	\begin{equation}\label{curl evo est}
		\abs{\dv{t}\norm{\zeta}_{H^\sigma(\Omt)}^2} \lesssim_{\Lambda_s} \qty(1 + \norm{\vw}_{H^{2+}(\Omt)}^2 + \norm{\vw}_{H^\sigma(\Omt)}^2 )\norm{\zeta}_{H^\sigma(\Omt)}^2 + \norm{g}_{H^\sigma(\Omt)}^2,
	\end{equation}
	for any real number $\sigma > 0$.
	
	\section{A Priori Estimates}\label{sec a priori est}
	In this section, we shall utilize the results in \S\ref{sec lin est} to derive the energy estimates for solutions to \eqref{PV problem}. We still implicitly assume that $\Omt$ is a bounded domain with it boundary $ \Gmt = \pd\Omt \in \Lambda_s$.
	
	First recall the following notation:
	\begin{equation*}
		\Dt \coloneqq \pdv{t} + \grad_\vv.
	\end{equation*}
	Then, one may infer from \eqref{DtDt k new} and \eqref{lin eqn kp} that the boundary part of higher-order energies can be expressed as (here $ m\ge 0$ is an integer):
	\begin{equation}
		\begin{split}
			E^m_{\text{bdry}} \coloneqq\, &\int_{\Gmt} \abs\big{(\slashed{\lap}^m\cN)^{\frac{1}{2}}\Dt\varkappa}^2 + \alpha \abs\big{\slashed{\grad}^{1+m}\cN\varkappa}^2 +(\grad_\vn\wt{q} - \grad_\vn q) \abs\big{\slashed{\grad}^{m}\cN\varkappa}^2 \dd{\ell}\\
			&+ \int_{\Gmt} \abs\big{(\slashed{\lap}^m\cN)^{\frac{1}{2}}(\grad_\vh\varkappa)}^2 + \abs\big{(\slashed{\lap}^m\cN)^{\frac{1}{2}}(\grad_\vH\varkappa)}^2 \dd{\ell}.
		\end{split}
	\end{equation}
	where $\varkappa$ is the signed curvature of $\Gmt$, and $q$ is the multiplier-type effective pressure given by \eqref{def q}.
	
	Due to the relation \eqref{Dt kappa} and the fact that $\vh \vdot \vn \equiv 0$ on $\Gmt$, the above boundary energy contains the information about the regularities of $\Gmt$ and the Neumann-type boundary data of $\vv$ and $\vh$. By utilizing the div-curl estimates, it remains to consider the following interior energy:
	\begin{equation}
		E_{\text{int}}^m \coloneqq \norm{\Curl (\vv\pm \vh)}_{H^{m+2}(\Omt)}^2.
	\end{equation}
	In summary, we define the total energy as:
	\begin{equation}
		E^m_\alpha \coloneqq 1 + \norm{\vv}_{L^2(\Omt)}^2 + \norm{\vh}_{L^2(\Omt)}^2 + 2\alpha\abs{\Gmt} + \norm{\vH}_{L^2(\cV_t)}^2 + E^m_{\text{bdry}} + E^m_{\text{int}},
	\end{equation}
	where $\abs{\Gmt}$ is the (Euclidean) length of $\Gmt$.
	
	It follows from \eqref{Dt kappa}, \eqref{def q}, \eqref{def wtq}, and the technical results in \S\ref{sec tech preparation} that the energy $E^m$ satisfies the bounds:
	\begin{equation}
		\begin{split}
			\abs{E^m_\alpha} \lesssim_{\Lambda_s}\, & \cP{\abs{\kappa}_{H^{m+\frac{1}{2}}}}\norm{\vv}_{H^{m+3}(\Omt)}^2 + \cP{\abs{\kappa}_{H^m}}\cdot\qty(\alpha\abs{\kappa}_{H^{m+2}(\Gmt)}^2 +\abs{\grad_\vh\kappa}_{H^{m+\frac{1}{2}}(\Gmt)}^2+\abs{\grad_\vH\varkappa}_{H^{m+\frac{1}{2}}(\Gmt)}^2) + \\
			&+ \cP{\norm{\vv}_{H^{2+}(\Omt)}, \norm{\vh}_{H^{2+}(\Omt)}, \abs{\kappa}_{H^m}, \abs{\cJ}_{H^{\frac{3}{2}+}(\cS)}} \abs{\kappa}_{H^{m+1}(\Gmt)}^2 + 1 + \alpha\abs{\Gmt} \\
			\lesssim_{\Lambda_s}\, &\cP{\abs{\kappa}_{H^{m+\frac{1}{2}}}, \norm{\vv}_{H^{2+}(\Omt)}, \norm{\vh}_{H^{2+}(\Omt)}, \abs{\cJ}_{H^{\frac{3}{2}+}(\cS)}} \times \\
			&\quad \times \qty(1 + \norm{\vv}_{H^{m+3}(\Omt)}^2 + \alpha\abs{\kappa}_{H^{m+2}(\Gmt)}^2 + \abs{\grad_\vh\kappa}_{H^{m+\frac{1}{2}}(\Gmt)}^2 + \abs{\grad_\vH\varkappa}_{H^{m+\frac{1}{2}}(\Gmt)}^2 + \abs{\kappa}_{H^{m+1}(\Gmt)}^2),
		\end{split}
	\end{equation}
	where $\mathcal{P}$ represents generic polynomials determined by $\Lambda_s$ and $m$. Moreover, div-curl and product estimates yield that
	\begin{equation*}
		\abs{\grad_\vH\varkappa}_{H^{m+\frac{1}{2}}(\Gmt)} \lesssim_{\Lambda_s} \abs{\vH}_{H^{m+1}(\Gmt)} \abs{\varkappa}_{H^{m+\frac{3}{2}}(\Gmt)} \lesssim_{\Lambda_s} \abs{\cJ}_{H^{m+1}(\cS)} \abs{\varkappa}_{H^{m+\frac{3}{2}}(\Gmt)}.
	\end{equation*}
	
	To simplify the notations, we denote by
	\begin{equation}
		\begin{split}
			M^m_\alpha \coloneqq\, &\norm{\vv}_{H^{m+3}(\Omt)}^2 + \norm{\vh}_{H^{m+3}(\Omt)}^2 + \qty(\abs{\pd_t \cJ}_{H^{m+\frac{3}{2}}(\cS)}^2 + \abs{\cJ}_{H^{m+\frac{5}{2}}(\cS)}^2)\cdot \qty(1 + \abs{\varkappa}_{H^{m+\frac{3}{2}}(\Gmt)}^2) + \\
			&+ \alpha\abs{\varkappa}_{H^{m+2}(\Gmt)}^2 + \abs{(\vh\vdot\grad)\varkappa}_{H^{m+\frac{1}{2}}(\Gmt)}^2 + \abs{\varkappa}_{H^{m+1}(\Gmt)}^2,
		\end{split}
	\end{equation}
	The main goal of \S\ref{sec a priori est} is to show the following a priori energy estimates:
	\begin{prop}\label{prop energy est}
		Suppose that $(\vv, \vh, \vH, \Omt)$ is a solution to \eqref{PV problem}, and $\pd\Omt \in \Lambda_s$ for some reference curve $\Gms$ and fixed constant $s > 2$. (Note here that $\Omt$ is not necessarily being simply-connected). Then, the solution admits the following energy estimates:
		\begin{itemize}
			\item {\bf Propagation estimate:}
			\begin{equation*}
				\abs{\dv{t}E^m_\alpha} \lesssim_{\Lambda_s} \cQ(M^m_\alpha),
			\end{equation*}
			where $\cQ$ is a generic polynomial determined by $\Lambda_s$ and $m$.
			\item {\bf Energy coercivity:}
			\begin{itemize}
				\item Case 1: $\alpha > 0$, there holds
				\begin{equation*}
					M^m_\alpha \le \cP{E^m_\alpha, \abs{\pd_t \cJ}_{H^{m+\frac{3}{2}}(\cS)}, \abs{\cJ}_{H^{m+\frac{5}{2}}(\cS)}},
				\end{equation*}
				where $\mathcal{P}$ is a generic polynomial determined $\alpha$, $\Lambda_s$, $\cS$, and $m$.
				\item Case 2: $\alpha = 0$ and $\abs{\vh} + \abs{\vH} \ge \lambda_0 > 0$ on $\Gmt$, there holds
				\begin{equation*}
					M^m_{\alpha=0} \le \cP{E^m_{\alpha=0}, \abs{\pd_t \cJ}_{H^{m+\frac{3}{2}}(\cS)}, \abs{\cJ}_{H^{m+\frac{5}{2}}(\cS)}},
				\end{equation*}
				where $\mathcal{P}$ is a generic polynomial determined by $\Lambda_s$,  $\lambda_0$, $\cS$, and $m$.
				\item Case 3: $\alpha = 0$, $\cJ\equiv 0$, $\abs{\vh}$ degenerates on $\Gmt$, and $(-\grad_\vn p)\ge c_0 > 0$ on $\Gmt$, there holds
				\begin{equation*}
					M^m_{\alpha=0} \le \cP{E^m_{\alpha=0}},
				\end{equation*}
				where $\mathcal{P}$ is a generic polynomial determined by $\Lambda_s$,  $c_0$, and $m$.
			\end{itemize}
		\end{itemize} 
	\end{prop} 
	
	\subsection{Propagation Estimate}\label{sec a pri est}
	
	First, it follows from the arguments in \S\ref{sec curl lin est} that
	\begin{equation}
		\abs{\dv{t}E_{\text{int}}^m} \lesssim_{\Lambda_s} \cQ\qty(\norm{\vv}_{H^{m+3}(\Omt)}, \norm{\vh}_{H^{m+3}(\Omt)}).
	\end{equation}
	For the $L^2$ energies, one can derive from \eqref{def Ephy} and \eqref{dt Ephy} that
	\begin{equation*}
		\abs{\dv{t} E_\text{phy}} \lesssim_{\Lambda_s} \abs{\cJ}_{L^2(\cS)}\qty( \norm{\pd_t\vH}_{L^2(\cV_t)} + \norm{\vH}_{H^1(\cV_t)} + \norm{\vh}_{H^1(\Omt)} + \norm{\vv}_{H^1(\Omt)}).
	\end{equation*}
	Particularly, \eqref{pre maxwell}-\eqref{BC} yield that
	\begin{equation*}
		\abs{\dv{t} E_\text{phy}} \lesssim_{\Lambda_s} M^m_{\alpha}.
	\end{equation*}
	Thus, it only remains to show the propagation estimate for $E^m_{\text{bdry}}$.
	
	Now, we consider the remainder term $\fkR_2$ appeared in \eqref{DtDt k new}. Indeed, it is direct to calculate that
	\begin{equation*}
		\begin{split}
			\fkR_2 =\, &  -\alpha\abs{\grad_\vtau\varkappa}^2 + \alpha\varkappa^2\cN\varkappa +   2\grad_\vtau(\abs{\vh}^2)\grad_\vtau\varkappa + \grad_\vtau(\abs{\vH}^2)\grad_\vtau\varkappa + \\
			&+ \varkappa^3\abs{\vh}^2 + \tfrac{1}{2}\varkappa^2\cN(\abs{\vH}^2) + 2\varkappa\grad_\vtau\grad_\vtau(\abs{\vh}^2) + \varkappa\grad_\vtau\grad_\vtau\abs{\vH}^2 + (\grad_\vn\wt{q})(\cN - \wt{\cN})\varkappa  + \\
			&-(\grad_\vn q)\varkappa(\varkappa + \vn\vdot\cN\vn)   - (\grad_{\cN\vn}q)\varkappa- (\grad_\vn\wt{q})\varkappa(\varkappa-\vn\vdot\wt{\cN}\vn) + (\grad_{\wt{\cN}\vn}\wt{q})\varkappa + \\
			&+  2\grad\vn_\cH\vb{\colon}\grad^2q + \grad_\vn \varrho  + 2\grad\vn_{\wt{\cH}}\vb{\colon}\grad^2\wt{q} + \grad_\vn\wt{\varrho}+ \grad_\vtau\grad_\vtau(\cN-\wt{\cN})\qty(\tfrac{1}{2}\abs{\vH}^2) +  \\
			&+ 2 (\grad_\vtau\vv \vdot\vn)(\vtau\vdot\grad_\vtau\grad_\vtau\vv) + 4(\grad_\vtau\vv\vdot\vtau)(\vn\vdot\grad_\vtau\grad_\vtau\vv) +  6\varkappa\abs{\grad_\vtau\vv\vdot\vtau}^2 - 3\varkappa\abs{\grad_\vtau\vv\vdot\vn}^2,
		\end{split}
	\end{equation*}
	where $q$, $\wt{q}$, $\varrho$, and $\wt{\varrho}$ are given by \eqref{def q}, \eqref{def wtq}, \eqref{def varrho}, and \eqref{def wt varrho}, respectively. It follows from routine but somehow tedious calculations that
	\begin{equation*}
		\begin{split}
			\abs{\fkR_2}_{H^{m+\frac{1}{2}}(\Gmt)} \lesssim_{\Lambda_s, \abs{\Gmt}_{H^{m+\frac{3}{2}}}}\, & \alpha\abs{\kappa}_{H^{m+\frac{3}{2}}(\Gmt)}\qty( \abs{\kappa}_{H^{\frac{3}{2}+}(\Gmt)} + \abs{\kappa}_{H^{m+\frac{3}{2}}(\Gmt)} + \abs{\varkappa}_{H^{\frac{1}{2}+}(\Gmt)}^2 + \abs{\varkappa}_{H^{m+\frac{1}{2}}(\Gmt)}^2) + \\ 
			&+ \abs{\vh}_{H^{m+2}(\Gmt)}\abs{\grad_\vh\kappa}_{H^{m+\frac{1}{2}}(\Gmt)} + \abs{\vH}_{H^{m+2}(\Gmt)}\abs{\grad_\vH\varkappa}_{H^{m+\frac{1}{2}}(\Gmt)} + \\
			&+\cQ\left(\abs{\kappa}_{H^{\frac{1}{2}+}}, \abs{\kappa}_{H^{m+\frac{1}{2}}}, \abs{\vh}_{H^{m+\frac{5}{2}}(\Gmt)}, \abs{\vv}_{H^{m+\frac{5}{2}}(\Gmt)}, \abs{\vH}_{H^{m+\frac{5}{2}}(\Gmt)}, \norm{\grad q}_{H^{m+2}(\Omt)}, \right. \\
			&\quad\quad \left. \norm{\grad q}_{H^{2+}(\Omt)}, \norm{\varrho}_{H^{m+2}(\Omt)}, \norm{\grad\wt{q}}_{H^{2+}(\cV_t)}, \norm{\grad\wt{q}}_{H^{m+2}(\cV_t)}, \norm{\wt{\varrho}}_{H^{m+2}(\cV_t)} \right),
		\end{split}
	\end{equation*}
	here $\cQ$ is a generic polynomial depending on $\Lambda_s$, and the dependence on $\abs{\Gmt}_{H^{m+\frac{3}{2}}}$ is of polynomial type. Recall that the multiplier $q$ is defined by \eqref{def q}, then the product and elliptic estimates yield
	\begin{equation*}
		\norm{\grad q}_{H^{m+2}(\Omt)} + \norm{\grad q}_{H^{2+}(\Omt)} \lesssim_{\Lambda_s, \abs{\Gmt}_{H^{m+\frac{3}{2}}}} \norm{\vv}_{H^{m+3}(\Omt)}^2 + \norm{\vh}_{H^{m+3}(\Omt)}^2.
	\end{equation*}
	Moreover, \eqref{def varrho} implies that
	\begin{equation*}
		\begin{split}
			\norm{\varrho}_{H^{m+2}(\Omt)} &\lesssim_{\Lambda_s} \norm{\lap\varrho}_{H^m(\Omt)} + \cP{\abs{\Gmt}_{H^{m+\frac{3}{2}}}}\norm{\varrho}_{L^2(\Omt)} \\
			&\lesssim_{\Lambda_s} \cQ\qty(\abs{\kappa}_{H^{m+\frac{1}{2}}(\Gmt)}, \norm{\vv}_{H^{m+3}(\Omt)}, \norm{\vh}_{H^{m+3}(\Omt)}),
		\end{split}
	\end{equation*}
	for some generic polynomial $\cQ$.
	Similarly, it follows that
	\begin{equation*}
		\norm{\grad\wt{q}}_{H^{2+}(\cV_t)} + \norm{\grad\wt{q}}_{H^{m+2}(\cV_t)} + \norm{\wt{\varrho}}_{H^{m+2}(\cV_t)} \lesssim_{\Lambda_s} \cQ\qty(\abs{\varkappa}_{H^{m+\frac{1}{2}}(\Gmt)}, \norm{\vH}_{H^{m+3}(\cV_t)}).
	\end{equation*}
	Uniform div-curl estimates yield that
	\begin{equation*}
		\norm{\vH}_{H^{m+3}(\cV_t)} \le \cP{\abs{\varkappa}_{H^{m+\frac{1}{2}}(\Gmt)}} \cdot \abs{\cJ}_{H^{m+\frac{5}{2}}(\cS)}
	\end{equation*}
	for some generic polynomial $\mathcal{P}$ determined by $\Lambda_s$, $\cS$, and $m$.
	
	Thus, one can simplify the estimate of $\fkR_2$ to be
	\begin{equation}
		\begin{split}
			\abs{\fkR_2}_{H^{m+\frac{1}{2}}(\Gmt)} \lesssim_{\Lambda_s, \abs{\Gmt}_{H^{m+\frac{3}{2}}}}\, & \alpha\abs{\kappa}_{H^{m+\frac{3}{2}}(\Gmt)}\qty( \abs{\kappa}_{H^{\frac{3}{2}+}(\Gmt)} + \abs{\kappa}_{H^{m+\frac{3}{2}}(\Gmt)} + \abs{\varkappa}_{H^{\frac{1}{2}+}(\Gmt)}^2 + \abs{\varkappa}_{H^{m+\frac{1}{2}}(\Gmt)}^2) + \\ 
			&+ \abs{\vh}_{H^{m+2}(\Gmt)}\abs{\grad_\vh\kappa}_{H^{m+\frac{1}{2}}(\Gmt)} + \abs{\vH}_{H^{m+2}(\Gmt)}\abs{\grad_\vH\varkappa}_{H^{m+\frac{1}{2}}(\Gmt)} + \\ &+\cQ\qty(\abs{\kappa}_{H^{\frac{1}{2}+}}, \abs{\kappa}_{H^{m+\frac{1}{2}}}, \norm{\vh}_{H^{m+3}(\Omt)}, \norm{\vv}_{H^{m+3}(\Omt)}, \abs{\cJ}_{H^{m+\frac{5}{2}}(\cS)}).
		\end{split}
	\end{equation}
	
	Moreover, \eqref{MHD} and \eqref{dt tau} yield that the $\cL_t\mu$ factor defined in \eqref{def lt mu} can be computed as
	\begin{equation}
		\cL_t (\vh\vdot\vtau) = \Dt(\vh\vdot\vtau) - (\grad_\vtau\vv\vdot\vtau)(\vh\vdot\vtau) = 0.
	\end{equation}
	Next, observe that
	\begin{equation}
		\cL_t (\vH\vdot\vtau) = \Dt(\vH\vdot\vtau) - (\grad_\vtau\vv\vdot\vtau)(\vH\vdot\vtau) = (\Dt\vH-\grad_\vH\vv)\vdot\vtau.
	\end{equation}
	Since $\pd_t \vH$ satisfies the div-curl system:
	\begin{equation*}
		\begin{cases*}
			\dive(\pd_t\vH) = 0 &in $\cV_t$, \\
			\Curl (\pd_t\vH) = 0 &in $\cV_t$, \\
			(\pd_t \vH)\vdot\vn = - \vH\vdot\Dt\vn - \grad_\vv\vH\vdot\vn &on $\Gmt$, \\
			\vN\smallwedge(\pd_t\vH) = \pd_t \cJ &on $\cS$,
		\end{cases*}
	\end{equation*}
	one can derive that
	\begin{equation}
		\abs{\cL_t (\vH\vdot\vtau)}_{H^{m+1}(\Gmt)} \le \cP{\abs{\varkappa}_{H^{m}(\Gmt)}} \cdot \norm{\vv}_{H^{m+3}(\Omt)}\cdot\qty(\abs{\cJ}_{H^{m+2}(\cS)} + \abs{\pd_t \cJ}_{H^{m+1}(\cS)})
	\end{equation}
	for some generic polynomial $\mathcal{P}$.
	
	The $\cDt(\grad_\vn q)$ term can now be calculated via \eqref{MHD}, \eqref{dt vn}, and \eqref{def q} that
	\begin{equation*}
		\begin{split}
			-\Dt(\grad_\vn q) &= (\grad_\vtau\vv \vdot \vn)\cdot \grad_\vtau q -\vn\vdot\grad\Dt q - \vn \vdot \comm{\Dt}{\grad} q \\
			&= \grad_\vn \comm{\Dt}{\lap^{-1}}\tr{(\grad\vv)^2 - (\grad\vh)^2} + \grad_\vn\lap^{-1}\Dt\tr{(\grad\vv)^2 - (\grad\vh)^2} + \grad\vv (\grad q, \vn).
		\end{split}
	\end{equation*}
	In view of the following commutator formula (cf. \cite{Shatah-Zeng2008-Geo})
	\begin{equation}
		\comm{\Dt}{\lap^{-1}}\psi = \lap^{-1}\qty(2\grad\vv\vb{\colon}\grad^2\lap^{-1}\psi + \lap\vv\vdot\grad\lap^{-1}\psi),
	\end{equation}
	one can conclude that
	\begin{equation}
		\abs{\Dt(\grad_\vn q)}_{L^\infty(\Gmt)} \le \cQ\qty(\norm{\vv}_{H^{2+}(\Omt)}, \norm{\vh}_{H^{2+}(\Omt)}).
	\end{equation}
	Similar estimates yield that
	\begin{equation}
		\abs{\Dt(\grad_\vn\wt{q})}_{L^\infty(\Gmt)} \le \cQ\qty(\abs{\cJ}_{H^{2}(\cS)}, \abs{\pd_t \cJ}_{H^1(\cS)}, \abs{\vv}_{H^2(\Gmt)}).
	\end{equation}
	
	Therefore, Proposition \ref{prop lin est theta} implies that
	\begin{equation}
		\abs{\dv{t}E^m_{\text{bdry}}} \lesssim_{\Lambda_s} \cQ(M^m_\alpha),
	\end{equation}
	which concludes the propagation estimate.

	\subsection{Energy Coercivity}
	When $\alpha > 0$, it is obvious from the ellipticity of Dirichlet-Neumann operators that
	\begin{equation}\label{bd eng coerc a>0}
		E^m_{\text{bdry}} \gtrsim_{\Lambda_s} \alpha\abs{\kappa}_{H^{m+2}}^2 + \cP{\abs{\Gmt}_{H^{m+2}}}\abs{\kappa}_{L^2(\Gmt)}^2, 
	\end{equation}
	which, together with interpolations and Lemma \ref{lem reg Gam}, yield that
	\begin{equation}
		\abs{\Gmt}_{H^{m+4}}^2 + \abs{\kappa}_{H^{m+2}(\Gmt)}^2 \le \cP{E^m_{\text{bdry}}},
	\end{equation}
	for some generic polynomial depending on $\Lambda_s$ and $\alpha$.
	
	Concerning the estimates for $\vv$ and $\vh$, one can utilize Lemma \ref{lem div-curl}. To confirm the boundary regularity, one can first derive from \eqref{Dt kappa} that
	\begin{equation}
		\begin{split}
			\abs{\vn\vdot\grad_\vtau\vv}_{H^{m+\frac{3}{2}}(\Gmt)} &\lesssim_{\Lambda_s, \abs{\Gmt}_{H^{m+\frac{3}{2}}}} \abs{\vn\vdot\grad_\vtau\vv}_{L^2(\Gmt)} + \abs{\grad_\vtau\vn\vdot\grad_\vtau\vv}_{H^{m+\frac{1}{2}}(\Gmt)} + \abs{\vn\vdot\grad_\vtau\grad_\vtau\vv}_{H^{m+\frac{1}{2}}(\Gmt)} \\
			&\lesssim_{\Lambda_s, \abs{\Gmt}_{H^{m+\frac{3}{2}}}} \abs{\grad_\vtau\vv}_{L^2(\Gmt)} + \abs{\kappa\grad_\vtau\vv\vdot\vtau}_{H^{m+\frac{1}{2}}(\Gmt)} + \abs{\Dt\kappa}_{H^{m+\frac{1}{2}}(\Gmt)} \\
			&\lesssim_{\Lambda_s, \abs{\Gmt}_{H^{m+\frac{3}{2}}}} \abs{\Dt\kappa}_{H^{m+\frac{1}{2}}(\Gmt)} + \cP{\abs{\kappa}_{H^{m+\frac{1}{2}}(\Gmt)}, \abs{\kappa}_{H^{\frac{1}{2}+}(\Gmt)}}\norm{\vv}_{H^{m+2}(\Omt)},
		\end{split}
	\end{equation}
	and
	\begin{equation}
		\begin{split}
			\abs{\vn\vdot\grad_\vtau\vh}_{H^{m+\frac{3}{2}}(\Omt)} &\lesssim_{\Lambda_s, \abs{\Gmt}_{H^{m+\frac{3}{2}}}} \abs{\vn\vdot\grad_\vtau\vh}_{L^s(\Gmt)} + \abs{\grad_\vtau\vn \vdot \grad_\vtau\vh}_{H^{m+\frac{1}{2}}(\Gmt)} + \abs{\vn\vdot\grad_\vtau\grad_\vtau\vh}_{H^{m+\frac{1}{2}}(\Gmt)}.
		\end{split}
	\end{equation}
	Since $\vh \vdot \vn \equiv 0$ on $\Gmt$, it holds that
	\begin{equation}
		\begin{split}
			\grad_\vtau\grad_\vtau\vh &= \grad_\vtau\grad_\vtau\qty\big[(\vh\vdot\vtau)\vtau] \\
			&= \qty\big[\grad_\vtau\grad_\vtau(\vh\vdot\vtau)-\kappa^2(\vh\vdot\vtau)]\vtau - \qty\big[2\kappa\grad_\vtau(\vh\vdot\vtau) + \grad_\vh\kappa]\vn,
		\end{split}
	\end{equation}
	which implies that
	\begin{equation}
		\abs{\vn\vdot\grad_\vtau\vh}_{H^{m+\frac{3}{2}}(\Omt)} \lesssim_{\Lambda_s, \abs{\Gmt}_{H^{m+\frac{3}{2}}}} \abs{\grad_\vh\kappa}_{H^{m+\frac{1}{2}}(\Gmt)} + \cP{\abs{\kappa}_{H^{m+\frac{1}{2}}(\Gmt)}, \abs{\kappa}_{H^{\frac{1}{2}+}(\Gmt)}}\norm{\vh}_{H^{m+2}(\Omt)}.
	\end{equation}
	
	Thus, Lemma \ref{lem div-curl} and interpolations lead to
	\begin{equation}
		\begin{split}
			\norm{\vv}_{H^{m+3}(\Omt)}^2 + \norm{\vh}_{H^{m+3}(\Omt)}^2 \lesssim_{\Lambda_s, \abs{\kappa}_{H^{m+\frac{1}{2}}(\Gmt)}}\, & \abs{\Dt\kappa}_{H^{m+\frac{1}{2}}(\Gmt)}^2 + \abs{\grad_\vh\kappa}_{H^{m+\frac{1}{2}}(\Gmt)}^2 + \\
			&+ \cQ\qty(\abs{\kappa}_{H^{m+\frac{1}{2}}(\Gmt)}, \norm{\vv}_{L^2(\Omt)}, \norm{\vh}_{L^2(\Omt)})
		\end{split}
	\end{equation}
	for some generic polynomial $\cQ$. These arguments conclude the coercivity for Case 1.
	
	When $\alpha = 0$ and $\abs{\vh} + \abs{\vH} \ge \lambda_0 > 0$ on $\Gmt$,  it follows that
	\begin{equation}\label{bd eng coerc a=0}
		\abs{\kappa}_{H^{m+\frac{3}{2}}(\Gmt)}^2 \lesssim_{\Lambda_s, \lambda_0, \abs{\Gmt}_{H^{m+\frac{3}{2}}}} 1 + \abs{\grad_\vh\kappa}_{H^{m+\frac{1}{2}}(\Gmt)}^2 + \abs{\grad_\vH\varkappa}_{H^{m+\frac{1}{2}}}^2 \le \cQ_{\lambda_0}\qty(E^m_{\alpha=0}),
	\end{equation}
	which finishes the proof for Case 2.
	
	Similarly, when $\alpha = 0 $, $\cJ \equiv 0$, and $\abs{\vh}$ degenerates on $\Gmt$, the Rayleigh-Taylor sign condition
	\begin{equation}
		-\grad_\vn q \ge c_0 > 0 \qq{on} \Gmt
	\end{equation}
	implies that
	\begin{equation}
		\abs{\kappa}_{H^{m+1}(\Gmt)}^2 \lesssim_{\Lambda_s, c_0, \abs{\Gmt}_{H^{m+1}}} E^m_{\alpha = 0},
	\end{equation}
	which leads to the coercivity for Case 3.
	
	\section{Construction of Iterates}\label{sec construct ite}
	We now turn to the construction of solutions to \eqref{PV problem}. As it would be much easier to construct iterations on fixed domains, we will pull back the free boundary problems into the reference frame through the coordinate map $\vPhi_{\Gmt}$ and its harmonic extension $\cX_{\Gmt}$ introduced in \S\ref{sec ref frame}. 
	
	For the simplicity of notations, from now on, we always assume that the index $s$ in $\Lambda_s$ satisfies $s \ge 3$.
	
	\subsection{Equations in the Reference Frame}
	
	\subsubsection{Relations between a curve and its curvature}
	Since the kinematics of $\Gmt$ is characterized through the evolution of $\varkappa$, it is necessary to quantitatively determine a curve $\Gmt \in \Lambda_s$ through its curvature $\varkappa$. On the other hand, $\Gmt$ is defined by a unique height function $\varphi_{\Gmt}\colon \Gms \to \R$, it suffices to find a one-to-one correspondence between $\varphi_{\Gmt}$ and $\varkappa\circ\vPhi_{\Gmt}$. Indeed, denote by
	\begin{equation}
		\mathscr{K}(\varphi_{\Gmt}) \coloneqq \varkappa_{\Gmt}\circ\vPhi_{\Gmt} \colon \Gms \to \R,
	\end{equation}
	which is just the pull-back of $\varkappa_{\Gmt}$ through the coordinate map. One can observe from \eqref{Dt kappa} that the variational derivative of $\mathscr{K}$ with respect to $\varphi_{\Gmt}$ is (here $\varkappa_*$ and $\slashed{\lap}_*$ are respectively the signed curvature of and the Laplace-Beltrami operator on $\Gms$)
	\begin{equation}
		\begin{split}
			\qty(\fdv{\mathscr{K}}{\varphi})_{\restriction_{\varphi = 0}} (\psi) &= -\vn_* \vdot \grad_{\vtau_*}\grad_{\vtau_*}(\psi\vn_*) - 2\varkappa_*\grad_{\vtau_*}(\psi\vn_*)\vdot\vtau_* \\
			&= \qty[-\varkappa_*^2 - \grad_{\vtau_*}\grad_{\vtau_*}]\psi = \qty(-\varkappa_*^2 - \slashed{\lap}_*)\psi,
		\end{split}
	\end{equation}
	which may have nontrivial kernels. Thus, one can consider the following modification (see also \cite{Shatah-Zeng2011})
	\begin{equation}\label{def ka}
		\fkk(\varphi_{\Gmt}) \coloneqq \underbrace{\varkappa_{\Gmt}\circ\vPhi_{\Gmt} + a^2\varphi_{\Gmt}}_{\eqqcolon \ka} \colon \Gms \to \R,
	\end{equation}
	where $a > 0$ is the ancillary parameter. (One simple candidate here is $a \coloneqq 1 + \abs{\varkappa_*}_{L^\infty}^2$).
	Then, it holds that
	\begin{equation}
		\qty(\fdv{\fkk}{\varphi})_{\restriction_{\varphi=0}}(\psi) = \qty[\qty(a^2 - \varkappa_*^2) - \slashed{\lap}_*]\psi,
	\end{equation}
	which is clearly a strictly positive self-adjoint operator on $L^2(\Gms)$. Namely, each curve in $\Lambda_s$ is uniquely determined by its ancillary curvature $\ka$ defined on the reference frame $\Gms$. More precisely, there holds the following lemma (cf. \cite[Lemma 2.2]{Shatah-Zeng2011}):
	\begin{lemma}\label{lem ka}
		Suppose that the bound $\delta_0$ in $\Lambda_s = \Lambda(\Gms, s-\frac{1}{2}, \delta_0)$ is small enough, and the ancillary constant $a \gg_{\Gms} 1$. Then, the map $\fkk$ is a $C^3$-diffeomorphism from $\Lambda_s \subset H^{s-\frac{1}{2}}(\Gms) \to H^{s-\frac{5}{2}}(\Gms)$. Denote by
		\begin{equation*}
			B_{\delta_1} \coloneqq \Set*{\ka \given \abs{\ka-\varkappa_*}_{H^{s-\frac{5}{2}}(\Gms)} < \delta_1}.
		\end{equation*}
		Then, so long as $\delta_1 \ll_{\Lambda_s} 1$, there holds the estimate:
		\begin{equation*}
			\abs{\fkk^{-1}}_{C^3\qty(B_{\delta_1} \to H^{s-\frac{1}{2}}(\Gms))} \lesssim_{\Lambda_s} 1.
		\end{equation*} 
		Moreover, if $\ka \in B_{\delta_1} \cap H^{\sigma-2}(\Gms)$ for some $s-\frac{1}{2} \le \sigma \le s+1$, it holds for $\max\{-\sigma, \sigma'-2\}\le\sigma''\le\sigma'\le\sigma$ that
		\begin{equation*}
			\abs{\qty(\var\fkk^{-1})_{\restriction_{\ka}}}_{\scL \colon H^{\sigma''}(\Gms) \to H^{\sigma'}(\Gms)} \lesssim_{\Lambda_s} a^{\sigma'-\sigma''-2}\qty(1+\abs{\ka}_{H^{\sigma-2}(\Gms)}).
		\end{equation*}
	\end{lemma}
	
	\subsubsection{Pull-back of evolution equations}
	In order to pull back \eqref{DtDt k new} into $\Gms$, one has to handle the pull-back of material derivatives. Indeed, the first boundary condition in \eqref{BC} implies that
	\begin{equation}\label{bc 1}
		(\vv\vdot\vn)\circ\vPhi_{\Gmt} = \pdv{\varphi_{\Gmt}}{t} \vn_* \vdot (\vn\circ\vPhi_{\Gmt}),
	\end{equation}
	which hints the decomposition:
	\begin{equation}\label{def vxi}
		\vv = \qty(\pdv{\varphi_{\Gmt}}{t} \vn_*)\circ \vPhi_{\Gmt}^{-1} + \underbrace{\qty[\vv -\qty(\pdv{\varphi_{\Gmt}}{t} \vn_*)\circ \vPhi_{\Gmt}^{-1} ]}_{\eqqcolon \vxi}.
	\end{equation}
	It follows from \eqref{bc 1} that the vector field $\vxi$ is tangential to $\Gmt$. Define $\vxi_* \in \textup{T}\Gms$ to be the push-out vector field of $\vxi$ through the inverse coordinate map, i.e.,
	\begin{equation}\label{def vxi*}
		\vxi_* \coloneqq \qty(\vPhi_{\Gmt}^{-1})_{\#} \vxi = \qty(\grad_{\vxi}\vPhi_{\Gmt}^{-1})\circ\vPhi_{\Gmt}.
	\end{equation}
	Then, it is clear that
	\begin{equation}\label{pull back dt}
		(\Dt\psi)\circ\vPhi_{\Gmt} = \pd_t(\psi\circ\vPhi_{\Gmt}) + \grad_{\vxi_*}(\psi\circ\vPhi_{\Gmt}).
	\end{equation}
	In particular, the vector field $\vxi_* \in \textup{T}\Gms$ can be regarded as an ancillary fields during the pull-back process.
	
	For the simplicity of notations, we denote by
	\begin{equation}
		\underline{\psi} \coloneqq \psi \circ \vPhi_{\Gmt}
	\end{equation}
	for any function $\psi$ defined on $\Gmt$. It follows from \eqref{pull back dt} that
	\begin{equation}\label{pull back evo eqn}
		\begin{split}
			\underline{\Dt\Dt\psi} &= \pd_t(\underline{\Dt\psi})+\grad_{\vxi_*}(\underline{\Dt\psi}) \\
			&= \pd_t\qty(\pd_t\underline{\psi} + \grad_{\vxi_*}\udl{\psi}) + \grad_{\vxi_*}(\pd_t\udl{\psi}+ \grad_{\vxi_*}\udl{\psi}) \\
			&= \pd_t\pd_t\udl{\psi} +2\grad_{\vxi_*}\pd_t\udl{\psi} + \grad_{\pd_t\vxi_*}\udl{\psi} + \grad_{\vxi_*}\grad_{\vxi_*}\udl{\psi}.
		\end{split}
	\end{equation}
	
	\subsubsection{Variations of ancillary vector fields}\label{sec var velocity}
	
	To analyze the vector field $\pd_t\vxi_*$, one may consider a more general situation. Assume that $\{\Omega_{t, \beta}\}_{t, \beta}$ is a bi-parameterized family of \emph{simply connected} domains with boundaries $\Gamma_{t, \beta} \in \Lambda_s$, and $\vv$ is the evolution velocity of $\Omega_{t, \beta}$ with respect to the parameter $t$, which is a vector field with constant divergence and satisfies the boundary condition \eqref{bc 1}.
	
	The main goal of \S\ref{sec var velocity} is to derive the relation between $\pd_\beta\vxi_*$ and $(\pd_\beta\ka, \pd_\beta\omega_*, \pd_t\pd_\beta\ka)$.
		
	Observe that, for a given function $\zeta_*(t, \beta)$ on $\Oms$, the following div-curl system
	\begin{equation}\label{div-curl prob}
		\begin{cases*}
			\dive \vw = \const\, \gamma &in $\Omega_{t, \beta}$,\\
			\Curl \vw = \zeta_* \circ \cX_{\Gamma_{t, \beta}}^{-1} &in $\Omega_{t, \beta}$, \\
			\vw \vdot \vn = \vn \vdot \qty(\pdv{\varphi_{\Gamma_{t, \beta}}}{t}\vn_*)\circ\vPhi_{\Gamma_{t, \beta}}^{-1} &on $\Gamma_{t, \beta}$
		\end{cases*}
	\end{equation}
	admits a unique solution satisfying the estimate in Lemma \ref{lem div-curl}. Here $\gamma(t, \beta)$ is a constant for each fixed $(t, \beta)$ ensuring the compatibility condition
	\begin{equation}\label{def gamma}
		\int_{\Gamma_{t, \beta}} \vn \vdot \qty(\pdv{\varphi_{\Gamma_{t, \beta}}}{t}\vn_*)\circ\vPhi_{\Gamma_{t, \beta}}^{-1} \dd{\ell} = \gamma\int_{\Omega_{t, \beta}} \dd{x}.
	\end{equation}
	Namely, the velocity field $\vv$ is uniquely determined by $\pd_t\varphi_{\Gamma_{t, \beta}}$ and $(\Curl \vv)\circ\cX_{\Gamma_{t, \beta}}$ through solving the div-curl problem \eqref{div-curl prob}. Note here that $\pd_t\varphi_{\Gamma_{t, \beta}}$ and $(\Curl\vv)\circ\cX_{\Gamma_{t, \beta}}$ are functions defined respectively on $\Gms$ and $\Oms$.
	
	For the simplicity of notations, we shall omit the subscript ``$\Gamma_{t, \beta}$'' for $\varphi$, $\vPhi$, and $\cX$. Furthermore, we denote by
	\begin{equation}
		\omega_* \coloneqq (\Curl \vv) \circ \cX \colon \Oms \to \R
	\end{equation}
	the pulled-back vorticity. 
	
	One may derive from \eqref{def vxi}-\eqref{def vxi*} that
	\begin{equation}
		\pdv{\beta}\vxi_* = \pdv{\beta}\qty[(\grad\vPhi)^{-1}\vdot\underline{\vxi}] = \pdv{\beta}\qty[(\grad\vPhi)^{-1}\vdot(\underline{\vv}-\pd_t\varphi\vn_*)].
	\end{equation}
	Recall that $\vPhi$ is defined through \eqref{def vPhi}, so it is routine to obtain that
	\begin{equation}\label{dbt vxi*}
		\begin{split}
			\pd_\beta\vxi_* &= -(\grad\vPhi)^{-1}\vdot\grad(\pd_\beta\varphi\vn_*)\vdot(\grad\vPhi)^{-1}\vdot(\underline{\vv}-\pd_t\varphi\vn_*) + (\grad\vPhi)^{-1}\vdot(\pd_\beta\underline{\vv}-\pd_t\pd_\beta\varphi\vn_*) \\
			&= - (\grad\vPhi)^{-1} \vdot \grad_{\vxi_*}(\pd_\beta\varphi\vn_*) + (\vPhi^{-1})_{\#}\qty[(\pd_\beta\underline{\vv}-\pd_t\pd_\beta\varphi\vn_*)\circ\vPhi^{-1}].
		\end{split}
	\end{equation}
	Since $\vv$ is assumed to satisfy \eqref{bc 1}, the vector field $(\pd_\beta\underline{\vv}-\pd_t\pd_\beta\varphi\vn_*)\circ\vPhi^{-1}$ is tangential to $\Gamma_{t, \beta}$, which implies that the usage of a ``push-out operator'' is legitimate. 
	
	Now, we turn to the decomposition of $\pd_\beta\underline{\vv}$. We shall extend the definition of $\underline{\vv}$ to the whole $\Oms$ through 
	\begin{equation*}
		\underline{\vv} \coloneqq \vv \circ \cX \colon \Oms \to \R^2.
	\end{equation*}
	Notice that $\cX$ can also be viewed as a vector field defined on $\R^2$. Then, by defining the vector field
	\begin{equation}
		\vb{b} \coloneqq (\pd_\beta \cX) \circ \cX^{-1} = \qty[\cH_{\Oms}(\pd_\beta\varphi\vn_*)] \circ \cX^{-1} \colon \Omega_{t, \beta} \to \R^2,
	\end{equation}
	it is obvious that
	\begin{equation}
		\pd_\beta (f \circ \cX) = (\pd_\beta f + \grad_{\vb{b}}f) \circ \cX
	\end{equation}
	for any function $f(t, \beta)$ defined in $\Omega_{t, \beta}$. For the simplicity of notations, we denote by
	\begin{equation}
		\Dbt \coloneqq \pd_\beta + \grad_{\vb{b}}
	\end{equation}
	the material derivative along the $\beta$-trajectory. Then, it is clear that
	\begin{equation}
		(\pd_\beta\underline{\vv})\circ\cX^{-1} = \Dbt\vv \qin \Omega_{t, \beta}.
	\end{equation}
	On the other hand, $\Dbt\vv$ is the unique solution to the following div-curl problem:
	\begin{equation}
		\begin{cases*}
			\dive\Dbt\vv = \comm{\dive}{\Dbt}\vv + \pd_\beta\gamma &in $\Omega_{t, \beta}$, \\
			\Curl\Dbt\vv = \comm{\Curl}{\Dbt}\vv + \Dbt\qty(\omega_*\circ\cX^{-1}) &in $\Omega_{t, \beta}$, \\
			\Dbt\vv \vdot \vn = \Dbt\qty[\vn\vdot(\pd_t\varphi\vn_*)\circ\vPhi^{-1}]-\vv\vdot\Dbt\vn &on $\Gamma_{t, \beta}$,
		\end{cases*}
	\end{equation}
	where $\gamma$ is the balancing constant defined by \eqref{def gamma}. Thus, commutator and transport formulae imply the following compatibility condition:
	\begin{equation}
		\int_{\Gamma_{t, \beta}} \Dbt\qty[\vn\vdot(\pd_t\varphi\vn_*)\circ\vPhi^{-1}]-\vv\vdot\Dbt\vn \dd{\ell} \equiv \int_{\Omega_{t, \beta}} \comm{\dive}{\Dbt}\vv + \pd_\beta\gamma \dd{x}.
	\end{equation}
	Standard div-curl estimates imply that, for $1 \le s' \le s-1$, there holds
	\begin{equation}
		\begin{split}
			\norm{\Dbt\vv}_{H^{s'}(\Omega_{t, \beta})} \lesssim_{\Lambda_s}\, &\norm{\pd_\beta\omega_*}_{H^{s'-1}(\Oms)} + \abs{\pd_t\pd_\beta\varphi}_{H^{s'-\frac{1}{2}}(\Gms)} + \\
			&+ \abs{\pd_\beta\varphi}_{H^{s'+\frac{1}{2}}(\Gms)}\qty(1 + \norm{\omega_*}_{H^{s-1}(\Oms)} + \abs{\pd_t\varphi}_{H^{s-\frac{1}{2}}(\Oms)}).
		\end{split}
	\end{equation}
	One can then obtain the estimate for $\pd_\beta\vxi_*$ through \eqref{dbt vxi*}.
	Recall that $\ka = \fkk(\varphi)$ and $\fkk$ is a $C^3$-diffeomorphism, which yield
	\begin{equation}\label{var varphi}
		\pd_\beta\varphi = \qty(\fdv{\fkk^{-1}}{\ka})_{\restriction_\ka}\!\!\pd_\beta\ka \qand \pd_t\pd_\beta\varphi = \qty(\fdv[2]{\fkk^{-1}}{\ka})_{\restriction_\ka} \!\!\pd_t\ka \pd_\beta\ka + \qty(\fdv{\fkk^{-1}}{\ka})_{\restriction_\ka}\!\!\pd_\beta\pd_t\ka.
	\end{equation}
	Namely, the linearity of previously discussed relations and the fact that $\ka(t, \beta)$ uniquely determines $\Gamma_{t, \beta}$ implies the existence of three linear operators, whose ranges are all $\textup{T}\Gms$, so that
	\begin{equation}\label{var vxi*}
		\pd_\beta\vxi_* = \underbrace{\vb{B}(\ka)\pd_\beta\pd_t\ka}_{\text{boundary evolution part}} + \underbrace{\vb{R}(\ka)\pd_\beta\omega_*}_{\text{rotation part}} + \underbrace{\vb{C}(\ka, \pd_t\ka, \omega_*)\pd_\beta\ka}_{\text{corrective part}}.
	\end{equation}
	Moreover, these three operators satisfies the following estimates (see also \cite[Appendix A]{Liu-Xin2023}):
	\begin{prop}\label{prop op BRC}
		With the hypothesis on constants $\delta_0$ and $a$ in Lemma \ref{lem ka}, suppose that $\ka \in B_{\delta_1} \subset H^{s-\frac{5}{2}}(\Gms)$. Then, there hold for $\frac{1}{2} \le \sigma \le s-\frac{3}{2}$ and $\sigma-2\le\sigma'\le\sigma$ that
		\begin{equation*}
				\abs{\vb{B}(\ka)}_{\scL \colon H^{\sigma'}(\Gms) \to H^{\sigma}(\Gms)} \lesssim_{\Lambda_s} a^{\sigma-\sigma'-2},
		\end{equation*}
			and
		\begin{equation*}
			\abs{\vb{R}(\ka)}_{\scL\colon H^{\sigma-\frac{1}{2}}(\Oms) \to H^{\sigma}(\Gms)} \lesssim_{\Lambda_s} 1.
		\end{equation*}
		Assume further that $\ka \in H^{s-\frac{3}{2}}(\Gms)$, $\pd_t\ka \in H^{s-\frac{5}{2}}(\Gms)$, and $\omega_* \in H^{s-1}(\Oms)$. Then, for $\frac{1}{2}\le\lambda\le s+\frac{1}{2}$, $\lambda-2\le \lambda'\le \lambda$, and $\sigma$ given previously, it follows that
			\begin{equation*}
				\abs{\vb{B}(\ka)}_{\scL\colon H^{\lambda'}(\Gms) \to H^{\lambda}(\Gms)} \lesssim_{\Lambda_s} a^{\lambda-\lambda'-2}\qty(1+\abs{\ka}_{H^{s-\frac{3}{2}}(\Gms)}),
			\end{equation*}
			and
			\begin{equation*}
				\abs{\vb{C}(\ka, \pd_t\ka, \omega_*)}_{\scL\colon H^{\sigma-1}(\Gms) \to H^{\sigma}(\Gms)} \lesssim_{\Lambda_s} 1 + \abs{\pd_t\ka}_{H^{s-\frac{5}{2}}(\Gms)} + \norm{\omega_*}_{H^{s-1}(\Oms)}.
			\end{equation*}
	\end{prop}
	
	\subsubsection{Evolution equation for $\ka$}
	With the help of \eqref{pull back evo eqn} and \eqref{var vxi*}, one may pull the evolution equation \eqref{DtDt k new} back to $\Gms$ through the coordinate map $\vPhi_{\Gmt}$.
	For the simplicity of notations, we introduce the following operator on $\Gms$:
	\begin{equation}\label{def scA}
		\scA(\ka)\psi \coloneqq \qty[(\grad_\vtau\grad_\vtau\cN)_{\restriction_{\Gamma}}(\psi\circ\vPhi_{\Gamma}^{-1})] \circ \vPhi_{\Gamma},
	\end{equation}
	which is exactly the pulled-back of operator $\slashed{\lap}\cN$ defined on $\Gamma$. For the vector field $\vh$, we define
	\begin{equation}\label{def vh*}
		\vh_* \coloneqq  \qty(\vPhi^{-1}_{\Gamma})_{\#} \vh = \qty(\grad_\vh\vPhi^{-1}_{\Gamma})\circ\vPhi_\Gamma \in \textup{T}\Gms
	\end{equation}
	to be the intrinsic pulled-back vector filed of $\vh \in \textup{T}\Gamma$. In particular, one has
	\begin{equation*}
		\vPhi_{\Gamma}^{\#}(\grad_\vh \phi) \equiv (\grad_\vh \phi)\circ\vPhi_{\Gamma} = \grad_{\vh_*}\underline{\phi} \qq{for any function} \phi \colon \Gamma \to \R.
	\end{equation*}
	Similarly, we define the pulled-back version of $\vH$ through
	\begin{equation}\label{def vH*}
		\vH_* \coloneqq  \qty(\vPhi^{-1}_{\Gamma})_{\#} \vH = \qty(\grad_\vH\vPhi^{-1}_{\Gamma})\circ\vPhi_\Gamma.
	\end{equation}
	
	With these preparations, one can derive the pulled-back version of \eqref{DtDt k new}:
	\begin{equation}\label{evo eqn kpudl}
		\begin{split}
			\pd_t\pd_t\udl{\varkappa} =\, &- 2\grad_{\vxi_*}\pd_t\udl{\varkappa} - \grad_{\pd_t\vxi_*}\udl{\varkappa} - \grad_{\vxi_*}\grad_{\vxi_*}\udl{\varkappa} +\alpha\scA(\ka)\udl{\varkappa} +\\
			& + \grad_{\vh_*}\grad_{\vh_*}\udl{\varkappa} + \grad_{\vH_*}\grad_{\vH_*}\udl{\varkappa} + \udl{(\grad_\vn q-\grad_\vn\wt{q})\cN\varkappa} + \udl{\fkR_2(\vv, \vh, \vH, \Gmt)}.
		\end{split}
	\end{equation}
	Note that \eqref{evo eqn kpudl} only holds for solutions to \eqref{PV problem}. For the sake of technical convenience, we would like to isolate the evolution equation for $\ka$ from the free boundary problem. One can consider the following ancillary vector field defined on $\Gmt$ (it is actually well-defined in the whole moving domain $\Omt$, but we presently use its restriction to the free boundary):
	\begin{equation}\label{def A}
		\vb{A}\coloneqq \Dt\vv - \grad_\vh\vh + \grad\qty[\alpha\cH\varkappa + q + \cH(\tfrac{1}{2}\abs{\vH}^2)],
	\end{equation}
	where $q$ is the multiplier-type pressure defined by \eqref{def q}. It is then clear that $\vb{A} \equiv \vb{0}$ whenever $(\vv, \vh, \vH, \Gmt)$ is a solution to \eqref{PV problem}.
	
	Thus, by plugging \eqref{def A} into \eqref{Dt2 kappa} and pulling back to $\Gms$, one obtains that
	\begin{equation}\label{pdt pdt kpudl}
		\begin{split}
			\pd_t\pd_t\udl{\varkappa} =\, &- 2\grad_{\vxi_*}\pd_t\udl{\varkappa} - \grad_{\pd_t\vxi_*}\udl{\varkappa} - \grad_{\vxi_*}\grad_{\vxi_*}\udl{\varkappa} +\alpha\scA(\ka)\udl{\varkappa} + \grad_{\vh_*}\grad_{\vh_*}\udl{\varkappa} + \grad_{\vH_*}\grad_{\vH_*}\udl{\varkappa} +\\
			&+\qty[\udl{-\grad_\vtau\grad_\vtau(\vb{A}\vdot\vn) + \grad_\vtau\varkappa (\vb{A}\vdot\vtau) - \varkappa^2 (\vb{A}\vdot\vn)}] +\udl{\fkR_3(\vv, \vh, \vH, \Gmt)},
		\end{split}
	\end{equation}
	where $\fkR_3$ is given by
	\begin{equation}\label{def R3}
		\begin{split}
			\fkR_3 =\, & -\alpha\abs{\grad_\vtau\varkappa}^2 + \alpha\varkappa^2\cN\varkappa + 2\grad_\vtau(\abs{\vh}^2)\grad_\vtau\varkappa + \varkappa^3\abs{\vh}^2 +2\varkappa\grad_\vtau\grad_\vtau(\abs{\vh}^2) + \\
			&+ \grad_\vtau\grad_\vtau\grad_\vn q + \varkappa^2\grad_\vn q + \grad_\vtau(\abs{\vH}^2)\grad_\vtau\varkappa + \tfrac{1}{2}\varkappa^2\cN(\abs{\vH}^2)+ \varkappa\grad_\vtau\grad_\vtau(\abs{\vH}^2) + \\
			& + \grad_\vtau\grad_\vtau\grad_\vn\wt{q} + \tfrac{1}{2}\grad_\vtau\grad_\vtau(\cN-\wt{\cN})(\abs{\vH}^2) +  6\varkappa\abs{\grad_\vtau\vv\vdot\vtau}^2 - 3\varkappa\abs{\grad_\vtau\vv\vdot\vn}^2 + \\
			& + 2 (\grad_\vtau\vv \vdot\vn)(\vtau\vdot\grad_\vtau\grad_\vtau\vv) + 4(\grad_\vtau\vv\vdot\vtau)(\vn\vdot\grad_\vtau\grad_\vtau\vv).
		\end{split}
	\end{equation}
	Furthermore, \eqref{var vxi*} implies that \eqref{pdt pdt kpudl} can be rewritten as
	\begin{equation*}
		\begin{split}
			\pd_t\pd_t\udl{\varkappa} =\, &- (\grad_{\vtau_*}\udl{\varkappa})\vtau_*\vdot \qty[\vb{B}(\ka)\pd_t\pd_t\ka + \vb{R}(\ka)\pd_t\omega_* + \vb{C}(\ka, \pd_t\ka, \omega_*)\pd_t\ka] + \\
			&- 2\grad_{\vxi_*}\pd_t\udl{\varkappa} - \grad_{\vxi_*}\grad_{\vxi_*}\udl{\varkappa} +\alpha\scA(\ka)\udl{\varkappa} + \grad_{\vh_*}\grad_{\vh_*}\udl{\varkappa} + \grad_{\vH_*}\grad_{\vH_*}\udl{\varkappa} + \\
			&+\qty[\udl{-\grad_\vtau\grad_\vtau(\vb{A}\vdot\vn) + \grad_\vtau\varkappa (\vb{A}\vdot\vtau) - \varkappa^2 (\vb{A}\vdot\vn)}] +\udl{\fkR_3(\vv, \vh, \vH, \Gmt)},
		\end{split}
	\end{equation*}
	where $\vtau_*$ is the unit tangent of $\Gms$.
	Next, recall that $\ka$ is defined to be $\udl{\varkappa} + a^2\varphi$, which yields that
	\begin{equation}\label{dtdt ka pre}
		\begin{split}
			\pd_t\pd_t \ka =\, &-\qty[(\grad_{\vtau_*}\udl{\varkappa})\vtau_*\vdot\vb{B}(\ka)]\pd_t\pd_t\ka -2\grad_{\vxi_*}\pd_t\ka - \grad_{\vxi_*}\grad_{\vxi_*}\ka + \alpha\scA(\ka)\ka + \\
			&+ \grad_{\vh_*}\grad_{\vh_*}\ka +\grad_{\vH_*}\grad_{\vH_*}\ka - (\grad_{\vtau_*}\udl{\varkappa})\vtau_*\vdot\qty[\vb{R}(\ka)\pd_t\omega_* + \vb{C}(\ka, \pd_t\ka, \omega_*)\pd_t\ka]  + \\
			&+ a^2\qty[\pd_t\pd_t\varphi +  2\grad_{\vxi_*}\pd_t\varphi + \grad_{\vxi_*}\grad_{\vxi_*}\varphi - \alpha\scA(\ka)\varphi - \grad_{\vh_*}\grad_{\vh_*}\varphi-\grad_{\vH_*}\grad_{\vH_*}\varphi]  + \\
			&+\udl{\fkR_3(\vv, \vh, \vH, \Gmt)} +\qty[\udl{-\grad_\vtau\grad_\vtau(\vb{A}\vdot\vn) + \grad_\vtau\varkappa (\vb{A}\vdot\vtau) - \varkappa^2 (\vb{A}\vdot\vn)}].
		\end{split}
	\end{equation}
	Although \eqref{var varphi} implies that $\pd_t\pd_t\varphi$ involve the factor $\pd_t\pd_t\ka$, it can be alternatively calculated through the relation \eqref{bc 1}
	that
	\begin{equation*}
		\pd_t\pd_t\varphi \vn_*\vdot \udl{\vn} + \pd_t\varphi\vn_*\vdot\pd_t\udl{\vn} = \pd_t(\udl{\vv\vdot\vn}).
	\end{equation*}
	On the other hand, \eqref{pull back dt} and \eqref{dt vn} yield that
	\begin{equation*}
		\begin{split}
			\pd_t\udl{\vv\vdot\vn} &= \udl{\Dt(\vv\vdot\vn)} - \grad_{\vxi_*}(\udl{\vv\vdot\vn}) \\
			&= \udl{\Dt\vv\vdot\vn} + \udl{\vv\vdot\Dt\vn} - \grad_{\vxi_*}\qty(\pd_t\varphi\vn_*\vdot\udl\vn) \\
			&= \udl{\vb{A}\vdot\vn} + \qty[\udl{\grad_\vh\vh\vdot\vn -\alpha\cN\varkappa - \grad_\vn q - \tfrac{1}{2}\cN(\abs{\vH}^2)}] - \udl{(\vv\vdot\vtau)(\grad_\vtau\vv\vdot\vn)} - \grad_{\vxi_*}\qty(\pd_t\varphi\vn_*\vdot\udl\vn),
		\end{split}
	\end{equation*}
	which implies
	\begin{equation*}
		\begin{split}
			\pd_t\pd_t\varphi \vn_*\vdot \udl{\vn} =\,  & \udl{\vb{A}\vdot\vn} +\qty[\udl{\grad_\vh\vh\vdot\vn -\alpha\cN\varkappa - \grad_\vn q - \tfrac{1}{2}\cN(\abs{\vH}^2)}] - \udl{(\vv\vdot\vtau)(\grad_\vtau\vv\vdot\vn)} \, + \\
			& - \grad_{\vxi_*}\qty(\pd_t\varphi\vn_*\vdot\udl\vn) + \pd_t\varphi\vn_*\vdot\qty[\udl{(\grad_\vtau\vv\vdot\vn)\vtau} + \grad_{\vxi_*}\udl\vn] \\
			=\, &\udl{\vb{A}\vdot\vn} - \qty[\udl{\abs{\vh}^2\varkappa +\alpha\cN\varkappa + \grad_\vn q + \tfrac{1}{2}\cN(\abs{\vH}^2)}]\, + \\
			&-\udl\vn\vdot\grad_{\vxi_*}(\pd_t\varphi\vn_*) + (\udl{\grad_\vtau\vv\vdot\vn})(\pd_t\varphi\vn_* - \udl\vv)\vdot\udl\vtau.
		\end{split}
	\end{equation*}
	Namely, one obtains the relation
	\begin{equation}\label{dtdt varphi}
		\begin{split}
			\pd_t\pd_t\varphi =\, &(\udl\vn\vdot\vn_*)^{-1}(\udl{\vb{A}\vdot\vn}) -(\udl\vn\vdot\vn_*)^{-1} \qty[\udl{\abs{\vh}^2\varkappa +\alpha\cN\varkappa + \grad_\vn q + \tfrac{1}{2}\cN(\abs{\vH}^2)}] \,+ \\
			 &-(\udl\vn\vdot\vn_*)^{-1}\qty[\udl\vn\vdot\grad_{\vxi_*}(\pd_t\varphi\vn_*) - (\udl{\grad_\vtau\vv\vdot\vn})(\pd_t\varphi\vn_* - \udl\vv)\vdot\udl\vtau].
		\end{split}
	\end{equation}
	In particular, by utilizing the ancillary vector field $\vb{A}$, the term $\pd_t\pd_t\varphi$ does not involve the $\pd_t\pd_t\ka$ factor.
	Combining \eqref{dtdt ka pre}-\eqref{dtdt varphi}, one arrives at
	\begin{equation}\label{dtdt ka pre2}
		\begin{split}
			&\hspace*{-2em}\qty[\textup{Id} + (\grad_{\vtau_*}\udl{\varkappa})\vtau_*\vdot\vb{B}(\ka)]\pd_t\pd_t\ka + 2\grad_{\vxi_*}\pd_t\ka + \grad_{\vxi_*}\grad_{\vxi_*}\ka  \\
			=\,&\alpha\scA(\ka)\ka + \grad_{\vh_*}\grad_{\vh_*}\ka + \grad_{\vH_*}\grad_{\vH_*}\ka + \\ &-a^2(\udl\vn\vdot\vn_*)^{-1}\qty{\qty[\udl{\abs{\vh}^2\varkappa +\alpha\cN\varkappa + \grad_\vn q + \tfrac{1}{2}\cN(\abs{\vH}^2)}] + \udl\vn\vdot\grad_{\vxi_*}(\pd_t\varphi\vn_*)} + \\
			&+ a^2(\udl\vn\vdot\vn_*)^{-1}\qty[(\udl{\grad_\vtau\vv\vdot\vn})(\pd_t\varphi\vn_* - \udl\vv)\vdot\udl\vtau] + \\
			& + a^2\qty[2\grad_{\vxi_*}\pd_t\varphi + \grad_{\vxi_*}\grad_{\vxi_*}\varphi - \alpha\scA(\ka)\varphi - \grad_{\vh_*}\grad_{\vh_*}\varphi - \grad_{\vH_*}\grad_{\vH_*}\varphi] + \\
			& - (\grad_{\vtau_*}\udl{\varkappa})\vtau_*\vdot\qty[\vb{R}(\ka)\pd_t\omega_* + \vb{C}(\ka, \pd_t\ka, \omega_*)\pd_t\ka] + \udl{\fkR_3(\vv, \vh, \vH, \Gmt)} +  \\
			&+ \qty[\udl{-\grad_\vtau\grad_\vtau(\vb{A}\vdot\vn) + \grad_\vtau\varkappa (\vb{A}\vdot\vtau) - \varkappa^2 (\vb{A}\vdot\vn)} +  a^2(\udl\vn\vdot\vn_*)^{-1}(\udl{\vb{A}\vdot\vn})].
		\end{split}
	\end{equation}
	To normalize the coefficient of $\pd_t\pd_t\ka$, one can consider the following linear operator defined on (some subspace of) $L^2(\Gms)$:
	\begin{equation}
		\scB(\ka) \coloneqq (\grad_{\vtau_*}\udl{\varkappa})\vtau_*\vdot\vb{B}(\ka),
	\end{equation}
	where $\varkappa$ is the curvature of the curve induced by $\ka$. In particular, Proposition \ref{prop op BRC} and product estimates imply that, for $\frac{1}{2} \le s' \le s-\frac{5}{2}$ (note here that $s \ge 3$), there holds
	\begin{equation}
		\abs{\scB(\ka)\psi}_{H^{s'}(\Gms)} \lesssim_{\epsilon, \Lambda_s} a^{\epsilon-2} \abs{\psi}_{H^{s'}(\Gms)} \abs{\ka}_{H^{s-\frac{3}{2}}(\Gms)} ,
	\end{equation}
	where $a$ is the ancillary constant in the definition of $\ka$, $0 < \epsilon \ll 1$ is a constant parameter, and one can simply take $\epsilon = 0$ when $s > 3$.
	
	Observe from the definition of $\Lambda_s = \Lambda(\Gms, s-\frac{1}{2}, \delta_0)$ that
	\begin{equation*}
		\abs{\ka}_{H^{s-\frac{3}{2}}(\Gms)} \lesssim_{\Lambda_s} \abs{\udl{\varkappa}}_{H^{s-\frac{3}{2}}(\Gms)} + a^2\abs{\varphi}_{H^{s-\frac{3}{2}}(\Gms)} \lesssim_{\Lambda_s} \abs{\varkappa}_{H^{s-\frac{3}{2}}(\Gmt)} + a^2\delta_0. 
	\end{equation*}
	Then,  by taking $a$ suitably large and $\delta_0 \ll_a 1$, one can obtain
	\begin{equation}
		\abs{\scB(\ka)}_{\scL\colon H^{s'}(\Gms) \to H^{s'}(\Gms)} \le \frac{1}{2} < 1
	\end{equation}
	for all $\ka \in B_{\delta_1} \cap H^{s-\frac{3}{2}}(\Gms)$ with the notations in Lemma \ref{lem ka}.	In other words, the operator $[\textup{Id}+\scB(\ka)] $ is an isomorphism in $H^{s'}(\Gms)$. 
	
	Recall that $\vv$ is uniquely determined by the boundary evolution and the preassigned ``vorticity'' $\omega_*$. In a similar manner, we denote by
	\begin{equation}
		\jmath_* \coloneqq (\Curl\vh)\circ\cX.
	\end{equation}
	Then, $\vh$ can be recovered from $(\ka, \jmath_*)$ through solving the div-curl problem
	\begin{equation}
		\begin{cases*}
			\dive\vh = 0 &in $\Omt$, \\
			\Curl \vh = \jmath_*\circ\cX_{\Gmt}^{-1} &in $\Omt$, \\
			\vh\vdot\vn = 0 &on $\Gmt$.
		\end{cases*}
	\end{equation}
	It is obvious that the vector field $\vH\colon\cV_t\to \R^2$ is uniquely determined by $\ka$ and $\cJ$. Thus, the remainder term $\fkR_3(\vv, \vh, \vH, \Gmt)$  can be equivalently expressed as $\widetilde{\fkR}_3(\ka, \pd_t\ka, \omega_*, \jmath_*, \cJ)$. We shall use notations $(\ka, \pd_t\ka, \omega_*, \jmath_*, \cJ)$ instead of $(\vv, \vh, \vH, \Gmt)$ in the following arguments.
	
	Applying the operator $[\textup{Id}+\scB(\ka)]^{-1} $ to both sides of \eqref{dtdt ka pre2} and reorganizing the notations, one will get
	\begin{equation}\label{dtdt ka}
		\begin{split}
			&\hspace*{-1.5em}\pd_t\pd_t\ka + 2\grad_{\vxi_*}\pd_t\ka + \grad_{\vxi_*}\grad_{\vxi_*}\ka - \alpha\scA(\ka)\ka - \grad_{\vh_*}\grad_{\vh_*}\ka - \grad_{\vH_*}\grad_{\vH_*}\ka \\
			=\, &\scF(\ka)\pd_t\omega_* + \scG(\ka, \pd_t\ka, \omega_*, \jmath_*, \cJ) + \\
			&+ [\textup{Id}+\scB(\ka)]^{-1}\qty[\udl{-\grad_\vtau\grad_\vtau(\vb{A}\vdot\vn) + \grad_\vtau\varkappa (\vb{A}\vdot\vtau) - \varkappa^2 (\vb{A}\vdot\vn)} +  a^2(\udl\vn\vdot\vn_*)^{-1}(\udl{\vb{A}\vdot\vn})].
		\end{split}
	\end{equation}
	This is the desired second order evolution equation for $\ka$, which is, in particular, independent of the free boundary problem \eqref{PV problem}. The relation between \eqref{dtdt ka} and solutions to \eqref{PV problem} is reflected by the ancillary vector field $\vb{A}$.
	
	The term $\scF(\ka)\pd_t\omega_* $ in \eqref{dtdt ka} is defined as
	\begin{equation}
		\scF(\ka)\pd_t\omega_* \coloneqq -\qty{[\textup{Id}+\scB(\ka)]^{-1}(\grad_{\vtau_*}\udl{\varkappa})\vtau_*\vdot\vb{R}(\ka)}\pd_t\omega_*,
	\end{equation}
	which indicates that $\scF(\ka)$ is a linear operator. Moreover, $\scF(\ka)$ satisfies the following estimates for $\frac{1}{2} \le s' \le s-\frac{5}{2}$:
	\begin{subequations}\label{est scF}
		\begin{equation}
			\abs{\scF(\ka)}_{\scL\colon H^{s'-\frac{1}{2}}(\Oms) \to H^{s'}(\Gms)} \lesssim_{\epsilon, \Lambda_s} \abs{\ka}_{H^{\frac{3}{2}+\epsilon}(\Gms)} + \abs{\ka}_{H^{s-\frac{3}{2}}(\Gms)}
		\end{equation}
		and
		\begin{equation}
			\abs{\scF(\ka)}_{\scL\colon H^{s'-\frac{1}{2}+\epsilon}(\Oms) \to H^{s'}(\Gms)} \lesssim_{\epsilon, \Lambda_s} \abs{\ka}_{H^{s-\frac{3}{2}}(\Gms)},
		\end{equation}
	\end{subequations}
	where $0 < \epsilon \ll 1$ is an arbitrarily small parameter and one can simply take $\epsilon = 0$ if $s > 3$.
	
	The term $\scG(\ka, \pd_t\ka, \omega_*, \jmath_*, \cJ)$ is the collection of all those miscellaneous terms. It follows from tedious but routine calculations that (see also \cite[Appendix. B]{Liu-Xin2023}), for $\frac{1}{2} \le s' \le s-\frac{5}{2}$, there holds
	\begin{equation}\label{est scG}
		\begin{split}
			&\hspace*{-2em}\abs{\scG(\ka, \pd_t\ka, \omega_*, \jmath_*, \cJ)}_{H^{s'}(\Gms)} \\
			 \le & \cQ_1\qty(\alpha\abs{\ka}_{H^{s'+ \frac{3}{2}}(\Gms)}, \abs{\ka}_{H^{s'+1}(\Gms)}, \abs{\pd_t\ka}_{H^{s'}(\Gms)}, \norm{(\omega_*, \jmath_*)}_{H^{s'+\frac{3}{2}}(\Oms)}, \abs{\cJ}_{H^{s'+2}(\cS)}) + \\
			&+ a^2\cQ_2\qty(\alpha\abs{\ka}_{H^{s'+1}(\Gms)}, \abs{\pd_t\ka}_{H^{s'-1}(\Gms)}, \norm{(\omega_*, \jmath_*)}_{H^{s'+\frac{1}{2}}(\Oms)}, \abs{\cJ}_{H^{s'+1}(\cS)}),
		\end{split}
	\end{equation}
	where $\cQ_1$ and $\cQ_2$ are both generic polynomials determined by $\Lambda_s$.
	
	\subsection{Iterative Process}\label{sec it}
	From now on, we always assume that one of the following conditions hold
	\begin{enumerate}
		\item $\alpha > 0$, i.e., there exists surface tension;
		\item $\alpha = 0$, but there is a constant $\lambda_0$ so that
		\begin{equation} \label{def lambda0}
			\abs{\vh} + \abs{\vH} \ge \lambda_0 > 0 \qq{on} \Gmt,
		\end{equation}
		i.e., there is no surface tension, but the total magnetic fields are non-degenerate on the free boundary.
	\end{enumerate}
	
	\subsubsection{Initial settings}
	For given positive constants $T, M_1, M_2$, and $M_3$, we define $\fkX = \fkX(\ka, \omega_*, \jmath_*)$ to be the collection of data in the spaces
	\begin{subequations}
		\begin{equation}\label{range ka}
			\ka \in C^1\qty([0, T]\to B_{\delta_1} \subset H^{s-\frac{5}{2}}(\Gms)),
		\end{equation}
		\begin{equation}
			\ka \in \begin{cases*}
				C^0\qty([0, T]\to H^{s-1}(\Gms)) \cap C^2\qty([0, T]\to H^{s-4}(\Gms)) &if $\alpha > 0$, \\
				C^0\qty([0, T]\to H^{s-\frac{3}{2}}(\Gms)) \cap C^2\qty([0, T]\to H^{s-\frac{7}{2}}(\Gms)) &if $\alpha = 0$,
			\end{cases*}
		\end{equation}
		and
		\begin{equation}
			\omega_*, \jmath_* \in C^0\qty([0, T]\to H^{s-1}(\Oms)) \cap C^1\qty([0, T] \to H^{s-2}(\Oms)).
		\end{equation}
	\end{subequations}
	Moreover, they are assumed to satisfy the bounds
	\begin{subequations}
		\begin{equation}
			\abs{\ka_{\restriction_{t=0}} - \varkappa_*}_{H^{s-\frac{5}{2}}(\Gms)} \le \frac{1}{2}\delta_1,
		\end{equation}
		\begin{equation}\label{def M1}
			\sup_{t\in[0, T]}\qty(\alpha\abs{\ka}_{H^{s-1}(\Gms)}^2, \abs{\ka}_{H^{s-\frac{3}{2}}(\Gms)}^2, \abs{\pd_t\ka}_{H^{s-\frac{5}{2}}(\Gms)}^2 , \norm{(\omega_*)}_{H^{s-1}(\Oms)}^2, \norm{(\jmath_*)}_{H^{s-1}(\Oms)}^2) \le M_1,
		\end{equation}
		\begin{equation}
			\sup_{t\in[0, T]}\qty(\norm{\pd_t\omega_*}_{H^{s-2}(\Oms)}^2, \norm{\pd_t\jmath_*}_{H^{s-2}(\Oms)}^2) \le M_2,
		\end{equation}
		and
		\begin{equation}
			a^2 M_3 \ge \sup_{t\in[0, T]} 
			\begin{cases*}
				\abs{\pd_t\pd_t\ka}_{H^{s-4}(\Gms)} & if $\alpha > 0$, \\
				\abs{\pd_t\pd_t\ka}_{H^{s-\frac{7}{2}}(\Gms)} & if $\alpha = 0$,
			\end{cases*}
		\end{equation}
	\end{subequations}
	where ``$a$'' is the ancillary constant in the definition of $\ka$.
	
	For the collection of initial data, we take two constants $0 < M_0 < M_1$ and $0 < \delta_2 \ll \delta_1$ and define $\fkI = \fkI(\widehat{\ka}, \widehat{\pd_t\ka}, \widehat{\omega_*}, \widehat{\jmath_*}) $ to be the collection of initial data satisfying the bounds
	\begin{subequations}
		\begin{equation}
		 \alpha\abs{\widehat{\ka} - \varkappa_*}_{H^{s-\frac{1}{2}}(\Gms)}^2 +	\abs{\widehat{\ka}-\varkappa_*}_{H^{s-\frac{3}{2}}(\Gms)}^2 \le \abs{\delta_2}^2
		\end{equation}
		and
		\begin{equation}
			\abs{\widehat{\pd_t\ka}}_{H^{s-\frac{5}{2}}(\Gms)}^2 + \norm{\widehat{\omega_*}}_{H^{s-1}(\Oms)}^2 + \norm{\widehat{\jmath_*}}_{H^{s-1}(\Oms)}^2 \le M_0.
		\end{equation}
	\end{subequations}
	
	Here we also assume that these bound constants are taken large enough so that
	\begin{equation}
		\sup_{t\in [0, T]} \abs{\cJ}_{H^{s-\frac{1}{2}}(\cS)}^2 \le M_1 \qc \sup_{t\in [0, T]}\abs{\pd_t \cJ}_{H^{s-\frac{3}{2}}(\cS)}^2 \le M_2, \qand \abs{\cJ_{\restriction_{t=0}}}_{H^{s-\frac{1}{2}}(\cS)}^2 \le M_0.
	\end{equation}
	
	\subsubsection{Recovery of domains, velocities, and magnetic fields}\label{sec recov}
	For a given family of states $\{(\ka, \omega_*, \jmath_*)(t)\}_t \subset \fkX$, one needs to re-define the moving domains and vector fields. First, one can notice that the family $\{\ka(t)\}$ induces a family of moving curves $\{\Gmt\} \subset \Lambda_s$ through the isomorphisms given by \eqref{def ka} and \eqref{def vPhi}. In particular, one has well-defined diffeomorphisms $\vPhi_{\Gmt}$ and $\cX_{\Gmt}$.
	
	The vector field $\vv$ defined in $\Omt$ is obtained through solving the following div-curl problem:
	\begin{equation}\label{def v}
		\begin{cases*}
			\dive\vv = \gamma &in $\Omt$, \\
			\Curl \vv =  \omega_*\circ\cX^{-1}_{\Gmt} &in $\Omt$, \\
			\vv\vdot\vn = \vn \vdot (\pd_t\varphi_{\Gmt}\vn_*) \circ \vPhi_{\Gmt}^{-1} &on $\Gmt$,
		\end{cases*}
	\end{equation}
	where $\gamma$ is the balancing constant to ensure the compatibility condition
	\begin{equation*}
		\gamma\int_{\Omt} \dd{x} = \int_{\Gmt} \vn \vdot (\pd_t\varphi_{\Gmt}\vn_*) \dd{\ell}.
	\end{equation*}
	Namely, for each fixed moment $t$, $\gamma$ is a constant given by
	\begin{equation}\label{def gm}
		\gamma \coloneqq \frac{1}{\abs{\Omt}} \int_{\Gmt} \vn \vdot (\pd_t\varphi_{\Gmt}\vn_*)\circ \vPhi_{\Gmt}^{-1} \dd{\ell},
	\end{equation}
	where $\abs{\Omt}$ is the volume (area) of the moving domain $\Omt$. We remark here that the simply-connectedness of $\Omt$ ensure the unique solvability of the div-curl problem.
	Similarly, the vector field $\vh$ is defined to be the unique solution to
	\begin{equation}\label{def h}
		\begin{cases*}
			\dive\vh = 0 &in $\Omt$, \\
			\Curl \vh = \jmath_* \circ \cX_{\Gmt}^{-1} &in $\Omt$, \\
			\vh \vdot \vn = 0 &on $\Gmt$.
		\end{cases*}
	\end{equation}
	The vector field $\vH$ is given by
	\begin{equation}
		\begin{cases*}
			\dive\vH = 0 &in $\cV_t$, \\
			\Curl\vH = 0 &in $\cV_t$, \\
			\vH\vdot\vn = 0 &on $\Gmt$, \\
			\vN\smallwedge\vH = \cJ &on $\cS$.
		\end{cases*}
	\end{equation}
	
	For the initial data, one clearly has
	\begin{equation}
		\qty(\wh{\ka}, \wh{\pd_t\ka}, \wh{\omega_*}, \wh{\jmath_*}) \leadstoext\leadsto \qty(\wh{\vv}, \wh{\vh}, \wh{\vH}, \wh{\Omega}).
	\end{equation}
	If $\alpha = 0$, we further require that the ``magnetic fields'' $\wh{\vh}$ and $\wh{\vH}$ satisfy
	\begin{equation}\label{req ini h}
		\abs*{\wh{\vh}} + \abs{\wh{\vH}} > 2\lambda_0 \qq{on} \wh{\Gamma},
	\end{equation}
	where $\lambda_0$ is the generic constant appeared in \eqref{def lambda0}.
	
	\subsubsection{Definition of iterates}
	Suppose now that $\{(\kaz, \omz, \jz)(t)\}_t \subset \fkX$ is the given collection of data, and $(\wh{\ka}, \wh{\pd_t\ka}, \wh{\omega_*}, \wh{\jmath_*}) \in \fkI$ is the given initial data. In particular, through the recovery process introduced in \S\ref{sec recov}, one obtains the family of data
	\begin{equation*}
		(\kaz, \omz, \jz)(t) \leadstoext\leadsto \qty(\vz(t), \hz(t), \Hz(t), \Omt^0).
	\end{equation*}
	
	Then, we define the one step iterates via solving the linear evolution problems:
	\begin{subequations}\label{ite eqns}
		\begin{equation}\label{it bdry}
			\left\{\begin{split}
				&\pd_t\pd_t\vartheta + 2\grad_{\xiz}\pd_t\vartheta + \grad_{\xiz}\grad_{\xiz}\vartheta - \alpha\scA\qty(\kaz)\vartheta - \grad_{\hz_*}\grad_{\hz_*}\vartheta - \grad_{\Hz_*}\grad_{\Hz_*}\vartheta \\
				&\quad = \scF\qty(\kaz)\pd_t\omz + \scG\qty(\kaz, \pd_t\kaz, \omz, \jz, \cJ) \qq{on} \Gms, \\
				&\vartheta_{\restriction_{t=0}} = \wh{\ka} \qc (\pd_t\vartheta)_{\restriction_{t=0}} = \wh{\pd_t\ka},
			\end{split}
			\right.
		\end{equation}
		where $\xiz$, $\hz_*$, and $\Hz_*$ are defined through \eqref{def vxi*}, \eqref{def vh*}, and \eqref{def vH*}, respectively, and $\scA(\kaz)$ is the operator given by \eqref{def scA}.
		\begin{equation}\label{it curl}
			\left\{
			\begin{split}
				&\pd_t (\varpi_-) + (\vz+\hz)\vdot\grad\varpi_- + \cB\qty[\pd(\vz+\hz), \pd(\vz-\hz)]=0 \qin \Omt^0, \\
				&\pd_t (\varpi_+) + (\vz-\hz)\vdot\grad\varpi_+ + \cB\qty[\pd(\vz-\hz), \pd(\vz+\hz)]=0 \qin \Omt^0, \\
				&(\varpi_+)_{\restriction_{t=0}} = (\wh{\omega_*} + \wh{\jmath_*})\circ \cX_{\Gamma^0_{t=0}}^{-1} \qc (\varpi_-)_{\restriction_{t=0}} = (\wh{\omega_*} - \wh{\jmath_*})\circ \cX_{\Gamma^0_{t=0}}^{-1},
			\end{split}
			\right.
		\end{equation}
		where $\cB$ is the bilinear form defined in \eqref{def cB}.
	\end{subequations}
	
	The existence of solutions to the linear initial value problem \eqref{it bdry} can be deduced from the estimates given by Proposition \ref{prop lin est theta} and finite element methods or standard semi-group theory (cf. \cite[Proposition 4.1]{Shatah-Zeng2011}). Indeed, although the estimates in \S\ref{sec evo lin k} are derived on moving curves, they can be transmitted to the reference frame through the uniformly bounded coordinate maps. The solvability of  \eqref{it curl} can be obtained by using the characteristic (ODE) methods. The precise estimates for solutions to \eqref{ite eqns} will be deferred to \S\ref{sec bound ite}.
	
	Now, assuming that the problem \eqref{ite eqns} is solvable. We define the one-step iterates by:
	\begin{equation}
		\kao \coloneqq \vartheta \qc \omo \coloneqq \frac{1}{2}\qty(\varpi_+ + \varpi_-) \circ \cX_{\Gamma^0_t}, \qand \jo \coloneqq \frac{1}{2}\qty(\varpi_+ - \varpi_-) \circ \cX_{\Gamma^0_t}.
	\end{equation}
	In summary, we have defined the iteration map $\fkT$ through solving \eqref{ite eqns}:
	\begin{equation}\label{def fkT}
		\fkT\qty[\qty(\kaz, \omz, \jz); \qty(\wh{\ka}, \wh{\pd_t\ka}, \wh{\omega_*}, \wh{\jmath_*})] \coloneqq	\qty(\kao, \omo, \jo).
	\end{equation}
	
	\subsection{Solutions to the Plasma-Vacuum Problems}
	Assume now that $\fkT$ admits a fixed point (not necessarily being unique) in $\fkX$ for given initial data in $\fkI$. Then, we need to show that the fixed point $(\ka, \omega_*, \jmath_*) \in \fkX$ can actually induce a solution to the free boundary problem \eqref{PV problem}.
	
	In view of the linear boundary evolution problem \eqref{it bdry} and the equation \eqref{dtdt ka}, one obtains the trivial identity
	\begin{equation*}
		-\grad_\vtau\grad_\vtau(\vb{A}\vdot\vn) + \grad_\vtau\varkappa (\vb{A}\vdot\vtau) - \varkappa^2 (\vb{A}\vdot\vn) + \frac{a^2(\vb{A}\vdot\vn)}{(\vn_*\circ\vPhi^{-1}) \vdot\vn} \equiv 0 \qq{on} \Gmt.
	\end{equation*}
	The definition of $\Lambda_s$ implies that the denominator $(\vn_*\circ\vPhi^{-1}) \vdot\vn \approx 1 $. In other words, when $a$ is suitably large, one can derive the following standard elliptic estimate:
	\begin{equation}\label{est A}
		\abs{\grad_\vtau(\vb{A}\vdot\vn)}_{L^2(\Gmt)}^2 + \frac{a^2}{2}\abs{\vb{A}\vdot\vn}_{L^2(\Gmt)}^2 \lesssim_{\Lambda_s} \abs{(\vb{A}\vdot\vtau)\grad_\vtau\varkappa}_{L^2(\Gmt)}^2,
	\end{equation}
	where the implicit constant is independent of $a$.
	
	On the other hand, by taking curl of \eqref{def A}, it can be deduced from \eqref{it curl} and \eqref{evo eqn curl} that $\vb{A}$ is actually irrotational. Similarly, by taking divergence of \eqref{def A} and invoking \eqref{def v}-\eqref{def h}, one obtains that
	\begin{equation}\label{div A}
		\dive \vb{A} = \pd_t \gamma = \pd_t(\dive\vv),
	\end{equation}
	which is a constant for each fixed moment $t$. For the simplicity of notations, we denote by
	\begin{equation}
		\Theta \coloneqq \vb{A} \vdot \vn.
	\end{equation}
	Then, it is clear that $\vb{A}$ satisfies the following div-curl system:
	\begin{equation}
		\begin{cases*}
			\dive\vb{A} = \frac{1}{\abs{\Omt}}\int_{\Gmt}\Theta \dd{\ell} &in $\Omt$, \\
			\Curl \vb{A} = 0 &in $\Omt$, \\
			\vb{A} \vdot \vn = \Theta &on $\Gmt$.
		\end{cases*}
	\end{equation} 
	In particular, there holds
	\begin{equation}
		\abs{\vb{A}\vdot\vtau}_{H^{\frac{1}{2}}(\Gmt)} \lesssim_{\Lambda_s} \abs{\Theta}_{H^{\frac{1}{2}}(\Gmt)}.
	\end{equation}
	Therefore, \eqref{est A} can be refined to
	\begin{equation}
		\abs{\grad_\vtau\Theta}_{L^2(\Gmt)}^2 + \frac{a^2}{2}\abs{\Theta}_{L^2(\Gmt)}^2 \lesssim_{\Lambda_s} \abs{\Theta}_{H^{\frac{1}{2}}(\Gmt)}^2\abs{\grad_\vtau\varkappa}_{H^{-\frac{1}{2}}(\Gmt)}^2 \lesssim_{\Lambda_s} \abs{\Theta}_{H^{\frac{1}{2}}(\Gmt)}^2\abs{\varkappa}_{H^{\frac{1}{2}}(\Gmt)}^2
	\end{equation}
	The largeness of $a$, the hypothesis \eqref{range ka}, and standard interpolations yield that
	\begin{equation}
		\Theta \equiv 0,
	\end{equation}
	which implies that
	\begin{equation}
		\vb{A} \equiv \vb{0} \implies \Dt\vv + \grad\qty[\alpha\varkappa_{\cH} + q + \cH(\tfrac{1}{2}\abs{\vH}^2)] - \grad_\vh\vh \equiv \vb{0}.
	\end{equation}
	Particularly, \eqref{div A} and the fact that $(\dive\vv)_{\restriction_{t=0}}=0$ leads to
	\begin{equation}
		\dive\vv \equiv 0.
	\end{equation}
	It remains to check that the evolution equation for $\vh$ in \eqref{MHD} holds. Indeed, if one denotes by
	\begin{equation}
		\vb*{\Xi} \coloneqq \Dt\vh-\grad_\vh\vv,
	\end{equation}
	then it is routine to check that
	\begin{equation}
		\begin{cases*}
			\dive\vb*{\Xi} = 0 &in $\Omt$, \\
			\Curl \vb*{\Xi} = 0 &in $\Omt$, \\
			\vb*{\Xi} \vdot\vn =0 &on $\Gmt$,
		\end{cases*}
	\end{equation}
	for all $t$.
	
	The above arguments conclude that each fixed point of $\fkT$ given by \eqref{def fkT} indeed induces a solution $(\vv, \vh, \vH, \Omt)$ to the free boundary problem \eqref{PV problem}.

	\section{Estimates for the Iterations}\label{sec est ite}
	
	In this section, we shall show that the iterate $\fkT$ is indeed a map from $\fkX$ to $\fkX$, at least for a short time period $T$, and we will show that $\fkT$ is a contraction in a larger space containing $\fkX$, which indicates the existence of unique fixed point for each legitimate input data.
	
	\subsection{Bounds for the Iterations}\label{sec bound ite}
	
	With the notations in \S\ref{sec it}, we first show that $\kao$ satisfies the bound \eqref{def M1}. Indeed, it follows from \eqref{est scF}-\eqref{est scG}, \eqref{var vxi*}, Proposition \ref{prop op BRC}, linear estimates presented in \S\S\ref{sec est alpha>0}-\ref{sec est h non-degenerate}, the energy coercivity estimates \eqref{bd eng coerc a>0} \& \eqref{bd eng coerc a=0}, and Gronwall's inequality that (see also \cite[Proposition 5.2]{Liu-Xin2023})
	\begin{equation}
		\begin{split}
			&\sup_{t\in[0, T]}\qty(\alpha\abs{\kao}_{H^{s-1}(\Gms)}^2 + \abs{\kao}_{H^{s-\frac{3}{2}}(\Gms)}^2 + \abs{\pd_t\kao}_{H^{s-\frac{5}{2}}(\Gms)}^2) \\
			 &\quad\le \exp{T\cQ(M_1, M_2, a^2 M_3)} \cdot \qty[\cP{M_0} + T(a^2 + M_2)\cP{M_1}],
		\end{split}
	\end{equation} 
	where $\mathcal{P}$ and $\cQ $ are generic polynomials depending on $\Lambda_s$ and either $\alpha$ (if $\alpha > 0$) or $\lambda_0$ in \eqref{def lambda0} (if $\alpha = 0$).
	
	By taking $M_1 \gg M_0$, $T \ll 1$, and $M_3 \gg (M_1 + M_2) $, one obtains that 
	\begin{equation}
		\sup_{t\in[0, T]}\qty(\alpha\abs{\kao}_{H^{s-1}(\Gms)}^2, \abs{\kao}_{H^{s-\frac{3}{2}}(\Gms)}^2, \abs{\pd_t\kao}_{H^{s-\frac{5}{2}}(\Gms)}^2) \le M_1
	\end{equation}
	and
	\begin{equation}
			a^2 M_3 \ge \sup_{t\in[0, T]} 
			\begin{cases*}
				\abs{\pd_t\pd_t\kao}_{H^{s-4}(\Gms)} & if $\alpha > 0$, \\
				\abs{\pd_t\pd_t\kao}_{H^{s-\frac{7}{2}}(\Gms)} & if $\alpha = 0$,
			\end{cases*}
	\end{equation}
	which follows from the evolution equation \eqref{it bdry}. These calculations are tedious but indeed standard, so we omit them here. One can refer to \cite{Liu-Xin2023} for the similar computations.			
	
	The estimates for $\omo$ and $\jo$ are much simpler, which can be directly derived from \eqref{curl evo est}, together with uniform product and elliptic estimates. More precisely, one can routinely get
	\begin{equation}
		\sup_{t\in[0, T]} \qty(\norm{\omo}_{H^{s-1}(\Oms)}^2, \norm{\jo}_{H^{s-1}(\Oms)}^2) \le M_1
	\end{equation}
	and
	\begin{equation}
		\sup_{t\in[0, T]} \qty(\norm{\pd_t\omo}_{H^{s-2}(\Oms)}^2, \norm{\pd_t\jo}_{H^{s-2}(\Oms)}^2)  \le T\cQ(M_1) \le M_2,
	\end{equation}
	provided that $M_2 \gg M_1$ and $T \ll 1$.

	Moreover, if $\alpha = 0$, the Sobolev embeddings and smallness of $T$ yields that \eqref{def lambda0} holds for $\vh^1$ on $t \in [0, T]$.
	
	In summary, with an appropriate choice of constants $T, M_0, M_1, M_2$, and $M_3$, the iteration map $\fkT$ defined by \eqref{def fkT} satisfies
	\begin{equation}
		\fkT \colon \fkI(\delta_2, M_0) \times \fkX(T, M_1, M_2, M_3) \to \fkX(T, M_1, M_2, M_3).
	\end{equation} 
	
	In particular, topological fixed point theorem implies the existence of fixed points for $\fkT$ when $s \ge 3$, which implies the existence of solutions to \eqref{PV problem}.
	
	\subsection{Variational Estimates}\label{sec var est}
	
	To show the uniqueness of fixed points, it suffices to check that the iteration map $\fkT$ is indeed a contraction in some larger spaces containing $\fkX$. For the sake of convenience, we assume that $\{(\wh{\ka}, \wh{\pd_t\ka}, \wh{\omega_*}, \wh{\jmath_*})(\beta)\}_\beta \subset \fkI(\delta_2, M_0)$ and $\{(\kaz, \omz, \jz)(\beta)\}_{\beta} \subset \fkX(T; M_1, M_2, M_3)$ are two families of data parameterized by $\beta$. For each fixed parameter $\beta$, we denote by $(\kao, \omo, \jo)(\beta)$ the output of iterates $\fkT$ given by \eqref{def fkT} with input $(\kaz, \omz, \jz)(\beta)$ and $(\wh{\ka}, \wh{\pd_t\ka}, \wh{\omega_*}, \wh{\jmath_*})(\beta)$. The main goal of \S\ref{sec var est} is to establish the estimates for variational derivatives of these data with respect to the parameter $\beta$.
	
	Indeed, applying $\pdv*{\beta}$ to the linear problem \eqref{it bdry} yields
	\begin{equation}\label{var it eqn ka}
		\left\{
		\begin{split}
			&\pd_t\pd_t\pd_\beta\kao + 2 \grad_{\xiz}\pd_t\pd_\beta\kao + \grad_{\xiz}\grad_{\xiz}\pd_\beta\kao - \alpha\scA(\kaz)\pd_\beta\kao - \grad_{\hz_*}\grad_{\hz_*}\pd_\beta\kao - \grad_{\Hz_*}\grad_{\Hz_*}\pd_\beta\kao\\
			&\quad = - 2\grad_{\pd_\beta\xiz}\pd_t\kao - \grad_{\pd_\beta\xiz}\grad_{\xiz}\kao - \grad_{\xiz}\grad_{\pd_\beta\xiz}\kao  - \alpha(\fdv*{\scA}{\ka})_{\restriction_{\kaz}}[\pd_\beta\kaz; \kao] + \\
			&\qquad \quad
			 + \grad_{\pd_\beta\hz_*}\grad_{\hz_*}\kao + \grad_{\hz_*}\grad_{\pd_\beta\hz_*}\kao + \grad_{\pd_\beta\Hz_*}\grad_{\Hz_*}\kao + \grad_{\Hz_*}\grad_{\pd_\beta\Hz_*}\kao + \scF(\kaz)\pd_t\pd_\beta\omz +  \\
			&\qquad \qquad + (\fdv*{\scF}{\ka})_{\restriction_{\kaz}}[\pd_\beta\kaz; \pd_t\omz] + \pd_{\beta}\qty[\scG(\kaz, \pd_t\kaz, \omz, \jz, \cJ)] \qq{on} \Gms, \\
			&(\pd_\beta\kao)_{\restriction_{t=0}} = \pd_\beta\wh{\ka} \qc (\pd_t\pd_\beta\kao)_{\restriction_{t=0}} = \pd_\beta\wh{\pd_t\ka}.
		\end{split}
		\right.
	\end{equation}
	In order to apply Proposition \ref{prop lin est theta}, one needs  precise estimates for the variational terms. Indeed, with the same notations as in \S\ref{sec var velocity}, we define
	\begin{equation}
		\vbz \coloneqq \qty(\pd_\beta\cX_{\Gamma^0_{t, \beta}}) \circ \cX_{\Gamma^0_{t, \beta}}^{-1} \qand \Dbt \coloneqq \pdv{\beta} + \grad_\vbz.
	\end{equation}
	Then, it is clear that
	\begin{equation}
		(\pd_\beta \psi) \circ \vPhi_{\Gamma^0_{t, \beta}}^{-1} = \Dbt\qty(\psi \circ \vPhi_{\Gamma^0_{t, \beta}}^{-1}) \qfor \psi \colon \Gms \to \R,
	\end{equation}
	and
	\begin{equation}
		(\pd_\beta f) \circ \cX_{\Gamma^0_{t, \beta}}^{-1} = \Dbt\qty(f \circ \cX_{\Gamma^0_{t, \beta}}^{-1}) \qfor f \colon \Oms \to \R.
	\end{equation}
	For the simplicity of notations, we denote by
	\begin{equation}
		\overline{\psi} \coloneqq \psi \circ \vPhi_{\Gamma^0_{t, \beta}}^{-1}\ \text{  for  }\ \psi \colon \Gms \to \R,
	\qand
		\overline{f} \coloneqq f \circ \cX_{\Gamma^0_{t, \beta}}^{-1}\  \text{  for  }\ f \colon \Oms \to \R.
	\end{equation}
	Namely, $\overline{\psi}$ and $\overline{f}$ are respectively the push-out functions defined on the moving boundary and domain through the coordinate maps. Thus, the definition \eqref{def scA} implies that
	\begin{equation}
		\overline{\pdv{\beta}\qty[\scA(\kaz)\kao]} = \Dbt\qty(\grad_{\vtau}\grad_\vtau\cN\overline{\kao}),
	\end{equation}
	where $\vtau$ is the unit tangent of $\Gamma^0_{t, \beta}$ and $\cN$ is the Dirichlet-Neumann operator on it. In other words, one obtains that
	\begin{equation}
		\pdv{\beta}\qty[\scA(\kaz)\kao]-\scA(\kaz)\pd_\beta\kao = \qty(\comm{\Dbt}{\grad_\vtau\grad_\vtau\cN}\overline{\kao})\circ\vPhi_{\Gamma^0_{t, \beta}} \eqqcolon \qty(\fdv{\scA}{\ka})_{\restriction_{\kaz}}\!\qty[\pd_\beta\kaz; \kao].
	\end{equation}
	Commutator formulae \eqref{comm Dt dtau}-\eqref{comm formula dt N} lead to
	\begin{equation}
		\abs{\comm{\Dbt}{\grad_\vtau\grad_\vtau\cN}\overline{\kao}}_{H^{s'}(\Gamma^0_{t, \beta})} \lesssim_{\Lambda_s} \cP{\abs{\Gamma_{t, \beta}^0}_{H^{s'+3}}}\abs{\vbz}_{H^{s'+3}(\Gamma_{t, \beta}^0)}\abs{\overline{\kao}}_{H^{s'+3}(\Gamma_{t, \beta}^0)} \qfor s'\ge - \frac{1}{2}.
	\end{equation}
	Note that the linear estimates in Proposition \ref{prop lin est theta} are valid for the source term $f \in H^{\ge\frac{1}{2}}(\Gmt)$, so we need to further require that
	\begin{equation}\label{req s ST}
		s-1 \ge s' + 3 \ge \frac{1}{2} + 3 \implies s \ge \frac{9}{2} \qq{when} \alpha > 0.
	\end{equation}
	
	Similarly, commutator, product, and elliptic estimates, together with quite tedious calculations, yield that (see also \cite[Appendix B]{Liu-Xin2023}), for $\frac{1}{2} \le\sigma \le s-\frac{7}{2}$ 
	\begin{equation}
		\abs{\qty(\fdv{\scF}{\ka})_{\restriction_{\kaz}}\!\![\pd_\beta\kaz; \pd_t\omz]}_{H^{\sigma}(\Gms)} \hspace*{-1em} \lesssim_{\Lambda_s} \abs{\pd_\beta\kaz}_{H^{\sigma+1}(\Gms)}\norm{\pd_t\omz}_{H^{s-2}(\Oms)}
	\end{equation}
	and
	\begin{equation}
		\begin{split}
			&\abs{\pdv{\beta}\qty[\scG(\kaz, \pd_t\kaz, \omz, \jz, \cJ)]}_{H^{\sigma}(\Gms)} \\ &\quad \le \cQ\qty(\abs{\pd_t\kaz}_{H^{s-\frac{5}{2}}(\Gms)}, \norm{\omz}_{H^{s-1}(\Oms)}, \norm{\jz}_{H^{s-1}(\Oms)}, \abs{\cJ}_{H^{s-\frac{1}{2}}(\cS)})\, \times \\
			&\qquad\quad \times \qty(\abs{\pd_\beta\kaz}_{H^{\sigma+1}(\Gms)} + \abs{\pd_\beta\pd_t\kaz}_{H^{\sigma}(\Gms)} + \norm{\pd_\beta\omz}_{H^{\sigma+\frac{3}{2}}(\Oms)} + \norm{\pd_\beta\jz}_{H^{\sigma+\frac{3}{2}}(\Oms)}),
		\end{split}
	\end{equation}
	where $\cQ$ is a generic polynomial determined by $\Lambda_s$.
	
	Thus, \eqref{var it eqn ka} can be estimated as Proposition \ref{prop lin est theta} by treating all terms on the R.H.S. as source terms.
	
	On the other hand, taking variational derivatives of \eqref{it curl}, one can obtain that
	\begin{equation}\label{eqn it var om j}
		\left\{
		\begin{split}
			&\pd_t\Dbt(\omega^1-\jmath^1) + (\vz+\hz)\vdot\grad\Dbt(\omega^1-\jmath^1) + \Dbt\cB\qty[\pd(\vz+\hz), \pd(\vz-\hz)] \\
			&\quad = - \comm{\Dbt}{\pd_t + \grad_{(\vz+\hz)}}(\omega^1-\jmath^1) \qin \Omega^0_{t, \beta}, \\
			&\pd_t\Dbt(\omega^1+\jmath^1) + (\vz-\hz)\vdot\grad\Dbt(\omega^1+\jmath^1) + \Dbt\cB\qty[\pd(\vz-\hz), \pd(\vz+\hz)] \\
			&\quad = - \comm{\Dbt}{\pd_t + \grad_{(\vz-\hz)}}(\omega^1+\jmath^1) \qin \Omega^0_{t, \beta}, \\
			&\qty(\Dbt\omega^1)_{\restriction_{t=0}} = (\pd_\beta\wh{\omega_*})\circ\cX_{\Gamma^0_{t=0, \beta}}^{-1} \qc \qty(\Dbt\jmath^1)_{\restriction_{t=0}} = (\pd_\beta\wh{\jmath_*})\circ\cX_{\Gamma^0_{t=0, \beta}}^{-1},
		\end{split}
		\right.
	\end{equation}
	where $\omega^1 \coloneqq \omega_*^1\circ\cX_{\Gamma^0_{t, \beta}}$ and $\jmath^1 \coloneqq \jmath_*^1 \circ \cX_{\Gamma^0_{t, \beta}}$. The energy estimate for \eqref{eqn it var om j} can be derived as \eqref{curl evo est}.
	
	\subsubsection{Estimates with surface tension}
	When $\alpha > 0$, we assume that $s \ge \frac{9}{2}$. Consider the following energy functionals (here $l = 0, 1$):
	\begin{equation}
		\begin{split}
			\fkE^l \coloneqq\, &\sup_{t\in[0, T]} \left(\abs{\pd_\beta\kal}_{H^{s-\frac{5}{2}}(\Gms)}^2 + \abs{\pd_t\pd_\beta\kal}_{H^{s-4}(\Gms)}^2 + \norm{\pd_\beta\jmath_*^l}_{H^{s-\frac{5}{2}}(\Oms)}^2 + \right. \\
			&\qquad \qquad \left. + \norm{\pd_\beta\omega_*^l}_{H^{s-\frac{5}{2}}(\Oms)}^2 + \norm{\pd_t\pd_\beta\omega_*^l}_{H^{s-\frac{9}{2}}(\Oms)}^2 \right),
		\end{split}
	\end{equation}
	and
	\begin{equation}
			\wh{\fkE} \coloneqq  \abs{\pd_\beta\wh{\ka}}_{H^{s-\frac{5}{2}}(\Gms)}^2 + \abs{\pd_\beta(\wh{\pd_t\ka})}_{H^{s-4}(\Gms)}^2 + \norm{\pd_\beta\wh{\omega_*}}_{H^{s-\frac{5}{2}}(\Oms)}^2 + \norm{\pd_\beta\wh{\jmath_*}}_{H^{s-\frac{5}{2}}(\Oms)}^2.
	\end{equation}
	Then, uniform elliptic, product, div-curl, and commutator estimates, together with the linear energy estimates and Gr\"{o}nwall's inequality yield that (see also \cite[(6.27)-(6.29)]{Liu-Xin2023})
	\begin{equation}
		\begin{split}
			\fkE^1 \lesssim_{\Lambda_s}\, &\exp{T\cQ\qty(M_1, M_2, a^2M_3)} \cdot \qty{\wh{\fkE} + T\qty[M_2 + (1+a^2)\cQ(M_1)]\fkE^0} + \\
			&\quad + \cQ(M_1)\exp{T\cQ(M_1)}\cdot\qty[\wh{\fkE} + T\cQ(M_1)\fkE^0],
		\end{split}
	\end{equation}
	where $\cQ$ represents a generic polynomial determined by $\Lambda_s$ and $\alpha$.
	Thus, if $T$ is taken sufficiently small, one obtains
	\begin{equation}
		\fkE^1 \le \frac{1}{2}\fkE^0 + \cQ(M_1)\wh{\fkE}.
	\end{equation}
	In particular, if the initial data are not varying with $\beta$, the iteration map $\fkT$ defined by \eqref{def fkT} is a contraction with respect to the functional $\fkE$, which implies the uniqueness of the fixed point.
	
	\subsubsection{Estimates without surface tension}
	When $\alpha = 0$, it is required that the induced magnetic fields $\vh$ and $\vH$ satisfy \eqref{def lambda0}. Note that this can be achieved if $T \ll 1$ and \eqref{req ini h} holds initially. If $\alpha = 0$, the range of $s$ can be relaxed to $s \ge 4$. Consider the energy functionals defined by (here $l = 0, 1$):
	\begin{equation}
		\begin{split}
			\wt{\fkE}^l \coloneqq\, &\sup_{t\in[0, T]} \left(\abs{\pd_\beta\kal}_{H^{s-\frac{5}{2}}(\Gms)}^2 + \abs{\pd_t\pd_\beta\kal}_{H^{s-\frac{7}{2}}(\Gms)}^2 + \norm{\pd_\beta\jmath_*^l}_{H^{s-2}(\Oms)}^2 + \right. \\
			&\qquad \qquad \left. + \norm{\pd_\beta\omega_*^l}_{H^{s-2}(\Oms)}^2 + \norm{\pd_t\pd_\beta\omega_*^l}_{H^{s-4}(\Oms)}^2 \right),
		\end{split}
	\end{equation}
	and
	\begin{equation}
			\wt{\fkE}_{\text{ini}} \coloneqq \abs{\pd_\beta\wh{\ka}}_{H^{s-\frac{5}{2}}(\Gms)}^2 + \abs{\pd_\beta(\wh{\pd_t\ka})}_{H^{s-\frac{7}{2}}(\Gms)}^2 + \norm{\pd_\beta\wh{\omega_*}}_{H^{s-2}(\Oms)}^2 + \norm{\pd_\beta\wh{\jmath_*}}_{H^{s-2}(\Oms)}^2.
	\end{equation}
	Then, it is routine to check that 
	\begin{equation}
		\begin{split}
			\wt{\fkE}^1 \lesssim_{\Lambda_s}\, &\exp{T\wt{\cQ}\qty(M_1, M_2, a^2M_3)} \cdot \qty{\wt{\fkE}_{\text{ini}} + T\qty[M_2 + (1+a^2)\wt{\cQ}(M_1)]\wt{\fkE}^0} + \\
			&\quad + \wt{\cQ}(M_1)\exp{T\wt{\cQ}(M_1)}\cdot\qty[\wt{\fkE}_{\text{ini}} + T\wt{\cQ}(M_1)\wt{\fkE}^0],
		\end{split}
	\end{equation}
	where $\wt{\cQ}$ represents a generic polynomial determined by $\Lambda_s$ and $\lambda_0$. Thus, if $\alpha=0$ and \eqref{def lambda0} holds, the iteration map $\fkT$ also has a fixed point whenever $T$ is taken sufficiently small.
	
	\appendix
	\fancyhead[RO,LE]{\sc{\biolinum Appendices}}
	\section{Linearized Problems around Circular Background Solutions}\label{sec App A}
	The linearized equations \eqref{lin MHD} follow from the routine arguments, so we first focus on the derivation of linearized boundary conditions \eqref{lin BC}. Indeed, assume that the solutions $(\vv, \vh, \Omt)$ is parameterized by $\beta$, and the variational velocity with respect to $\beta$ is denoted by $\vps$. For simplicity, we denote by
	\begin{equation*}
		\Dbt \coloneqq \pd_\beta + \grad_\vps
	\end{equation*}
	the material derivative along the $\beta$-path, and
	\begin{equation}\label{decomp psi}
		\vps = \psi^\top\vtau + \psi^\perp \vn \qq{on} \Gmt.
	\end{equation}
	Then, applying $\Dbt$ to \eqref{BC}\textsubscript{2} yields that (recall $\vH\equiv\vb{0}$)
	\begin{equation*}
		\pd_\beta p + \grad_\vps p = \alpha\Dbt\varkappa = -\vn\vdot\grad_\vtau\grad_\vtau\vps - 2\varkappa(\grad_\vtau\vps\vdot\vtau).
	\end{equation*}
	Denote by $\wt{p} \coloneqq \pd_\beta p$, it follows from \eqref{decomp psi} and \eqref{deriv gmt} that
	\begin{equation*}
		\wt{p} = -\alpha\grad_\vtau\grad_\vtau\psi^\perp - \qty[\alpha\varkappa^2 + \grad_\vn p]\psi^\perp.
	\end{equation*}
	Similarly, it can be derived from \eqref{BC}\textsubscript{3} that
	\begin{equation*}
		\Dbt\vh \vdot\vn + \vh \vdot\Dbt\vn = 0,
	\end{equation*}
	which, together with \eqref{dt vn} and the notation $\wt{\vh} \coloneqq \pd_\beta\vh$, imply that
	\begin{equation*}
		\wt{\vh} \vdot\vn = \grad_\vtau\qty[\psi^\perp\vh\vdot\vtau].
	\end{equation*}
	Next, one can observe the condition \eqref{BC}\textsubscript{1} means that the fluid particles lying on the free interface will never leave. Taking a $C^2$ defining function $\Psi$ of $\Gamma_{t, \beta}$, for which
	\begin{equation*}
		\Gamma_{t, \beta} = \qty{\Psi(t, \beta, x) = 0} \qand \det(\grad_x \Psi) > 0 \text{ in a neighborhood of } \Gamma_{t, \beta}.
	\end{equation*}
	Then, it follows that
	\begin{equation*}
		\Dt\Psi \equiv 0 \qand \Dbt\Psi \equiv 0 \qq{on} \Gamma_{t, \beta}.
	\end{equation*}
	Through extending $\vps$ into a neighborhood of $\Gamma_{t, \beta}$ in an appropriate sense, it is legitimate to calculate that
	\begin{equation*}
		\comm{\Dt}{\Dbt} \coloneqq \Dt\Dbt - \Dbt\Dt = \grad_{\Dt\vps-\Dbt\vv}.
	\end{equation*}
	On the other hand, it is obvious that
	\begin{equation*}
		\comm{\Dt}{\Dbt}\Psi = \grad_{\Dt\vps-\Dbt\vv} \Psi = 0 \qq{on} \Gamma_{t, \beta},
	\end{equation*}
	which implies
	\begin{equation*}
		\qty(\Dt\vps-\Dbt\vv) \vdot \vn = 0 \qq{on} \Gamma_{t, \beta}.
	\end{equation*}
	Particularly, there holds (here $\wt{\vv} \coloneqq \pd_\beta\vv$)
	\begin{equation*}
		\Dt\psi^\perp = \wt{\vv}\vdot\vn + \psi^\perp (\grad_\vn \vv \vdot \vn) \qq{on} \Gamma_{t}.
	\end{equation*}
	In summary, one obtains the boundary conditions \eqref{lin BC}, which is actually independent of the profile of background solutions.
	
	Now, we consider the linearized problems around the circular stationary solutions \eqref{cir back sol}. For the sake of generality, we assume now that the parameters $\fkV$ and $\fkH$ are functions of $r$. Namely,
	\begin{equation*}
		\vv = r\fkV(r) \etheta \qand \vh = r\fkH(r) \etheta.
	\end{equation*}
	Here we remark that $\fkV(r)$ represents the angular velocity of the rotating flows, and the vorticity and current are respectively given by:
	\begin{equation*}
		\omega \coloneqq \Curl\vv = 2\fkV + r\pd_r\fkV \qand \jmath \coloneqq \Curl\vh = 2\fkH + r\pd_r\fkH.
	\end{equation*}
	The linearized equations \eqref{lin MHD} can be expressed in the polar coordinates as:
	\begin{equation}\label{lin wave eqn}
		\begin{cases*}
			\pd_t\wt{v}^r + \fkV\pd_\theta\wt{v}^r-2\fkV\wt{v}^\theta+\pd_r\wt{p}=-2\fkH\wt{h}^\theta + \fkH\pd_\theta\wt{h}^r, \\
			\pd_t\wt{v}^\theta+\fkV\wt{v}^r+\fkV\pd_\theta\wt{v}^\theta+\wt{v}^r\pd_r(r\fkV)+r^{-1}\pd_\theta\wt{p} = \wt{h}^r\pd_r(r\fkH)+\fkH\wt{h}^r+\fkH\pd_\theta\wt{h}^\theta, \\
			\pd_t\wt{h}^r+\fkV\pd_\theta\wt{h}^r=\fkH\pd_\theta\wt{v}^r, \\
			\pd_t\wt{h}^\theta+\fkV\pd_\theta\wt{h}^\theta+r\wt{v}^r\pd_r\fkH = r\wt{h}^r\pd_r\fkV + \fkH\pd_\theta\wt{v}^\theta, \\
			\pd_r(r\wt{v}^r)+ \pd_\theta\wt{v}^\theta = \pd_r(r\wt{h}^r)+\pd_\theta\wt{h}^\theta = 0.
		\end{cases*}
	\end{equation}
	Plugging the perturbation profile \eqref{ptb profile} into \eqref{lin wave eqn} yields that
	\begin{equation*}
		\begin{cases*}
			(c-\fkV)\pd_r\qty[r\pd_r(r\wh{v}^r)]+\wh{v}^r r\pd_r(2\fkV+r\pd_r\fkV) + k^2(\fkV-c)\wh{v}^r \\
			\quad = - \fkH\pd_r\qty[r\pd_r(r\wh{h}^r)] + \wh{h}^r r\pd_r(2\fkH+r\pd_r\fkH) + k^2 \fkH\wh{h}^r, \\
			(\fkV-c)\wh{h}^r = \fkH\wh{v}^r,
		\end{cases*}
	\end{equation*}
	which implies \eqref{lin wave pur}, since $\fkV$ and $\fkH$ are both constant there. Indeed, the above equations can be reduced to \eqref{lin wave pur} whenever
	\begin{equation*}
		\pd_r\omega = \pd_r\jmath = 0,
	\end{equation*}
	i.e., the flow has constant vorticity and the magnetic field corresponds to a constant current.
	
	Concerning the dispersive relations, it follows from \eqref{cir back sol} that
	\begin{equation*}
		\grad_\vn p = r\fkV^2 - r\fkH^2 \qq{on} \{r = 1\}.
	\end{equation*}
	Moreover, it follows from \eqref{ptb profile}, \eqref{lin BC}\textsubscript{2}, and \eqref{renormalize s} that
	\begin{equation*}
		ik\wh{p}= \alpha\qty(k^2-1)+\fkH^2-\fkV^2 \qq{when} r=1.
	\end{equation*}
	Therefore, one can derive from \eqref{lin wave eqn}\textsubscript{2, 5} and \eqref{ptb profile} that
	\begin{equation*}
		(c-\fkV)z'(0) + 2\fkV z(0) + \alpha\qty(k^2-1)+\fkH^2-\fkV^2 = 2\fkH^2(\fkV-c)^{-1}z(0) - \fkH^2(\fkV-c)^{-1}z'(0),
	\end{equation*}
	where $z$ is defined through \eqref{def z}, and $z' \equiv \dv*{z}{s}$. Particularly, the solution formula \eqref{sol z} implies that
	\begin{equation*}
		k(c-\fkV)(\fkV-c) + 2\fkV(\fkV-c) + \alpha\qty(k^2-1) + \fkH^2-\fkV^2 = 2\fkH^2 - k\fkH^2,
	\end{equation*}
	which can be reduced to \eqref{disp}.

	\section*{Acknowledgments}
	Liu's research is supported by the UM Postdoctoral Fellow scheme under the UM Talent Programme at the University of Macau. Liu extends gratitude to Prof. Changfeng Gui at the University of Macau for his encouragement and support.
	
	Luo's research is supported by a grant from the Research Grants Council of the Hong Kong Special Administrative Region, China (Project No. 11310023).
	
	
	\subsection*{Data availability}
	This manuscript has no associated data.
	
	\subsection*{Conflict of interest}
	The authors have no conflict of interest to disclose.
	
	\bibliographystyle{amsplain0}
	\fancyhead[RO,LE]{\sc{\biolinum References}}
	{\small	\bibliography{ref}}
	\end{document}